\documentclass[letterpaper,12pt,reqno]{amsart}
\RequirePackage[utf8]{inputenc}
\usepackage[portrait,margin=3cm]{geometry}
\usepackage{mathrsfs,soul}
\usepackage{hyperref}
\usepackage[foot]{amsaddr}
\usepackage{amsmath,amssymb,amsthm,amsfonts,amsbsy,latexsym,dsfont,mathtools,nccmath,scalerel}

\makeatletter 
\@mparswitchfalse%
\makeatother
\normalmarginpar 
\usepackage[textsize=small]{todonotes}       
\usepackage{graphicx}
\usepackage[numeric,initials,nobysame]{amsrefs}
\usepackage{upref,setspace}
\usepackage{xcolor,colortbl}

\usepackage{fourier}

\usepackage{xcolor}

\usepackage{amsbsy,latexsym,dsfont}
\usepackage{mathrsfs,soul}

\usepackage{enumerate}

\newenvironment{enumeratea}{\begin{enumerate}[\upshape (a)]}{\end{enumerate}}

\usepackage{paralist}

\newenvironment{inparaenuma}{\begin{inparaenum}[\upshape (a) ]}{\end{inparaenum}}

\definecolor{refkey}{gray}{.75}
\definecolor{labelkey}{gray}{.75}

\newtheorem{thm}{Theorem}[section]
\newtheorem{lem}[thm]{Lemma}
\newtheorem{cor}[thm]{Corollary}
\newtheorem{prop}[thm]{Proposition}

\newtheorem{ass}[thm]{Condition}

\newtheorem{conj}[thm]{Conjecture}

\theoremstyle{remark}

\theoremstyle{definition}
\newtheorem{rem}{Remark}

\newcommand{\ind}{\mathds{1}}
\newcommand{\eps}{\varepsilon}

\newcommand{\set}[1]{\left\{#1\right\}}

\newcommand{\ie}{\emph{i.e.,}}

\newcommand{\equald}{\stackrel{\mathrm{d}}{=}}

\newcommand{\weakc}{\stackrel{\mathrm{d}}{\longrightarrow}}
\newcommand{\convas}{\stackrel{\mathrm{a.s.}}{\longrightarrow}}

\def\qed{ \hfill $\blacksquare$}

\newcommand{\cA}{\mathcal{A}}\newcommand{\cB}{\mathcal{B}}\newcommand{\cC}{\mathcal{C}}
\newcommand{\cD}{\mathcal{D}}\newcommand{\cE}{\mathcal{E}}\newcommand{\cF}{\mathcal{F}}
\newcommand{\cG}{\mathcal{G}}

\newcommand{\cM}{\mathcal{M}}
\newcommand{\cP}{\mathcal{P}}\newcommand{\cR}{\mathcal{R}}
\newcommand{\cS}{\mathcal{S}}\newcommand{\cT}{\mathcal{T}}
\newcommand{\cV}{\mathcal{V}}
\newcommand{\cY}{\mathcal{Y}}
\newcommand{\vone}{\mathbf{1}}

\newcommand{\vJ}{\mathbf{J}}
\newcommand{\vM}{\mathbf{M}}

\newcommand{\vX}{\mathbf{X}}
\newcommand{\vb}{\mathbf{b}}
\newcommand{\ve}{\mathbf{e}}

\newcommand{\vp}{\mathbf{p}}
\newcommand{\vt}{\mathbf{t}}\newcommand{\vu}{\mathbf{u}}
\newcommand{\vv}{\mathbf{v}}\newcommand{\vw}{\mathbf{w}}\newcommand{\vx}{\mathbf{x}}

\newcommand{\mvzero}{\boldsymbol{0}}
\newcommand{\mvone}{\boldsymbol{1}}

\newcommand{\mva}{\boldsymbol{a}}

\newcommand{\mvx}{\boldsymbol{x}}

\newcommand{\mvbeta}{\boldsymbol{\beta}}

\newcommand{\mvmu}{\boldsymbol{\mu}}

\newcommand{\mvxi}{\boldsymbol{\xi}}\newcommand{\mvXi}{\boldsymbol{\Xi}}

\newcommand{\fF}{\mathfrak{F}}
\newcommand{\fG}{\mathfrak{G}}\newcommand{\fH}{\mathfrak{H}}

\newcommand{\fT}{\mathfrak{T}}
\newcommand{\fX}{\mathfrak{X}}
\newcommand{\fZ}{\mathfrak{Z}}
\newcommand{\fb}{\mathfrak{b}}

\newcommand{\bE}{\mathbb{E}}
\newcommand{\bG}{\mathbb{G}}

\newcommand{\bN}{\mathbb{N}}
\newcommand{\bP}{\mathbb{P}}\newcommand{\bR}{\mathbb{R}}
\newcommand{\bT}{\mathbb{T}}

\newcommand{\bZ}{\mathbb{Z}}

\newcommand{\sP}{\mathscr{P}}

\newcommand{\sS}{\mathscr{S}}

\DeclareMathOperator{\E}{\mathbb{E}}
\DeclareMathOperator{\pr}{\mathbb{P}}

\DeclareMathOperator{\var}{Var}

\DeclareMathOperator{\dis}{dis}

\DeclareMathOperator{\GH}{GH}
\DeclareMathOperator{\GHP}{GHP}

\DeclareMathOperator{\ord}{ord}

\DeclareMathOperator{\con}{con}

\DeclareMathOperator{\ERRG}{ERRG}

\DeclareMathOperator{\dsc}{dsc}
\DeclareMathOperator{\scl}{scl}

\DeclareMathOperator{\mass}{mass}

\DeclareMathOperator{\spls}{spls}

\DeclareMathOperator{\crit}{Crit}

\DeclareMathOperator{\pois}{Pois}

\DeclareMathOperator{\osc}{Osc}

\newcommand{\sss}{\scriptscriptstyle}
\newcommand{\erdos}{Erd\H{o}s-R\'enyi }
\newcommand{\ldown}{l^2_{\downarrow}}

\newcommand{\diam}{\mathrm{diam}}

\newcommand{\vCrit}{\mathbf{Crit}}

\usepackage[textsize=small]{todonotes}       

\definecolor{jrnl}{rgb}{0.0, 0.5, 0.0}
\definecolor{jrnl1}{rgb}{0.0, 0.72, 0.6}
\definecolor{aqua}{rgb}{0.0, 1.0, 1.0}
\definecolor{webbrown}{rgb}{.6,0,0}
\definecolor{pinegreen}{rgb}{0.0, 0.47, 0.44}
\definecolor{ultramarineblue}{rgb}{0.25, 0.4, 0.96}
\definecolor{jrnl}{rgb}{0.0, 0.5, 0.0}
\definecolor{lincolngreen}{rgb}{0.11, 0.35, 0.02}
\definecolor{green(html/cssgreen)}{rgb}{0.0, 0.5, 0.0}
\definecolor{airforceblue}{rgb}{0.36, 0.54, 0.66}
\definecolor{azure}{rgb}{0.0, 0.5, 1.0}
\definecolor{bleudefrance}{rgb}{0.19, 0.55, 0.91}
\definecolor{cobalt}{rgb}{0.0, 0.28, 0.67}

\hypersetup{colorlinks=true, linktocpage=true, pdfstartpage=1, pdfstartview=FitV,
	breaklinks=true, pdfpagemode=UseNone, pageanchor=true, pdfpagemode=UseOutlines,
	plainpages=false, bookmarksnumbered, bookmarksopen=true, bookmarksopenlevel=1,
	hypertexnames=true, pdfhighlight=/O,
	urlcolor=webbrown, linkcolor=cobalt, citecolor=jrnl
}

\newcommand{\bmu}{\boldsymbol{\mu}}

\newcommand{\chsen}[1]{\textcolor{black}{{#1}}}

\newcommand{\ch}[1]{{#1}}

\newcommand{\chh}[1]{{#1}}

\newcommand{\graphon}{\mathcal{G}^{\sss \text{gphn}}}
\newcommand{\tildegraphon}{\widetilde{\mathcal{G}}^{\sss \text{gphn}}}
\newcommand{\graphonn}{G^{\sss \text{gphn}}}
\newcommand{\rgiv}{\mathcal{G}^{\sss \text{rgiv}}}

\def\Xint#1{\mathchoice
	{\XXint\displaystyle\textstyle{#1}}%
	{\XXint\textstyle\scriptstyle{#1}}%
	{\XXint\scriptstyle\scriptscriptstyle{#1}}%
	{\XXint\scriptscriptstyle\scriptscriptstyle{#1}}%
	\!\int}
\def\XXint#1#2#3{{\setbox0=\hbox{$#1{#2#3}{\int}$ }
		\vcenter{\hbox{$#2#3$ }}\kern-.6\wd0}}

\def\dashint{\Xint-}

\makeatletter
\providecommand{\leftsquigarrow}{%
	\mathrel{\mathpalette\reflect@squig\relax}%
}
\newcommand{\reflect@squig}[2]{%
	\reflectbox{$\m@th#1\rightsquigarrow$}%
}
\makeatother

\newcommand{\connects}{\leftrightsquigarrow}
\newcommand{\connectsell}{\stackrel{\ell}{\leftrightsquigarrow}}

\begin{document}

	\title[Universality: $L^3$ graphons and critical percolation]{Scaling limits and universality: Critical percolation on weighted graphs converging to an $L^3$ graphon}
	
	\date{}

\author[Baslingker]{Jnaneshwar Baslingker$^1$}
\address{$^1$Department of Mathematics, Indian Institute of Science}	
	\author[Bhamidi]{Shankar Bhamidi$^2$}
	\address{$^2$Department of Statistics and OR, University of North Carolina, Chapel Hill}
	\author[Broutin]{Nicolas Broutin$^3$}
	\address{$^3$LSPM, Sorbonne Universit\'{e}, France, and Institut Universitaire de France (IUF)}
	\author[Sen]{Sanchayan Sen$^4$}
	\address{$^4$Department of Mathematics, Indian Institute of Science}
	\author[Wang]{Xuan Wang$^5$}
	\address{$^5$Databricks}
	\maketitle
	\begin{abstract}
We develop a general universality technique for establishing metric scaling limits of critical random discrete structures exhibiting mean-field behavior that requires four ingredients: (i) from the barely subcritical regime to the critical window, components merge approximately like the multiplicative coalescent, (ii) asymptotics of the susceptibility functions are the same as that of the \erdos random graph, (iii) asymptotic negligibility of the maximal component size and the diameter in the barely subcritical regime, and (iv) macroscopic averaging of distances between vertices in the barely subcritical regime.

As an application of the general universality theorem, we establish, under some regularity conditions, the critical percolation scaling limit of graphs that converge, in a suitable topology, to an $L^3$ graphon. 
In particular, we define a notion of the critical window in this setting.
The $L^3$ assumption ensures that the model is in the \erdos universality class and that the scaling limit is Brownian.
Our results do not assume any specific functional form for the graphon. 
As a consequence of our results on graphons, we obtain the metric scaling limit for Aldous-Pittel's RGIV model \cite{aldous2000random} inside the critical window.

Our universality principle has applications in a number of other problems including in the study of noise sensitivity of critical random graphs \cite{lubetzky-peled-noise}.
In \cite{bhamidi-broutin-sen-wang-bsr}, we use our universality theorem to establish the metric scaling limit of critical bounded size rules. 
Our method should yield the critical metric scaling limit of Ruci\'nski and Wormald's random graph process with degree restrictions \cite{rucinski-wormald} provided an additional technical condition about the barely subcritical behavior of this model can be proved.
\end{abstract}

\tableofcontents

\section{Introduction}\label{sec:intro}
Over the last two decades, motivated both by questions in fields such as combinatorics and statistical physics as well as the explosion in data on real networks, an array of random graph models have been proposed. 
Many of these models undergo a phase transition in the size of the largest component.
The canonical example of such phenomenon is the \erdos random graph $\ERRG(n,t/n)$ where the phase transition occurs at the critical time $t_c(\ERRG)=1$.
One fundamental question in the study of random discrete structures is understanding the nature of this phase transition, and see how they compare to the phase transition in the \erdos random graph.

In the early 2000s, based on simulations, it was predicted in the statistical physics literature 
\cites{BraBulCohHavSta03,wu2006transport,braunstein2007optimal,chen2006universal}
that asymptotics in the critical regime are ``universal.''
More precisely, under moment assumptions on the average degree, a wide array of models are conjectured to exhibit the same behavior in the critical regime as the \erdos random graph in the sense that for any fixed $i\geq 1$, the size of the $i$-th maximal component scales like $n^{2/3}$, and the typical distance in the component scales like $n^{1/3}$. 
Denoting the $i$-th maximal component of the \erdos random graph $\ERRG(n, n^{-1}+\lambda n^{-4/3})$ by $\cC_i(\lambda)$, it was proved in \cite{aldous1997brownian} that the sequence of rescaled component sizes
$
\big(n^{-2/3}|\cC_i(\lambda)|\, ; \, i\geq 1\big)
$
conervges in distribution to a limit that will be defined in Section~\ref{sec:smc-def}.
Turning to the intrinsic geometry of the components, it was shown in \cite{addario2012continuum} that rescaling edge lengths in the $i$-th maximal component $\cC_i(\lambda)$ of the \erdos random graph $\ERRG(n , 1/n+\lambda/n^{4/3})$ and denoting the resulting metric space by $n^{-1/3}\cC_i(\lambda)$, one has
\begin{equation}\label{eqn:lim-add-br-go}
	\big(n^{-1/3}\cC_i(\lambda);\, i\geq 1\big) 
	\weakc 
	\vCrit(\lambda) := \big(\crit_i(\lambda);\, i\geq 1 \big),
\end{equation}
for a sequence of limiting random fractals that are described in more detail in Section~\ref{sec:prelim}.
This scaling limit is expected to be the universal scaling limit for a large class of critical random graph models exhibiting mean-field behavior.
Further, existing literature suggests that the components of the high-dimensional discrete torus \cite{hofstad-sapozhnikov,heydenreich-hofstad-1,heydenreich-hofstad-2} and the hypercube \cite{hofstad-nachmias} under critical percolation also share the \erdos scaling limit.

There are two standard techniques for establishing scaling limits as above:
(a) {\it The exploration technique:}
Here, the random graph under consideration is explored, often in a breadth-first or a depth-first manner, and alongside one keeps track of certain associated processes, e.g., the exploration-walk process, the height process, the contour process, or the local time field associated with the height/contour process.
One then proves distributional limits for these associated processes and the metric scaling limit is derived from these limit theorems.
(b) {\it The spanning tree technique:}
Here, one first establishes the scaling limit of the rescaled component sizes of the model being considered.
Then conditional on the component sizes, one has to find a way to represent each connected component as a well-behaved spanning tree with certain additional random edges added to it.
Then the metric scaling limit of each connected component can be derived using the scaling limit of the spanning tree and also incorporating the location of the additional random edges.
When applicable, these techniques can be powerful tools in obtaining precise results about the geometry of the model being considered.
However, the model being studied has to have some nice properties in order to be amenable to analysis via these methods.
For the first method to be applicable, one has to be able to set up a suitable exploration process for which it is possible to obtain the distributional limit in the Skorohod topology along with the associated height/contour process.
Similarly, for the second technique, each connected component needs to have a spanning tree whose metric scaling limit can be established, which is usually done by connecting its distribution to a well-studied model of random trees.
These two methods have been successfully applied to the following two models of random graphs:
the inhomogeneous random graph with a rank-one kernel (including the \erdos random graph), and the configuration model and (simple) random graphs with a given degree sequence.
See 
\cite{riordan2012phase,joseph2010component,goldschmidt-conchon-stable-graph,goldschmidt-haas-senizergues,bhamidi-sen-vacant-set,dhara-hofstad-leeuwaarden-sen-cm-finite-third-moment,sd-rvdh-jvl-ss-cm-component-scaling-infinite-third-moment,SB-SD-vdH-SS-cm-metric-infinite-third-moment,SB-SD-vdH-SS-global-lower-mass,dhara-hofstad-leeuwaarden,miermont-sen,SBSSXW14,bhamidi-hofstad-sen-multiplicative,broutin-duquesne-wang,addario2012continuum,aldous1997brownian,addario2010critical}
for various results on the critical scaling limits for these two models proved under different sets of assumptions.
In general, it is not clear how the above two methods can be applied to the many other standard models of random discrete structures studied in the literature.

\subsection{Main contributions of this paper}
We prove a general universality principle for systematically approaching the problem of establishing metric scaling limits of critical random discrete structures exhibiting mean-field behavior.
Applying this technique requires the verification of four properties: 
(i) From the barely subcritical regime to the critical window, components merge approximately like the multiplicative coalescent.
(ii) The second and the third susceptibility functions at certain suitable points in the barely subcritical regime exhibit the same asymptotics as those of the \erdos random graph.
(iii)  The size of the maximal component and the diameter in the barely subcritical regime is asymptotically negligibile in an appropriate sense.  
(iv) The average distances between vertices in the barely subcritical regime is well-behaved.
Often, these four criteria can be verified for models where it is not clear if the exploration process technique or the spanning tree technique can be used directly, thus allowing for broader applicability.

We apply our general theorem to random graphs obtained under critical percolation on edge-weighted graphs that converge, in a suitable topology, to an $L^3$ graphon.
The theory of graphons is now a standard tool for sampling and studying convergence properties of dense graph sequences, and has found applications in many branches of applied science, especially in statistics and machine learning.
The precise definition of the model we are considering is given in Section~\ref{sec:def-graphon-percolation}.
The phase transition in this model follows from the results in \cite{BBSJOR07,bollobas-borgs-chayes-riordan,bollobas-janson-riordan-cut-norm}.
Note that in all these works, only the transition from the purely subcritical regime to the purely supercritical regime was studied--the near critical behavior was not considered.

As a consequence of our universality principle, we will show that under some regularity conditions, the metric scaling limit in this model is the same as that of the \erdos random graph (given in \eqref{eqn:lim-add-br-go}) modulo a constant.
To our knowledge, this gives the first result where the existence of a critical window is established and the metric scaling limit is obtained inside the critical window without any assumptions on the specific functional form of the graphon.
(As mentioned above, the metric scaling limit for this model was only known in the special case of the rank-one setting where the kernel has a multiplicative form: $W(x, y) = \phi(x)\phi(y)$, $x, y\in [0, 1]$.)
Here, the $L^3$ condition ensures that the scaling limit is of a Brownian nature, i.e., the model belongs to the \erdos universality class.
There are other regimes where we expect scaling limits of a different type to arise;
we will revisit this point in Remark~\ref{rem:4}.
As a further application of our result on graphons, we obtain the metric scaling limit of Aldous-Pittel's RGIV model \cite{aldous2000random} inside the critical window.

The main universality principle of this paper can be applied to several other problems.
In \cite{bhamidi-broutin-sen-wang-bsr}, we apply our universality theorem to 
dynamic random graph models including bounded size rules  \cite{spencer2007birth} and establish the metric scaling limit inside the critical window.
In \cite{lubetzky-peled-noise}, the authors use this universality result as one of the ingredients to study noise sensitivity of critical random graphs.
We expect that our general method can be used  to obtain the critical metric scaling limit of Ruci\'nski and Wormald's random graph process with degree consraints \cite{rucinski-wormald}. 
At present, we are missing a technical condition about the barely subcritical behavior of this model.
We will briefly discuss this in Section~\ref{sec:disc}.
Further, as shown in \cite{addario2013scaling} 
(see also \cite{addarioberry-sen,bhamidi-sen-heavy-tailed-mst}), 
obtaining the metric structure of maximal components at criticality is a key tool in analyzing more complicated objects such as the minimal spanning tree (MST) on the giant component in the supercritical regime. 
Establishing the critical scaling limits such as those obtained in this paper would be one of the steps in establishing universality for such objects as well.

\medskip

{\textit{\textbf {Organization of the paper:}}
In Section~\ref{sec:models} we define the specific random graph models considered in this paper.
In Section~\ref{sec:prelim}, we define mathematical constructs needed to state the results. 
Then Section~\ref{sec:res} describes the general universality result.  
Section~\ref{sec:res-models} contains results on applications of the universality result to specific models.
The proofs of these results are given Sections~\ref{sec:proofs-universality} and \ref{sec:graphon-proofs}.
Finally, in Section~\ref{sec:disc}, we state a result on the scaling limit of the critical stochastic block model and discuss possible extensions of our results.

\subsection{Models}
\label{sec:models}
In this section, we define the two specific models to which our general universality result (Theorem~\ref{thm:gen-2} stated below) is applied in this paper.

\subsubsection{Graphons and critical percolation on edge-weighted graphs}\label{sec:def-graphon-percolation}
For $p\geq 1$, a Borel measurable function $W:[0, 1]^2 \to [0,\infty)$ is called an $L^p$ graphon if
$W$ is symmetric (i.e., $W(x, y)=W(y, x)$ for $x, y\in [0, 1]$), and 
$\int_0^1\int_0^1 W(x, y)^p dx dy<\infty$.
Clearly, an $L^p$ graphon is also an $L^q$ graphon if $p>q$.

Now, to any $K\in L^2([0, 1]^2)$ we can associate a bounded linear operator 
$T_K: L^2[0, 1]\to L^2[0, 1] $ via
\begin{align}\label{eqn:3}
	T_K f(\cdot):=\int_0^1 K(\cdot, y)f(y) dy\, ,\ \ f\in L^2[0, 1]\, .
\end{align}
Further, for $p, q\in [1, \infty]$ and a linear operator $T^{\bullet}: L^p[0, 1]\to L^q[0, 1]$, we write
\begin{align}\label{eqn:54}
	\|T^{\bullet}\|_{p,q}
	=
	\sup\big\{ \|T^{\bullet} f\|_q\, :\, \| f\|_p =1 \big\}
\end{align}
for the operator norm of $T^{\bullet}$.
We call an $L^2$ graphon $W$ critical if 
$\|T_W\|_{2,2}=1$.
The reason for this terminology comes from the criterion for an inhomogeneous random graph to be critical \cite[Theorem~6.1]{BBSJOR07}.

Consider a sequence of graphs with vertex sets $[n] = \{1, 2, \ldots, n \}$ and edge weights 
$\mvbeta_n =\big(\beta_{ij}^{\sss(n)}\, ;\, i, j\in [n]\big)$,
where $\beta_{ii}^{\sss(n)}=0$ for $i\in [n]$ and 
$\beta_{ij}^{\sss(n)}=\beta_{ji}^{\sss(n)}\geq 0$ for $1\leq i\neq j\leq n$.
Using the weight sequence $\mvbeta_n$, we will formulate a notion (see Conditions~\ref{ass:graphon-1}, \ref{ass:graphon-2}, and \ref{ass:graphon-3} below)
for these graphs to converge to a critical $L^3$ graphon.
We write
\begin{align}\label{eqn:222}
p_{ij}^{\sss (n)} 
= 
1\wedge~(n^{-1}\beta_{ij}^{\sss(n)})
\ \ \text{ and }\ \
\widetilde p_{ij}^{\sss (n)} 
= 
1 - \exp\big(-\beta_{ij}^{\sss(n)}/n\big)
\, , \ \ i\leq i\neq j\leq n\, .
\end{align}
Let $\graphon_n = \graphon_n(\mvbeta_n)$ (resp. $\tildegraphon_n = \tildegraphon_n(\mvbeta_n)$) denote the random graph with vertex set $[n]$, where the edges $\{i, j\}$ appear independently with respective probabilities
$p_{ij}^{\sss (n)} $ (resp. $\widetilde p_{ij}^{\sss (n)} $), $1\leq i<j\leq n$.
We will use the notation $\graphonn_n$ to refer to either of the graphs $\graphon_n$ and $\tildegraphon_n$, as all our results will hold for both these graphs.
As a consequence of our general universality result (Theorem~\ref{thm:gen-2}), we will obtain the metric scaling limit of this model (Theorem~\ref{thm:scaling-limit-graphon}).

\subsubsection{Aldous and Pittel's RGIV process \cite{aldous2000random}:}
\label{sec:def-rgiv}
The random graph with immigrating vertices (RGIV) process is a dynamic random graph model that evolves as follows:
At time $0$, the RGIV process is the empty graph.
Then vertices immigrate at rate $n$, and for each pair of existing vertices, an edge connecting them appears at rate $1/n$.

A deterministic version of this process was studied in the physics literature \cite{spouge}, where it was observed that a phase transition occurs in this process at time $\pi/2$; see also \cite{capobianco-frank}.
Aldous and Pittel \cite{aldous2000random} studied the behavior of this process inside the critical window, i.e., at time 
\begin{align}\label{eqn:1a}
t_{n,\lambda}=\pi/2 +\lambda n^{-1/3}\, ,
\end{align}
and they established the scaling limit of the component sizes.
For $\lambda\in\bR$, write $\rgiv_n(\lambda)$ for the random graph process at time $t_{n,\lambda}$ as given by \eqref{eqn:1a}.
As an application of Theorem~\ref{thm:scaling-limit-graphon}, we will obtain the metric scaling limit of $\rgiv_n(\lambda)$ (Theorem~\ref{thm:rgiv-scaling-limit}).

\section{Preliminaries}
\label{sec:prelim}


\subsection{GHP convergence of compact metric measure  spaces}\label{sec:gh-mc}

We primarily follow \cite{EJP2116,addario2013scaling,burago2001course}. 
For a compact metric space $(X, d)$ and $A_1, A_2\subseteq X$, we define the Hausdorff distance between $A_1$ and $A_2$ to be
\[
d_{\text{H}}(A_1, A_2):=\inf\big\{\eps>0\ :\ A_1\subseteq A_2^{\eps}\ \text{ and }\ A_2\subseteq A_1^{\eps}\big\},
\]
where $A_1^{\eps}:=\bigcup_{x\in A_1}\{y\in X\, :\, d(y, x)\leq \eps\}$.
For two metric spaces $X_1 = (X_1,d_1)$ and $X_2 = (X_2, d_2)$ and $\cR\subseteq X_1 \times X_2$, define the distortion of $\cR$ with respect to $X_1$ and $X_2$ as
\begin{equation}\label{eqn:def-distortion}
\dis(\cR) =\dis(\cR; X_1, X_2)
:= \sup \big\{|d_1(x_1,y_1) - d_2(x_2, y_2)|: (x_1,x_2) , (y_1,y_2) \in \cR\big\}.
\end{equation}
A correspondence $\cR$ between $X_1$ and $X_2$ is a measurable subset of $X_1 \times X_2$ such that for every $x_1 \in X_1$ there exists at least one $x_2 \in X_2$ such that $(x_1,x_2) \in \cR$ and vice-versa. Let $\text{Corr}(X_1,X_2)$ be the set of all such correspondences.  The Gromov-Hausdorff distance between $(X_1,d_1)$ and $(X_2, d_2)$ is defined as
\begin{equation}
\label{eqn:dgh}
d_{\GH}(X_1, X_2) = \frac{1}{2}\inf \big\{\dis(\cR): \cR \in \text{Corr}(X_1, X_2)\big\}.
\end{equation}

The Gromov-Hausdorff-Prokhorov distance extends the above to keep track of measures on the corresponding spaces.
A compact metric measure space is a compact metric space equipped with a finite measure on its Borel sigma algebra.
Let $X_1 = (X_1, d_1, \mu_1)$ and $X_2=(X_2,d_2, \mu_2)$ be compact metric measure spaces.
Suppose $\pi$ is a measure on the product space $X_1\times X_2$ with marginals $\pi_1, \pi_2$. 
The discrepancy of $\pi$ with respect to $X_1$ and $X_2$ is defined as
\begin{equation}
\label{eqn:def-discrepancy}
\dsc(\pi)=\dsc(\pi; X_1, X_2):= ||\mu_1-\pi_1|| + ||\mu_2-\pi_2||,
\end{equation}
where $||\cdot||$ denotes the total variation of signed measures. 
Now define
\begin{equation}\label{eqn:dghp}
\bar d_{\GHP}(X_1, X_2):= \inf\Big\{ \max\big(\dis(\cR)/2, ~\dsc(\pi),~\pi(\cR^c)\big) \Big\},
\end{equation}
where the infimum is taken over all correspondences $\cR$ and measures $\pi$ on $X_1 \times X_2$. The function $\bar d_{\GHP}$ is a pseudometric and defines an equivalence relation $X_1\sim X_2$ if and only if $\bar d_{\GHP}(X_1,X_2) = 0$ on. Let $\sS $ be the space of all equivalence classes of compact metric measure spaces and let $d_{\GHP}$ be the induced metric. Then by \cite{EJP2116}, $(\sS, d_{\GHP})$ is a complete separable metric space. 
We will continue to use $X = (X, d, \mu)$ to denote both the metric space and {the} corresponding equivalence class. 
In this paper we will typically deal with a sequence of infinite sequences of compact metric measure spaces; we will view it as a sequence in $\sS^{\bN}$.
When talking about distributional convergence for a sequence of $\sS^{\bN}$-valued random objects, it will be understood that the underlying topology on $\sS^{\bN}$ is the product topology inherited from $d_{\GHP}$.

\medskip

\textbf{The scaling operator:} We will need to rescale both the metric as well as the associated measures on the components in the critical regime. For $a, b> 0$, let $\scl(a, b):  \sS \to \sS$ be the scaling operator
\begin{align*}
\scl(a, b)[(X , d , \mu)]:= (X, d',\mu'),
\end{align*}
where $d'(x,y) := a d(x,y)$ for all $x,y \in X$, and $\mu'(A) := b \mu(A)$ for all $A \subseteq X$. 
For simplicity, write $\scl(a, b) X := \scl(a, b)[(X , d , \mu)]$ and $a X := \scl(a, 1) X$.
\chsen{We will write $aX$ and $a\cdot X$ interchangeably.}

\subsection{Notation}
\label{sec:gr-constr}
We write $|A|$ or $\# A$ for the cardinality of a set $A$.
For any $n\in\bZ_{>0}$, we will write $[n]$ for the set $\{1, 2, \ldots, n\}$.
For a finite graph $G$, we write $V(G)$ and $E(G)$ for the set of vertices and the set of edges in $G$ respectively.
We write $|G|$ to mean $|V(G)|$, and
$\spls(G) := |E(G)| - |G|+1$ for the number of surplus edges in $G$.
We view a connected component $\cC$ of $G$ as a metric space using the graph distance $d_{G}$.
When the graph $G$ is clear from the context, we will supress the subscript and simply write $d$ for the graph distance. 
Often in our analysis, associated to a graph $G$ there will be a collection of vertex weights 
$\vw = \big\{w_v: v\in V(G)\big\}$.
Natural measures on $G$ are
\begin{inparaenuma}
	\item {\it Counting measure:} $\mu_{\text{ct}}(A): = |A|$, for $A \subseteq V(G)$; and
	\item {\it Weighted measure:} $\mu_\vw(A) :=  \sum_{v \in A} w_v$, for $A \subseteq V(G)$. 
\end{inparaenuma}
For $G$ finite and connected, we will often use $G$ for both the graph and its associated metric measure space.	
If the measure to be considered is not the counting measure then that will be explicitly mentioned.
If nothing is specified, the measure will be understood to be the counting measure.
For two graphs $G_1$ and $G_2$, we write $G_1\subseteq G_2$ to mean that $G_1$ is a subgraph of $G_2$.
For a metric space $(\fX, d)$, write $\diam(\fX):=\sup_{u, v\in\fX}d(u, v)$.
For a finite graph $G$, we let 
$\diam(G):=\max\big\{\diam(\cC)\, :\, \cC\text{ connected component of }G\big\}$. 

We will say that a sequence $\cE_n,\, n\geq 1$, of events occurs with high probability (whp) if $\pr(\cE_n)\to 1$.
We will freely omit ceilings and floors when there is no confusion in doing so.
Throughout this paper, $C, C'$ etc. will denote positive constants whose values may change from line to line; the values of these constants do not depend on the variable (usually `$n$' or `$\ell$') used to index sequences.
For two sequences of real numbers $(a_n)$ and $(b_n)$ 
(resp. two sequences of real-valued functions $(f_n)$ and $(g_n)$ defined on a domain $D$), 
we will write $a_n\lesssim b_n$ (resp. $f_n\lesssim g_n$) to express the relation $a_n\leq Cb_n$ for all $n$
(resp. $f_n(x)\leq C g_n(x)$ for all $x\in D$ and all $n$), 
where $C$ is a constant as above.

For two sequences of real numbers $(a_n)$ and $(b_n)$, we will write $a_n\sim b_n$ to mean $a_n/b_n\to 1$ as $n\to\infty$.
For two sequences $(X_n)$ and $(Y_n)$ of random variables where $X_n$ and $Y_n$ are defined on the same probability space for each $n$, we will write $X_n\sim Y_n$ to mean that $X_n/Y_n\weakc 1$ as $n\to\infty$.
(We also write $X\sim\mu$ to mean that the random variable $X$ has law $\mu$; this, however, should not cause any confusion.)
We write $X_n=O_P(Y_n)$ if the sequence $(X_n/Y_n)$ is tight.
We will write $X_n=o_P(Y_n)$ to mean $X_n/Y_n\weakc 0$ as $n\to\infty$.

We will write `pmf' for `probability mass function.'
We will denote the uniform probability measures on $[0, 1]$ and $[0, 1]^2$ by $\mu_{\circ}$ and $\mu_{\sss\square}$ respectively.
For $x, y\in\bR$, we will write $x\wedge y$ and $x\vee y$ for $\min\{x, y\}$ and $\max\{x, y\}$ respectively.
For a function $f$ taking values in $\bR$, we will denote by $f^+$ and $f^-$ the functions $f\vee 0$ and $-(f\wedge 0)$ respectively.
For $f, g\in L^2([0, 1])$, we will write $\langle f, g\rangle$ for $\int_0^1 f(x)g(x)dx$.

\subsection{Real trees and quotient spaces}\label{sec:cont-limit-descp}
\chsen{A compact metric space $(X,d)$ is called a \emph{real tree} \cite{legall-book,evans-book} if for any two points $x, y\in X$,
(i) there is a unique isometry $f$ from $[0, d(x, y)]$ into $X$ with $f(0)=x$ and $f(d(x, y)) =y$, and
(ii) whenever a function $g:[0, 1]\to X$ is continuous and injective with $g(0)=x$ and $g(1)=y$, we have 
$g([0, 1])=f([0, d(x, y)])$.}

\chsen{For $0 \leq a < b <\infty$, an \emph{excursion} on $[a,b]$ is a continuous function $h \in C([a,b], \bR)$ with $h(a)=0=h(b)$ and $h(t) > 0$ for $t \in (a,b)$. The length of such an excursion is $b-a$. For $l \in(0,\infty)$, let $\cE_l$ be the space of all excursions on the interval $[0,l]$. Given an excursion $h \in \cE_l$, one can construct a real tree as follows. Define the pseudometric $d_h$ on $[0,l]$ as
\begin{equation*}
	d_h(s,t):= h(s) + h(t) - 2 \inf_{u \in [s,t]}h(u)\, , \; \mbox{ for } s,t  \in [0,l].
\end{equation*}
Define the equivalence relation $s \sim t$ if and only if $d_h(s,t) = 0$. Let $[0,l]/\sim$ denote the corresponding quotient space and consider the metric space $\cT(h):= ([0,l]/\sim \, ,\, \bar d_h)$, where $\bar d_h$ is the metric on the equivalence classes induced by $d_h$. Then $\cT(h)$ is a real tree \cite{legall-book,evans-book}. Let $q_h:[0,l] \to \cT(h)$ be the canonical projection and write $\mu_h$ for the pushforward of the Lebesgue measure on $[0,l]$ onto $\cT(h)$ via $q_h$. 
Equipped with $\mu_h$, $\cT(h)$ is now a compact metric measure space.}

\chsen{Since our limit objects will not necessarily be trees, we define a procedure that incorporates ``shortcuts'' (more precisely identification of points) on a real tree. Let $h,g \in \cE_l$ be two excursions, and $\cP \subseteq \bR_+\times \bR_+$ be a locally finite set, and define
\begin{equation*}
	g \cap \cP := \big\{(x,y) \in \cP: 0 \leq x \leq l, \; 0 \leq y < g(x)  \big\} \, .
\end{equation*}
Using these three ingredients, construct a metric space $\cG(h,g,\cP)$ as follows. Let $\cT(h)$ be the real tree associated with $h$ and $q_h: [0,l] \to \cT(h)$ be the canonical projection.  
Note that $|g \cap \cP| < \infty$, and suppose $g \cap \cP  = \{ (x_i, y_i)  : 1\leq i\leq k \}$ for some $k< \infty$. 
For each $i \leq k$, define
$
r(x_i,y_i) := \inf\set{x: x \geq x_i, \; g(x) \leq y_i}
$.
For $1 \leq i \leq k$, identify the points $q_h(x_i)$ and $q_h(r(x_i,y_i))$ in $\cT(h)$. Call the resulting metric space $\cG(h,g,\cP)$. 
Equipping this metric space with the pushforward of the measure $\mu_h$ on $\cT(h)$ makes $\cG(h,g,\cP)$ a compact metric measure space. 
Thus, $\cG(h,g,\cP)$ can be viewed as the metric measure space obtained by adding a finite number of  shortcuts in $\cT(h)$.}

\vskip5pt

\noindent \textbf{Tilted Brownian excursions:} 
Let $(\ve_l(s), s \in [0,l])$ be a Brownian excursion of length $l$. 
Fix $\theta >0$ and let $\tilde \ve_l^\theta$ be an $\cE_l$-valued random variable such that for any bounded continuous function $f : \cE_l \to \bR$,
\begin{equation}
\label{eqn:tilt-exc-def}
\E[f(\tilde \ve_l^\theta)] = { \E\Big[f(\ve_l) \exp\Big(\theta\medint\int_0^l \ve_l(s)ds\Big) \Big] }\big/{\E\Big[\exp\Big(\theta\medint\int_0^l \ve_l(s)ds \Big)\Big]}.
\end{equation}
For simplicity, write $\ve(\cdot) := \ve_1(\cdot)$, $\tilde \ve^\theta (\cdot) := \tilde \ve_1^{\theta} (\cdot)$, and $\tilde \ve_l (\cdot) := \tilde \ve_l^1 (\cdot)$.
In Section \ref{sec:1}, we will use the random objects defined in this section to describe a construction of the limiting random metric measure spaces of interest in this paper.

\subsection{Multiplicative coalescent and the random graph $\cG(\vx,q)$}
\label{sec:smc-def}
A multiplicative coalescent $\big(\vX(t);\, t\in \cA\big)$ is a Markov process with state space 
$
\ldown := \{(x_1,x_2,\ldots): x_1\geq x_2\geq \cdots \geq 0, \sum_i x_i^2< \infty\}
$
endowed with the metric inherited from $l^2$;
here, either $\cA=\bR$ or $\cA=[t_0, \infty)$ for some $t_0\in\bR$.
Its evolution can be described in words as follows:
Fix $\vx\in \ldown$.
Then conditional on $\vX(t) = \vx$, each pair of clusters $i$ and $j$ merge at rate $x_i x_j$ to form a new cluster of size $x_i+x_j$. While this description makes sense for a finite collection of clusters (i.e., $x_i = 0$ for $i> K$ for some finite $K$), Aldous \cite{aldous1997brownian} showed that \chh{it also} makes sense for $\vx\in \ldown$.
We refer the reader to \cite{aldous1997brownian} for a more detailed study of how the multiplicative coalescent is connected to the evolution of the \erdos random graph.

An object closely related to the multiplicative coalescent is the random graph $\cG(\vx,q)$:  Fix vertex set $[n]$, a collection of positive vertex weights $\vx=(x_i, i\in [n])$, and parameter $q>0$. Construct the random graph $\cG(\vx, q)$ by placing an edge between $i\neq j\in [n]$ with probability
$1 - \exp(-q x_i x_j)$,
independently across pairs $\{i, j\}$ with $i\neq j$.
For a connected component $\cC$ of $\cG(\vx, q)$, define 
$\mass (\cC) = \sum_{i \in \cC} x_i={\mu_{\vx}(\cC)}$ 
{(recall the notation from Section \ref{sec:gr-constr})}. 
Rank the components in terms of their masses and let $\cC_i$ be the $i$-th maximal component.

\begin{ass}\label{ass:aldous-basic-assumption}
Consider a sequence of vertex weights $\vx = \vx^{\sss(n)}$ and a sequence $q = q^{\sss(n)}$.
Let $\sigma_k : = \sum_{i \in [n]} x_i^k$ for $k = 2, 3$ and $x_{\max} := \max_{i\in [n]} x_i$. 
Assume that there exists a constant $\lambda \in \bR$ such that as $n \to \infty$,
\begin{equation*}
{\sigma_3}/{(\sigma_2)^3} \to 1, 
\ \ 
q - {(\sigma_2)^{-1}} \to \lambda, 
\ \ \text{ and }\ \ 
{x_{\max}}/{\sigma_2} \to 0.
\end{equation*}
\end{ass}
Note that these conditions imply that $\sigma_2 \to 0$. 
Fix $\lambda \in \bR$ and let $B$ be a standard Brownian motion.  
Define the processes $W_\lambda$ and $\tilde{W}_\lambda$ via
\begin{equation}
\label{eqn:parabolic-bm}
W_\lambda (t) := B(t) + \lambda t - t^2/2\, ,\ \ \text{ and }\ \  
\tilde W_\lambda (t) := W_\lambda (t) - \inf_{s \in [0,t]} W_\lambda (s)\, ,\ \ t\geq 0\, .
\end{equation}
An excursion of $\tilde W_\lambda$ is an interval $(l,r) \subset \bR^+$ such that $\tilde W_\lambda(l)=\tilde W_\lambda(r) = 0$ and $\tilde W_\lambda(t)>0$ for all $t \in (l,r)$. Write $r-l$ for the length of such an excursion. 
Aldous \cite{aldous1997brownian} showed that the lengths of the excursions of $\tilde W_\lambda$ can be arranged in decreasing order as
$\gamma_1(\lambda) > \gamma_2(\lambda) > \ldots > 0$.
Conditional on $\tilde W_\lambda$, let $\pois_i(\lambda)$, $i\geq 1$, be independent random variables where $\pois_i(\lambda)$ is Poisson distributed with mean equal to the area underneath the excursion with length $\gamma_i(\lambda)$. 

\begin{thm}[\cite{aldous1997brownian}]\label{thm:aldous-review}
Under Condition \ref{ass:aldous-basic-assumption}, 	$\left(\mass( \cC_i); i \geq 1 \right) \weakc \mvxi(\lambda):= (\gamma_i(\lambda); i\geq 1)$,  as $n \to \infty$, with respect to the topology on $\ldown$.
Further, 
$\big(\big(\mass( \cC_i),\spls(\cC_i) \big);\, i \geq 1 \big) 
\weakc 
\mvXi(\lambda):=\big( \big(\gamma_i(\lambda), \pois_i(\lambda)\big);\, i\geq 1)$,  as $n \to \infty$, with respect to the product topology.
\end{thm}

\begin{rem}
In \cite[Proposition 4]{aldous1997brownian}, Aldous only proves the first assertion in Theorem \ref{thm:aldous-review}, and the second assertion in Theorem \ref{thm:aldous-review} is proved for the special case of the \erdos random graph \cite[Corollary 2]{aldous1997brownian}.
However, the second convergence in Theorem \ref{thm:aldous-review} follows easily from \cite[Proposition~10]{aldous1997brownian}.
\qed
\end{rem}

\subsection{Scaling limits of components in the critical \erdos random graph}\label{sec:1}
Now we can define the scaling limit of the maximal components in
$\ERRG(n, 1/n+\lambda/n^{4/3})$ derived in \cite{addario2012continuum}. 
Recall the definitions of tilted Brownian excursions and the metric measure space $\cG(h,g,\cP)$ from Section \ref{sec:cont-limit-descp}. 
Let $\mvxi(\lambda)= (\gamma_i(\lambda), i\geq 1)$ be as in Theorem \ref{thm:aldous-review}. Conditional on $\mvxi(\lambda)$, let $\tilde \ve_{\gamma_i(\lambda)}$, $i\geq 1$, be independent tilted Brownian excursions with $\tilde \ve_{\gamma_i(\lambda)}$ having length $\gamma_i(\lambda)$.  Let $\cP_i$, $i\geq 1$, be independent rate one Poisson processes on $\bR_+^2$ that are also independent of $\big(\tilde \ve_{\gamma_i(\lambda)};\,  i\geq 1\big)$. 
Then the random metric measure spaces $\crit_i(\lambda)$ appearing in \eqref{eqn:lim-add-br-go} are given by
\begin{equation}\label{eqn:limit-metric-def}
\crit_i(\lambda)
:= 
\cG\big( 2\tilde \ve_{\gamma_i(\lambda)}, \tilde \ve_{\gamma_i(\lambda)}, \cP_i \big) \, ,  
\ \ \
i\geq 1\, .
\end{equation}
As mentioned in \eqref{eqn:lim-add-br-go}, we will write 
$\vCrit(\lambda) = \big( \crit_i(\lambda) \, ;\,  i\geq 1 \big)$.

\subsection{The BK inequality for connection events}\label{sec:bk-inequality}
\chsen{In this section we recall an application of the van den Berg-Kesten (BK) inequality to connection events;
see \cite{van-den-berg-kesten-bk-inequality} or \cite[Section~2.3]{grimmett-percolation-book} for the statement and the proof of the BK inequality, and \cite{reimer-bkr-inequality} for a more general result due to Reimer.
}

\chsen{Consider $n\geq 1$ and $p_{ij}\in [0, 1]$, $1\leq i<j\leq n$.
Let $K_n$ denote the complete graph on $[n]$, and let $H_n$ denote the random subgraph of $K_n$ where each edge $\{i, j\}$ is deleted independently with probability $1-p_{ij}$, $1\leq i<j\leq n$.
Fix $k\geq 2$ and $i_1, j_1,\ldots, i_k, j_k \in [n]$.
For $1\leq s\leq k$, let $\Pi_s$ be a collection of paths in $K_n$ connecting $i_s$ and $j_s$.
While applying the BK inequality in this paper, $\Pi_s$ will either be the set of all paths in $K_n$ connecting $i_s$ and $j_s$, or the set of self-avoiding paths in $K_n$ with endpoints $i_s$ and $j_s$ having a specified number (say $\ell$) of edges.
For $1\leq s\leq k$, let $A_s$ denote the event that $H_n$ contains a path $\pi_s\in\Pi_s$.
We write $A_1\circ A_2\circ \ldots\circ A_k$ for the event that $H_n$ contains edge-disjoint paths 
$\pi_1\in\Pi_1,$ $\pi_2\in\Pi_2$, \ldots, $\pi_k\in \Pi_k$.
(The symbol `$\circ$' represents disjoint occurrence of events.)
Then the following result is a consequence of the BK inequality:
\begin{lem}\label{lem:bk-inequality}
We have,
$
\pr\big(A_1\circ \ldots \circ A_k \big)
\leq
\prod_{s=1}^k \pr(A_s)
$.
\end{lem}
}

\section{Results: Universality}
\label{sec:res}
Our plan is to extend Theorem \ref{thm:aldous-review} in two stages.
In the first stage, we consider the random graph $\cG(\vx, q)$, and for $i\geq 1$, we view the component $\cC_i$ as a metric measure space endowed with the measure $\mu_{\vx}$.
Then \cite[Theorem~7.3]{SBSSXW14} implies that under Condition \ref{ass:aldous-basic-assumption} and Condition \ref{ass:gen-1} stated below, the sequence $(\cC_i, i\geq 1)$ properly rescaled converges to $\vCrit(\lambda)$ as in \eqref{eqn:lim-add-br-go}.
In the second stage of the extension of Theorem \ref{thm:aldous-review}, we replace each vertex $i \in [n]$ in the graph $\cG(\vx,q)$ with a metric measure space $(M_i, d_i, \mu_i)$; 
we refer to the spaces $(M_i, d_i, \mu_i)$, $i\in [n]$, as `blobs.'
In Theorem \ref{thm:gen-2}, we show that the metric measure space $\bar\cC_i$ now associated with $\cC_i$, under Conditions \ref{ass:aldous-basic-assumption} and \ref{ass:gen-2}, converges to $\crit_i(\lambda)$ after proper rescaling, owing to macroscopic averaging of distances within blobs.

Before moving on to the statements of these results, let us explain the strategy to apply this general theorem to specific random graph models.
Consider a random graph process that undergoes a phase transition in the size of the largest component, and let $\bG_n$ denote the random graph at a point inside the critical window.
We couple $\bG_n$ with a suitable barely subcritical random graph model $\bG_n^{\bullet}$, and we take the components of $\bG_n^{\bullet}$ to be the `blobs.'
This can be seen as a kind of renormalization technique.
We then study the connections that form in going from $\bG_n^{\bullet}$ to $\bG_n$, and we show that the connectivity pattern between blobs can be \emph{approximately} described through the graph $\cG(\vx,q)$ for a suitable $q$, where $x_i$ is a model dependent functional of the $i$-th blob.
We emphasize that at this stage, we are viewing each blob as a single vertex ignoring their internal structure. 
Thus, $\cG(\vx,q)$ should be thought of as approximately describing the blob-level superstructure owing to edges created while going from $\bG_n^{\bullet}$ to $\bG_n$.
Now, once we verify that Condition~\ref{ass:gen-2} is satisfied by the blobs, Theorem~\ref{thm:gen-2} applies to the random graph $\overline\bG_n$ obtained by replacing the vertices of $\cG(\vx,q)$ by the blobs (where the internal structure of the blobs have been incorporated).
Then the final step is to show that $\bG_n$ is well-approximated by $\overline\bG_n$.

\subsection{Blob-level superstructure}
\label{sec:stage-one-blob-level}
Recall the random graph $\cG(\vx,q)$ from Section \ref{sec:smc-def}, and the notation used in Condition \ref{ass:aldous-basic-assumption}.
View $(\cC_i, i \geq 1)$ as metric measure spaces equipped with the graph distance and the weighted measure $\mu_{\vx}$. 
As mentioned in Section~\ref{sec:gh-mc}, when talking about distributional convergence of $\cS^{\bN}$-valued random variables, we will work with the product topology inherited from $d_{\GHP}$.

\begin{ass}\label{ass:gen-1}
Assume that there exist $\eta_0 \in (0,\infty)$ and $r_0 \in (0,\infty)$ such that ${x_{\max}}/{\sigma_2^{3/2+\eta_0}} \to 0$ and  ${\sigma_2^{r_0}}/{x_{\min}} \to 0$
as $n \to \infty$, where $x_{\min}:=\min_{i\in [n]}x_i$.
\end{ass}

\begin{thm}\label{thm:gen-1}
	Under Conditions \ref{ass:aldous-basic-assumption} and \ref{ass:gen-1}, 
	$\big(\scl(\sigma_2, 1 ) \cC_i^{\sss(n)}, i \geq 1 \big) \weakc \vCrit(\lambda)$ as $n\to\infty$. 
\end{thm}

\begin{rem} 
Theorem \ref{thm:gen-1} can be derived as a special case of Theorem \ref{thm:gen-2} stated below.
However, we state Theorem \ref{thm:gen-1} in the above form as this explains the motivation behind us proceeding in two stages as discussed above.
Note also that the scaling limit of the critical \erdos random graph $\ERRG(n, n^{-1}+\lambda n^{-4/3})$ derived in \cite{addario2012continuum} can be recovered from Theorem \ref{thm:gen-1} by taking 
$q=n^{1/3}+\lambda$ and $x_i = n^{-2/3}$ for $i \in [n]$
(which results in $\sigma_2 = n^{-1/3}$). 
\qed
\end{rem}

\subsection{Incorporating the internal structure of the blobs}
\label{sec:inter-blob-distance}
We need the following ingredients:

\noindent{\upshape (a)} {\bf Blob level superstructure: }  A simple finite graph $\cG$ with vertex set $[n]$ and vertex weight sequence $\vx:=(x_i, i \in [n])$.

\vskip3pt

\noindent{\upshape (b)} {\bf Blobs:} A family of compact metric measure spaces 
$\vM := \{(M_i, d_i, \mu_i): i \in [n] \}$, 
where $\mu_i$ is a probability measure for all $i$.

\vskip3pt

\noindent{\upshape (c)} {\bf Blob to blob junction points:} A collection of points $\vX := \{X_{i,j}: i \in [n], j \in [n]\}$ such that $X_{i,j} \in M_i$ for all $i,j$.

\vskip3pt

Using these three ingredients, we define a metric measure space 
$\Gamma(\cG,\vx,\vM,\vX) = (\bar M, \bar d, \bar \mu)$ as follows: Let $\bar M := \bigsqcup_{i \in [n]} M_i$. Define the measure $\bar \mu$ as
\begin{equation}
\bar \mu( A ) = \sum_{i \in [n]} x_i \mu_i(A \cap M_i), \mbox{ for } A \subseteq \bar M.
\label{eqn:barmu-on-full-met}
\end{equation}
The metric $\bar d$ is obtained by using the intra-blob
distance functions $\big(d_i;\, i\in [n]\big)$ together with the graph distance on $\cG$ by putting an edge of length one between the pairs of vertices
$\{ \{X_{i,j}, X_{j,i} \}: \{i,j \} \mbox{ is an edge in } \cG \}$.
Thus, for $x, y \in \bar M$ with $x \in M_{j_1}$ and $y \in M_{j_2}$, 
\begin{equation}\label{eqn:44}
\bar d(x,y) = \inf\big\{ k +  d_{j_1}(x, X_{j_1,i_1}) + \sum_{\ell=1}^{k-1} d_{i_\ell}(X_{i_\ell, i_{\ell-1}}, X_{i_\ell, i_{\ell+1}}) + d_{j_2}(X_{j_2, i_{k-1}}, y) \big\},
\end{equation}
where the infimum is taken over $k\geq 1$ and all paths $(i_0, i_1,\ldots,i_{k-1}, i_k)$ in $\cG$ with $i_0 =j_1$ and $i_k = j_2$. 
(Here, the infimum of an empty set is understood to be $+\infty$.
This corresponds to the case where $j_1$ and $j_2$ do not belong to the same component of $\cG$.)

The above is a deterministic procedure for creating a new metric measure space. 
Now assume that we are provided with a parameter sequence $q^{\sss(n)}$, weight sequence $\vx^{\sss(n)} := (x_i^{\sss(n)}, i \in [n])$, and the family of metric measure spaces $\vM^{\sss(n)} := \big((M_i^{\sss(n)}, d_i^{\sss(n)}, \mu_i^{\sss(n)}), i \in [n]\big)$.
As before, we will suppress the dependence on $n$. 
Let $\cG(\vx, q)$ be the random graph defined in Section~\ref{sec:smc-def} constructed using the weight sequence $\vx$ and the parameter $q$.   Let $(\cC_i, i \geq 1)$ denote the  connected components of $\cG(\vx,q)$ ranked in decreasing order of their masses. 
Let $X_{i,j}$, $i,j \in [n]$,  be independent random variables (that are also independent of $\cG(\vx,q)$) such that for each fixed $i$, $X_{i,j}$, $j\in [n]$, are i.i.d. $M_i$-valued random variables that are $\mu_i$ distributed. 
Let $\vX = (X_{i,j},\ i,j \in [n])$.

\chsen{We define $\bar \cG(\vx,q,\vM) := \Gamma\big( \cG(\vx, q), \vx, \vM, \vX \big)$.}
Let $\bar \cC_i$ be the component in $\bar \cG(\vx, q, \vM)$ corresponding to the $i$-th maximal component $\cC_i$ in $\cG(\vx,q)$, viewed as a compact metric measure space as explained in \eqref{eqn:barmu-on-full-met} and \eqref{eqn:44}. 
Define the quantities
\begin{gather}
u_{i} := \int_{M_i \times M_i} d_i(x,y)  \mu_i (dx) \mu_i (dy), \qquad i\in [n],
\label{eqn:uik-def}\\
\mbox{$\tau := \sum_{i \in [n]} x_i^2 u_i$ and ${\diam_{\max}} := \max_{i \in [n]} \diam(M_i)$.}
\label{eqn:tau-dmax-def}
\end{gather}
\chh{Note that $u_i$ is the expectation of {the distance between} two blob-to-blob junction points in blob $i$ and $\tau/\sigma_2$ is the weighted average of these ``typical'' distances. }
\begin{ass}
	\label{ass:gen-2}
	There exist $\eta_0 \in (0,\infty)$ and $r_0 \in (0,\infty)$ such that
	{\upshape (a)} Condition \ref{ass:gen-1} holds, and
	{\upshape (b)}  as $n \to \infty$, $ (\diam_{\max}\cdot \sigma_2^{3/2-\eta_0})/{ (\tau + \sigma_2)} \to 0$ and  $ {\sigma_2 x_{\max} {\diam_{\max}}}/{\tau} \to 0$. 	
\end{ass}
\begin{thm}
	\label{thm:gen-2}
	Under Conditions~\ref{ass:aldous-basic-assumption} and \ref{ass:gen-2},
	\begin{equation*}
	\Big(\scl\Big(\frac{\sigma_2^2}{\sigma_2 + \tau}, 1 \Big) \bar \cC_i\, ;\, i \geq 1 \Big) \weakc \vCrit(\lambda), \mbox{ as } n \to \infty.
	\end{equation*}
\end{thm}

\begin{rem} \ch{Depending on whether $\lim_{n \to \infty} {\tau}/{\sigma_2}\in [0, \infty)$ or equals $\infty$, the above theorem deals with different scales.
If each blob is a fixed connected graph $G$, then $\lim_{n\to\infty}\tau/\sigma_2\in[0,\infty)$.
(The critical \erdos random graph corresponds to the case where $G$ is a single vertex.)
In \chh{the} applications below, $\sigma_2\sim n^{\delta-1/3}$ while $\tau \sim n^{2\delta-1/3}$ which corresponds to the case where the above limit equals $\infty$.}
\qed	
\end{rem}

\section{Results for the associated random graph models}
\label{sec:res-models}
Consider a sequence of edge weights $\mvbeta_n =\big(\beta_{ij }^{\sss(n)}\, ;\, i, j\in [n]\big)$, $n\geq 1$, as in  Section~\ref{sec:def-graphon-percolation}. 
Let $Q_{ij}^{\sss(n)} =\big((i-1)/n\, , \, i/n \big)\times \big((j-1)/n\, , \, j/n \big)$, $1\leq i, j\leq n$.
Let $W_n:[0, 1]^2\to [0, \infty)$ be given by
\begin{align}\label{eqn:2}
W_n(x, y)
=
\left\{
\begin{array}{l}
	\beta_{ij }^{\sss(n)}\, , \text{ if } (x, y)\in Q_{ij}^{\sss(n)},\ \ 1\leq i,j\leq n,\\
	0,\text{ otherwise}.
\end{array}
\right.
\end{align}

\subsection{The main result for graphons}
Recall the notation introduced in Section~\ref{sec:def-graphon-percolation}. 

\begin{ass}\label{ass:graphon-1}
The sequence $(W_n\, ;\, n\geq 1)$ is uniformly bounded in $L^3([0, 1]^2)$, i.e.,
\begin{equation}\label{eqn:20}
\sup_{n\geq 1}	\medint\int_{[0, 1]^2} W_n^3
=
\sup_{n\geq 1}\, \frac{1}{n^2}\sum_{i, j\in [n]}(\beta_{ij}^{\sss (n)})^3
<
\infty\, .
\end{equation}
Further, there exists an $L^3$ graphon $W$ that is critical (i.e., $\| T_W \|_{2, 2} = 1$) and a symmetric $H\in L^2([0, 1]^2)$ such that 
\begin{equation}\label{eqn:21}
\big\| 
T_{H-n^{1/3}(W_n-W)}
\big\|_{2,2}
=
\big\| 
T_H - n^{1/3}\big(T_{W_n}-T_W\big)
\big\|_{2,2}
\to 0
\ \ \text{ as }\ n\to\infty\, .
\end{equation}
\end{ass}

Let $D_n:=\bigcup_{i=1}^n Q_{ii}^{\sss(n)}$.
Then it is easy to show that 
$n^{1/3}\|T_{W\cdot\ind_{D_n}}\|_{2,2}\to 0$ as $n\to\infty$ for any $W\in L^3([0, 1]^2)$.
Further, 
$\|T_{H\cdot\ind_{D_n}}\|_{2,2}
\leq
\big(\int_{D_n} H^2\big)^{1/2}
\to 0$ as $n\to\infty$ for any $H\in L^2([0, 1]^2)$.
Consequently, the following relatively easier to check condition implies \eqref{eqn:21}: 
\begin{equation}\label{eqn:21-prime}
	\int_{[0, 1]^2\setminus D_n}
	\bigg[H(x, y) - n^{1/3}\bigg(W_n(x, y)-W(x,y)\bigg) \bigg]^2 dx dy 
	\to 0
	\ \ \text{ as }\ n\to\infty\, .
\end{equation}
We will need some additional regularity conditions on the sequence $(\mvbeta_n)$.

\begin{ass}\label{ass:graphon-2}
There exists $\theta_1\in (0, 1/8)$ such that
\[
\max_{i\in [n]}\, \frac{1}{n}\sum_{j\in [n]}\big( \beta_{ij }^{\sss(n)} \big)^{3/2}
=
o(n^{\theta_1})\, ,
\]
or equivalently, there exists $\theta_0\in(0, 1/12)$ such that
\[
\max_{x\in [0, 1]}\Big(\medint\int_0^1 W_n^{3/2}(x , y )\, dy\Big)^{2/3}
=
o(n^{\theta_0})\, .
\]
\end{ass}
		
\begin{ass}\label{ass:graphon-3}
Let $\theta_0$ be as in Condition~\ref{ass:graphon-2}.
Then there exist $\delta_0\in (0\, ,\, 1/3-2\theta_0)$ and $\varpi_0\in [0\, ,\, 2\delta_0+1/3)$ such that the following hold:
For each $n\geq 1$, there exists $B_n\subseteq [n]$ satisfying
\begin{gather}
|B_n|
=
o(n^{\varpi_0}) \, , \ \ \ \
\sum_{i\in B_n}\sum_{j\in [n]}\big( \beta_{ij }^{\sss(n)} \big)^2 
= 
o(n^{4/3}) \, ,\ \text{ and} 
\label{eqn:31}\\
\#\big\{
(i, j)\, :\, 
i, j\in [n]\setminus B_n \, \text{ and }\,
\beta_{ij }^{\sss(n)}\leq n^{-\delta_0}
\big\}
=
o(n^{1+\delta_0})\, .
\label{eqn:33}
\end{gather}
\end{ass}

Note that Condition~\ref{ass:graphon-3} is trivially satisfied if there exists $\delta_0\in (0\, ,\, 1/3-2\theta_0)$ such that for all large $n$,
$\beta_{ij}^{\sss (n)} > n^{-\delta_0}$ for $1\leq i\neq j \leq n$.

\begin{rem}\label{rem:4}
Before proceeding to the statements of our results, let us explain Conditions~\ref{ass:graphon-1}, \ref{ass:graphon-2}, and \ref{ass:graphon-3}.
We first observe two simple facts:
While studying the metric space structure of the graphs, the labels of the vertices do not play any special role.
Hence, for a sequence of edge weights $(\mvbeta_n)$ satisfying Conditions \ref{ass:graphon-2} and \ref{ass:graphon-3}, if we can find permutations $\pi_n$ on $[n]$ such that Condition~\ref{ass:graphon-1} is satisfied by the sequence of weight sequences 
$\big(\beta_{\sss \pi_n(i), \pi_n(j)  }^{\sss(n)}\, ;\, i,j\in[n]\big)$, then the conclusion in Theorem~\ref{thm:scaling-limit-graphon} stated below continues to hold.
Next, note that if $W$ is not a critical graphon, then we can simply consider the edge weights 
$c\beta_{ij}^{\sss(n)}$ where $c=1/\|T_W\|_{2,2}$.
So, the assumption $\|T_W\|_{2,2}=1$ enables us to avoid having this normalizing constant all throughout the paper.

Now, the \erdos random graph $\ERRG(n, t_n/n)$ exhibits critical behavior if there exists $\lambda\in\bR$ such that $n^{1/3}(t_n -1)\to\lambda$ as $n\to\infty$.
Analogously, \eqref{eqn:21} can be written as
\begin{align}\label{eqn:43}
	n^{1/3}(W_n-W)\to H\ \ \text{ as }\ \ n\to\infty\, ,
\end{align}
where the convergence is in the sense of \eqref{eqn:21}.
This ensures that the graphs $\graphonn_n$ are in the critical window.
Here, the critical graphon $W$ plays the role of the critical threshold $t_c=1$ for the \erdos random graph, and the function $H\in L^2([0, 1]^2)$ replaces $\lambda$--the location inside the critical window.

There is a huge body of work surrounding the convergence of graphons. 
Two important topologies with respect to which convergence is studied are the ones induced by the cut norm and the cut metric \cite{frieze-kannan, borgs-et-al-I, janson-cut-norm}.
It is easy to see that \eqref{eqn:21} implies that the convergence in \eqref{eqn:43} holds with respect to the cut norm, and thus the latter induces a weaker topology than the one considered in \eqref{eqn:21};
see \cite[Lemma~E.6]{janson-cut-norm} (this result is stated for bounded graphons; however, this particular implication does not require the boundedness assumption).
It is not clear to us if convergence in the cut norm in \eqref{eqn:43} would, in general, be a strong enough assumption to deduce the metric scaling limit from.
Note however that by \cite[Lemma~E.6]{janson-cut-norm}, in the special case where 
$H\in L^{\infty}([0, 1]^2)$ and $\sup_{n, x, y}n^{1/3}|W_n(x, y)-W(x, y)|<\infty$,
convergence in \eqref{eqn:43} with respect to the cut norm is equivalent to \eqref{eqn:21}.

The assumption that $W$ is an $L^3$ graphon and the condition \eqref{eqn:20} are somewhat similar to the `finite third moment condition' imposed on the degree distribution while studying the configuration model and uniform (simple) graphs with given degree sequence;
it is known 
\cite{dhara-hofstad-leeuwaarden-sen-cm-finite-third-moment, sd-rvdh-jvl-ss-cm-component-scaling-infinite-third-moment,
SB-SD-vdH-SS-cm-metric-infinite-third-moment, 
bhamidi-sen-vacant-set, 
goldschmidt-conchon-stable-graph,
SB-SD-vdH-SS-global-lower-mass}
that for the latter two models, under the finite third moment condition, one gets a Brownian scaling limit inside the critical window, whereas under finite $p$-th moment with $p\in (2, 3)$ but infinite third moment condition, the scaling limit inside the critical window is of a different nature.
Theorems~\ref{thm:scaling-limit-graphon} and \ref{thm:component-limit-graphon} stated below establish Brownian scaling limits for $\graphonn_n$ when $W\in L^3([0, 1]^2)$ and \eqref{eqn:20} holds.
It would be interesting to see if the $L^3$-boundedness condition on the graphons can be relaxed to just the leading eigenfunctions of $T_W$ and $T_{W_n}$ being uniformly bounded in $L^3([0, 1])$ (together with weakening the topology of convergence in \eqref{eqn:21} as discussed above).
Note however that even the more basic problem of whether the linear growth of the giant holds slightly above criticality in this model under the assumption of the leading eigenfunction being in $L^3([0, 1])$ remains open to date; 
see \cite[Conjecture~15.7]{BBSJOR07}.

Similarly, we expect that if the leading eigenfunctions are in $L^p([0, 1])\setminus L^3([0, 1])$ for some $p\in (2, 3)$, then the metric scaling limit in this model will be rescaled versions of the scaling limits arising in \cite{bhamidi-hofstad-sen-multiplicative,SB-SD-vdH-SS-cm-metric-infinite-third-moment}. 
We will pursue this in future work.

Condition~\ref{ass:graphon-2} is similar in spirit to the assumption in \cite[Display~(3.11)]{BBSJOR07}.
We work with the smaller exponent `$3/2$' as opposed to the exponent `$2$' in \cite[Display~(3.11)]{BBSJOR07}.
Also, we allow the possibility $\sup_{x\in [0, 1]}\int_0^1 W^{3/2}(x ,y)dy = \infty$ provided 
$\max_{x\in [0, 1]}\int_0^1 W_n^{3/2}(x ,y)dy$ grows in a controlled way.

Condition~\ref{ass:graphon-3} is a technical assumption that allows us to couple the random graph $\tildegraphon_n$ with the multiplicative coalescent and apply Theorem~\ref{thm:gen-2}.
In applications, $B_n$ is to be taken as the set of vertices $i\in [n]$ such that there are too many $j\in [n]$ for which $\beta_{ij}^{\sss (n)}$ is too small.
The assumption in \eqref{eqn:31} will ensure that the vertices in $B_n$ do not have a significant influence on the intrinsic geometry of the maximal components in $\tildegraphon_n$.
The condition \eqref{eqn:33} says that there are not many pairs $(i, j)$ with both $i$ and $j$ from outside of $B_n$ such that $\beta_{ij}^{\sss (n)}$ is too small.
\qed
\end{rem}

Since $W$ is a symmetric function in $L^3([0, 1]^2)\subseteq L^2([0, 1]^2)$, $T_W$ is a compact and self-adjoint operator on $L^2([0, 1])$. 
Now, $\|T_W\|_{2, 2}=1$ by assumption, and consequently, either $+1$ or $-1$ is an eigenvalue of $T_W$.
In the present setting, we have the following in addition:

\begin{lem}\label{lem:1}
The following hold under Conditions~\ref{ass:graphon-1}, \ref{ass:graphon-2}, and ~\ref{ass:graphon-3}:
\begin{enumeratea}
\item\label{item:21}
$\|T_W\|_{2, 2}=1$ is an eigenvalue of $T_W$.
\item\label{item:22}
The eigenspace corresponding to $1$ is of the form $\{c\psi\, :\, c\in\bR\}$ for an eigenfunction $\psi$ that satisfies $\psi(x)>0$ for $\mu_{\circ}$ a.e. $x$ and $\|\psi\|_2=1$.
(Clearly, such a $\psi$ is unique.)
Further, $\int_0^1\psi^3 <\infty$.
\item\label{item:23}
There exists $\Delta_W\in(0, 1]$ such that $|\lambda|\leq 1-\Delta_W$ for each eigenvalue $\lambda\neq 1$ of $T_W$.
In other words, $-1$ is not an eigenvalue of $T_W$.
\end{enumeratea}
\end{lem}

We defer the proof of Lemma~\ref{lem:1} to Section~\ref{sec:g1}.
Define
\begin{align}\label{eqn:def-alpha-chi-zeta}
\alpha:=\frac{1}{(\int_0^1\psi)^2}\, ,\;\;
\chi:=\frac{\int_0^1 \psi^3}{(\int_0^1\psi)^3}\, , \ \text{ and }\
\zeta:=\frac{\langle \psi, T_H\psi \rangle}{(\int_0^1\psi)^2} \, .
\end{align}

\begin{thm}\label{thm:scaling-limit-graphon}
For $i\geq 1$, let $\cC_i(\graphonn_n)$ be the $i$-th largest component in $\graphonn_n$, and view it as a metric measure space using the counting measure. 
Under Conditions~\ref{ass:graphon-1}, \ref{ass:graphon-2}, and \ref{ass:graphon-3},
\begin{equation}\label{eqn:65}
\Big(
\scl\Big( \frac{\chi^{2/3} }{\alpha n^{1/3}}\, ,\, \frac{\chi^{1/3}}{n^{2/3}} \Big) 
\cC_i(\graphonn_n)\, ,
\ i \geq 1 \Big) 
\weakc 
\vCrit\Big(\frac{\zeta}{\chi^{2/3}}\Big)
\ \ \text{ as }\ \ n\to\infty .
\end{equation}
\end{thm}

Theorem~\ref{thm:scaling-limit-graphon} implies the distributional convergence of the scaled sizes of the  maximal components.
However, it does not immediately show that the number of cycles in each maximal component forms a tight sequence, since there may be cycles of length $o(n^{1/3})$ that do not show up in the GHP limit.
The next result shows the aforementioned tightness, and in fact, its proof will show that there are no cycles of length  $o(n^{1/3})$ in each maximal component whp.
Let $\Xi(\lambda)$ be as in the statement of Theorem~\ref{thm:aldous-review}.

\begin{thm}\label{thm:component-limit-graphon}
Under Conditions~\ref{ass:graphon-1}, \ref{ass:graphon-2}, and \ref{ass:graphon-3}, 
\begin{equation}\label{eqn:66}
\Big(
\bigg(
n^{-2/3}\chi^{1/3}\big|\cC_i(\graphonn_n)\big|,\, \spls\big(\cC_i(\graphonn_n)\big)
\bigg) \ ;\ 
i \geq 1 
\Big)
\weakc
\Xi\big(\chi^{-2/3}\zeta\big) 
\ \ \text{ as } \ \ n \to \infty\, ,
\end{equation}
with respect to the product topology.
\end{thm}

\begin{rem}
From our proofs, it will follow that Theorem~\ref{thm:component-limit-graphon} remains true if we weaken Conditions~\ref{ass:graphon-2} and \ref{ass:graphon-3} to allow 
$\theta_0\in (0, 1/6)$ (equivalently, $\theta_1\in(0, 1/4)$) and $\delta_0\in(0, 1/3 - \theta_0)$.
The stronger assumptions in Conditions~\ref{ass:graphon-2} and \ref{ass:graphon-3} are really needed for the metric scaling limit in Theorem~\ref{thm:scaling-limit-graphon}.
However, we state both results under the same set of conditions for simplicity.
\end{rem}

\begin{rem}
It is known that $\gamma_1(\lambda)/\lambda \weakc 2$ as $\lambda\to\infty$;
see \cite[Lemma 5.6]{addario2013scaling} for a proof
(see also  \cite{addarioberry-bhamidi-sen-leader} for the analogue of this result for the multiplicative coalescent in the regime where the scaling limit is a pure-jump process).
Thus, in the special case $H = \lambda W$ where $\lambda\in\bR$, \eqref{eqn:66} implies that
\begin{align}\label{eqn:66-a}
\lim_{\lambda\to\infty}\lambda^{-1}\bigg(\lim_{n \to \infty} n^{-2/3} |\cC_1(\graphonn_n)|\bigg)
=
2\big(\medint{\int}_0^1 \psi\big) \cdot \big(\medint{\int}_0^1 \psi^3\big)^{-1}\, ,
\end{align}
where both limits are to be interpreted in the distribution sense.
The analogue of \eqref{eqn:66-a} in the purely supercritical regime was proved in \cite[Theorem~3.17~(ii)]{BBSJOR07}.
(The constant $c_0$ appearing in \cite[Display~(3.12)]{BBSJOR07} is missing in \eqref{eqn:66-a} since we work with the normalization $\|T_W\|_{2,2}=1$.
Similarly, the factor $\int \psi^2$ in \cite[Display~(3.12)]{BBSJOR07} is absent in \eqref{eqn:66-a} because we use the normalization $\|\psi\|_2 = 1$.)
Thus, \eqref{eqn:66-a} generalizes \cite[Theorem~3.17~(ii)]{BBSJOR07} to the right side of the critical window.
\end{rem}

\subsection{Some special cases of edge weights}\label{sec:special-cases}
We discuss four special cases of edge weights below:

\subsubsection{The stochastic block model}\label{sec:stochastic-block-model}
Theorem~\ref{thm:scaling-limit-graphon} can be applied to the important special case of the stochastic block model (SBM).
However, for the SBM, instead of working with compact operators, one can just work with finite dimensional matrices.
So, we can rework the proof of Theorem~\ref{thm:scaling-limit-graphon} in the simpler setting of matrices, and prove a result for the SBM that is slightly more general than the one we get as a corollary of Theorem~\ref{thm:scaling-limit-graphon}.
We will state this result without proof in Section~\ref{sec:disc}.

\subsubsection{The edge weights $W(i/n,\, j/n)$}\label{sec:w-i/n-j/n}
Suppose $K\in L^3([0, 1]^2)$ is a nonnegative, symmetric kernel, and assume that
\begin{align}\label{eqn:37}
	W=\|T_K\|_{2,2}^{-1}\cdot K\, .
\end{align}
When $K$ has reasonable continuity properties, a natural class of edge weights to study would be $W(i/n,\, j/n)$, $1\leq i\neq j\leq n$.
Incorporating the location in the critical window, we consider the more general setting where the edge weights are given by
\begin{align}\label{eqn:444}
	\beta_{ii}^{\sss(n)}=0
	\ \ \text{ for }\ i\in[n]
	\ \text{ and }\
	\beta_{ij}^{\sss (n)}
	=
	\Big(
	W\big(\tfrac{i}{n}\, ,\, \tfrac{j}{n}\big)
	+
	n^{-1/3}
	H\big(\tfrac{i}{n}\, ,\, \tfrac{j}{n}\big)
	\,\Big)\vee 0
	\ \ \text{ for }\ 1\leq i\neq j\leq n\, .
\end{align}
For simplicity, we assume that either $H=\lambda W$ for some $\lambda\in\bR$, or $H$ is continuous on $[0, 1]^2$.

Clearly, if $K$ is bounded away from zero and H\"{o}lder$(\gamma)$ on $[0, 1]^2$ with $\gamma>1/3$, then Conditions~\ref{ass:graphon-1}, \ref{ass:graphon-2}, and \ref{ass:graphon-3} are satisfied by the edge weights in \eqref{eqn:444}, and consequently, Theorems~\ref{thm:scaling-limit-graphon} and \ref{thm:component-limit-graphon} are applicable in such cases.
Some specific examples where Theorems~\ref{thm:scaling-limit-graphon} and \ref{thm:component-limit-graphon} are applicable are
\begin{align}\label{eqn:5}
	\eta+(x\vee y)^a\, , \ \  \text{ or }\ \
	\eta\vee(x\vee y)^a\, , \ \  \text{ or }\ \
	(\eta+ x+ y)^a\, , \ \ \text{ or }\ \
	\eta+(x+y)^a\, ,  \ \ \text{ for any }\ \ a, \eta>0.
\end{align}
(Note that the above H\"{o}lder condition is not met in some of these examples.)

Now, suppose $K$ is either unbounded, or is not bounded away from zero, or both.
We briefly discuss applications of Theorems~\ref{thm:scaling-limit-graphon} and \ref{thm:component-limit-graphon} to such graphons via a few simple examples.
Consider the family of unbounded kernels
\begin{align}\label{eqn:30}
	K(x, y)= 
	(x\vee y)^{-a} , 
	\ \ \text{ or } \ \
	(x+y)^{-a},
	\ \ \text{ or } \ \
	\big(x^a+y^a\big)^{-1},
	\ \ \ 0<a<2/3\, .
\end{align}
See Sections~3.6, 16.1, 16.3, and 16.6 of \cite{BBSJOR07} for a survey of the works related to inhomogeneous random graphs with kernels that are functions of $(x\vee y)$.
Note that the kernels in \eqref{eqn:30} are not in $L^3([0, 1]^2)$ if $a\geq 2/3$.
It can be easily checked that Conditions~\ref{ass:graphon-1}, \ref{ass:graphon-2}, and \ref{ass:graphon-3} are satisfied by the edge weights in \eqref{eqn:444} when $K$ is as in \eqref{eqn:30}, and consequently, Theorems~\ref{thm:scaling-limit-graphon} and \ref{thm:component-limit-graphon} are applicable in these cases.
Next, take the unbounded kernels $|x-y|^{-a}$, which are in $L^3([0, 1]^2)$ if $0<a<1/3$.
Then Theorems~\ref{thm:scaling-limit-graphon} and \ref{thm:component-limit-graphon} apply when
\begin{align}\label{eqn:55}
	K(x, y):= |x-y|^{-a} ,\ \ 0<a<1/6\, .
\end{align}

Some examples of kernels $K$ that are not bounded away from zero to which Theorems~\ref{thm:scaling-limit-graphon} and \ref{thm:component-limit-graphon} are applicable are
\begin{gather}
K(x, y)=(x\wedge y)^a\, ,\ \ a>0\, ,
\label{eqn:451}
\\
K(x, y)=(x\vee y)^a\, ,\ \ 0\leq a<1\, ,\ \  \text{ and }
\label{eqn:452}
\\
K(x, y)=(x + y)^a\, ,\ \ 0\leq a<1\, .
\label{eqn:453}
\end{gather}
Now, unlike \eqref{eqn:5} or \eqref{eqn:451}, Theorem~\ref{thm:scaling-limit-graphon} cannot be used directly to get the scaling limit for $K(x, y) = (x\vee y)^a$ or $(x+y)^a$ when $a\geq 1$ and the edge weights are as in \eqref{eqn:444}.
However, Theorem~\ref{thm:scaling-limit-graphon} is applicable to these kernels when the edge weights are 
$\beta_{ij}^{\sss(n)}\vee n^{-\delta}$ where $\beta_{ij}^{\sss(n)}$ is as in \eqref{eqn:444} and $\delta\in (1/6, 1/3)$.
Thus, to get the metric scaling limit in such cases, it would suffice to establish the scaling limit of the maximal component sizes and their surplus edges.
Take, for example, the kernel $x\vee y$.
Note that $\|T_{x\vee y}\|_{2, 2}=:z_0$ is the unique positive solution to 
$\tanh \big(1/\sqrt{z}\big)=\sqrt{z}$, 
and the corresponding (one-dimensional) eigenspace is spanned by
\begin{align}\label{eqn:67}
	\tilde\psi(x)
	=
	\frac{\sqrt{2}}{\cosh\big(1/z_0^{1/2}\big)}\cdot\cosh\big(x/z_0^{1/2}\big)\, ,\ \ \ x\in [0, 1].
\end{align}
Recall the definition of $\Xi(\lambda)$ from the statement of Theorem~\ref{thm:aldous-review}.
\begin{conj}\label{conj:x-max-y}
	Fix $\lambda\in\bR$, and consider $\graphonn_n$ constructed using the edge weights
	\begin{align}\label{eqn:77}
		\tilde\beta_{ij}^{\sss(n)}
		=
		\frac{1}{z_0}\Big(1+\frac{\lambda}{n^{1/3}}\Big)\cdot \bigg( \frac{i}{n}\vee \frac{j}{n} \bigg)
		\, ,\ \ \
		1\leq i\neq j\leq n\, .
	\end{align}
	Then
	\begin{equation}\label{eqn:91}
		\Big(\,
		\Big(
		\frac{\tilde\chi^{1/3}}{n^{2/3}}\big|\cC_i(\graphonn_n)\big|\, ,\
		\spls\big(\cC_i(\graphonn_n)\big)
		\Big)\,
		;\ i \geq 1 
		\Big)
		\weakc
		\Xi\left(\frac{\tilde\zeta}{\tilde\chi^{2/3}}\right) 
		\ \ \text{ as } \ \ n \to \infty
	\end{equation}
	with respect to the product topology, where 
	$\tilde\chi=(\int_0^1\tilde\psi^3)(\int_0^1\tilde\psi)^{-3}$ and 
	$\tilde\zeta=\lambda(\int_0^1\tilde\psi)^{-2}$
	with $\tilde\psi$ as in \eqref{eqn:67}.
\end{conj}
Now, for any $\delta\in (1/6,\, 1/3)$, Theorem~\ref{thm:scaling-limit-graphon} and Theorem~\ref{thm:component-limit-graphon} show that the convergences in \eqref{eqn:65} and \eqref{eqn:66} hold when the edge weights are $\tilde\beta_{ij}^{\sss (n)}\vee n^{-\delta}$, where $\tilde\beta_{ij}^{\sss (n)}$ is as in \eqref{eqn:77}.
This can be used in conjuction with the convergence in \eqref{eqn:91} and a coupling argument to get the metric scaling limit when the edge weights are as in \eqref{eqn:77}.
Thus, the task of deriving the metric scaling limit in such cases can be reduced to the easier problem of establishing the scaling limit of the component sizes and their surplus edges.
We will briefly revisit this point in Section~\ref{sec:disc}.

Finally, as examples of kernels that are neither bounded nor bounded away from zero, consider
\begin{gather}
	K(x, y)= 
	\frac{1}{(x\vee y)^{a}} -1\, ,
	\ \ \text{ or }\ \ \
	\Big(\frac{1}{x\vee y} - 1\Big)^a,
	\ \ \text{ or }\ \ \
	\Big(\frac{2}{x+y} - 1\Big)^a , 
	\ \ \ 0<a<2/3\, .\label{eqn:30-a}
\end{gather}
As discussed in \cite[Section~16.3]{BBSJOR07}, the first two kernels in \eqref{eqn:30-a} with $a=1$ correspond to the CHKNS model \cite{callaway-hopcroft-kleinberg-newman-strogatz} and a related process due to Durrett \cite{durrett-chkns-model}.
Note also that the kernels in \eqref{eqn:30-a} are not in $L^3([0, 1]^2)$ when $a\geq 2/3$.
Here, again it is routine to verify that Conditions~\ref{ass:graphon-1}, \ref{ass:graphon-2}, and \ref{ass:graphon-3} are satisfied.

\subsubsection{Edge weights sampled using uniform random variables}
For a graphon $W$, a natural collection of edge weights that is widely studied in the statistics literature is  $W(U_i, U_j)$, $1\leq i\neq j\leq n$, where $U_1,\ldots, U_n$ are i.i.d. $\text{Uniform[0, 1]}$ random variables.
As explained at the beginning of Remark~\ref{rem:4}, a permutation of the vertex labels does not make any difference, so we can instead take the edge weights as 
$W(V^{\sss (n)}_i, V^{\sss(n)}_j)$, $1\leq i\neq j\leq n$, 
where $V^{\sss(n)}_1\leq\cdots\leq V^{\sss(n)}_n$ are $U_1,\ldots, U_n$ arranged in nondecreasing order.
We will work in a slightly more general setting and take the edge weights as
\begin{align}\label{eqn:23}
	\beta_{ii}^{\sss(n)}=0
\ \ \text{ for }\ i\in[n]
\ \text{ and }\
\beta_{ij}^{\sss(n)}
=
\Big(
W\big( V^{\sss (n)}_i, V^{\sss(n)}_j \big) 
+ 
\frac{1}{n^{1/3}}\cdot H\big( V^{\sss (n)}_i, V^{\sss(n)}_j \big)
\Big)\vee 0\, ,\ \ \
1\leq i\neq j\leq n\, ,
\end{align}
where $H:[0, 1]^2\to\bR$ is symmetric and Borel measurable.

\subsubsection{Edge weights that provide the best $L^2$-approximation}
For $f\in L^1([0, 1]^2)$, let $\dashint_{\sss Q_{ij}^{\sss (n)}}f:=n^2\int_{\sss Q_{ij}^{\sss (n)}} f$ 
denote the average of $f$ over $Q_{ij}^{\sss (n)}$, $i, j\in [n]$.
Then, in view of Condition~\ref{ass:graphon-1} and \eqref{eqn:21-prime}, another natural choice for the edge weights would be 
\begin{align}\label{eqn:56}
	\beta_{ii}^{\sss(n)}=0
	\ \ \text{ for }\ i\in[n]
	\ \text{ and }\
	\beta_{ij}^{\sss (n)}
	=
	\Big(\
	\dashint_{\sss Q_{ij}^{\sss (n)}}W
	+
	\frac{1}{n^{1/3}}
	\dashint_{\sss Q_{ij}^{\sss (n)}}H
	\,\Big)\vee 0
	\ \ \text{ for }\ 1\leq i\neq j\leq n\, ,
\end{align}
since for given $W$ and $H$, the function $W_n$ as given by \eqref{eqn:2} corresponding to the edge weights in \eqref{eqn:56} gives the best $L^2$-approximation for $(W+n^{-1/3}H)$ among all functions that are zero on $D_n$, and constant and nonnegative on $Q_{ij}^{\sss(n)}$ for $1\leq i\neq j\leq n$.

In the next two sections, we formulate regularity conditions for $W$ and $H$ under which Theorems~\ref{thm:scaling-limit-graphon} and \ref{thm:component-limit-graphon} are applicable to the collections of edge weights in \eqref{eqn:23} and \eqref{eqn:56}.
In these regularity assumptions, we eschew complete generality in favor of relatively simpler conditions.

\subsection{Special case: the edge weights in \eqref{eqn:23}}\label{sec:2}
For $\eps>0$ and $K:[0, 1]^2\to[0, \infty)$, let $\osc_{\eps; K}:[0, 1]^2\to[0, \infty]$ be given by
\[
\osc_{\eps; K}(x, y)
:=
\sup\big\{
|K(x', y')-K(x, y)|\, :\,
(x', y')\in [0, 1]^2\ \text{ and }\ 
|x-x'|+|y-y'| \leq |\log\eps|\sqrt{\eps}
\big\}\, .
\]

\begin{ass}\label{ass:W-Ui-Uj}
Assume that $W$ is an $L^3$ graphon satisfying the following:
\begin{enumeratea}
\item\label{item:1}
For each $\eps>0$, there exists a Borel $\cA_{\eps}\subseteq [0, 1]^2$ such that
\[
\medint\int_{\widetilde\cA_{\eps}} W^2(x, y) dx dy
+
\medint\int_{[0, 1]^2\setminus\cA_{\eps}} \osc_{\eps; W}^2(x, y) dx dy
=
o(\eps^{2/3})\ \text{ as }\ \eps\to 0\, ,
\]
where $\widetilde\cA_{\eps}$ is the set of points in $[0, 1]^2$ such that there is at least one point in $\cA_{\eps}$ at distance at most $(|\log\eps|\sqrt{\eps})$ from them.
\item\label{item:2}
There exists $p_1>8$ such that
\[
\medint\int_0^1
\Big( 
\medint\int_0^1 W^{3/2}(x, y)dy
\Big)^{p_1} dx
<\infty\, .
\]
\item\label{item:3}
There exist $\delta_0\in (0,\, 1/3-4/(3p_1))$ and $\varpi_0\in [0\, ,\, 2\delta_0+1/3)$ such that the following hold:
For each $\eps>0$, there exists a Borel $\cB_{\eps}\subseteq [0, 1]$ satisfying
\begin{gather}
 \mu_{\circ}(\cB_{\eps})
=
o(\eps^{1-\varpi_0}) \, ,\ \ \
\medint\int_{\cB_{\eps}\times [0, 1]} W^2(x, y)dx dy 
= 
o(\eps^{2/3})\, ,\ \text{ and}  
\label{eqn:34}
\\
\mu_{\sss \square}\big(
\big\{
(x, y) \in \big([0, 1]\setminus\cB_{\eps}\big)^2\, :\, 
W(x, y)\leq \eps^{\delta_0}
\big\}
\big)
=o(\eps^{1-\delta_0})\ \ \text{ as }\eps\to 0\, .
\label{eqn:35}
\end{gather}
\end{enumeratea}
\end{ass}

Note that if $\inf_{x, y\in [0, 1]}W(x, y)>0$, then Condition~\ref{ass:W-Ui-Uj}~\eqref{item:3} is satisfied trivially with 
$\cB_{\eps}=\emptyset$ for all $\eps>0$, 
any $\delta_0\in (0,\, 1/3-4/(3p_1))$, and any $\varpi_0\in [0\, ,\, 2\delta_0+1/3)$.
Condition~\ref{ass:W-Ui-Uj}~\eqref{item:1} is only needed while verifying that the convergence in \eqref{eqn:21} holds in distribution when $H\in L^p([0, 1]^2)$ for some $p>2$ and the edge weights are given by \eqref{eqn:23}; 
see the proof of Theorem~\ref{thm:graphon-Ui-Uj} given in Section~\ref{sec:g13}.
For certain kernels that do not satisfy Condition~\ref{ass:W-Ui-Uj}~\eqref{item:1}, it may still be possible to directly verify that \eqref{eqn:21} holds in distribution; see e.g., the discussion below \eqref{eqn:45}.
We expect that Condition~\ref{ass:W-Ui-Uj}~\eqref{item:1} can be substantially weakened, and we will come back to this point in Section~\ref{sec:disc}.

\begin{thm}\label{thm:graphon-Ui-Uj}
Consider $\graphonn_n$ constructed using edge weights of the form \eqref{eqn:23}.
If $W$ is a critical $L^3$ graphon that satisfies Condition~\ref{ass:W-Ui-Uj} and either
\begin{enumeratea}
\item\label{item:4}
$H=\lambda W$ for some $\lambda\in\bR$, or
\item\label{item:5}
$H^+\in L^3([0, 1]^2)$ and 
$H^-\in L^q([0, 1]^2)$ where $q=q(\delta_0):=3(1-\delta_0)/(1-3\delta_0)$,
\end{enumeratea}
then the convergences in \eqref{eqn:65} and \eqref{eqn:66} hold with the same constants as in \eqref{eqn:def-alpha-chi-zeta}.

If further $\inf_{x, y\in [0, 1]}W(x, y)>0$, then the convergences in \eqref{eqn:65} and \eqref{eqn:66} hold with the same constants as in \eqref{eqn:def-alpha-chi-zeta} for any $H\in L^3([0, 1]^2)$.
\end{thm}

The next result follows easily from Theorem~\ref{thm:graphon-Ui-Uj}.

\begin{cor}\label{cor:holder-graphon}
Suppose $K:[0, 1]^2\to [0, \infty)$ is a symmetric, Borel measurable function that is bounded, bounded away from zero (i.e., $\inf_{x, y\in [0, 1]}K(x, y)>0$), and satisfies
\begin{align}\label{eqn:78}
\medint\int_{[0, 1]^2} \osc_{\eps; K}^2(x, y) dx dy=o(\eps^{2/3})\text{ as }\ \eps\to 0.
\end{align}
Let $W=\|T_K\|_{2,2}^{-1}\cdot K$, and consider $\graphonn_n$ constructed using edge weights of the form \eqref{eqn:23}, where $H\in L^3([0, 1]^2)$.
Then the convergences in \eqref{eqn:65} and \eqref{eqn:66} hold with the same constants as in \eqref{eqn:def-alpha-chi-zeta}.
\end{cor}

In particular, Corollary~\ref{cor:holder-graphon} is applicable to any $K$ that is H\"{o}lder$(\gamma)$ on 
$[0, 1]^2$ with $\gamma>2/3$ and is bounded away from zero.
However, \eqref{eqn:78} may hold even when the above H\"{o}lder condition does not.
For example, Corollary~\ref{cor:holder-graphon} applies to the kernels in \eqref{eqn:5}.

Now, as we did in Section~\ref{sec:w-i/n-j/n}, consider $W$ as in \eqref{eqn:37}, where 
$K:[0, 1]^2\to [0, \infty)$ is a symmetric function in $L^3([0, 1]^2)$ that is either unbounded, or is not bounded away from zero, or both.
We briefly discuss applications of Theorem~\ref{thm:graphon-Ui-Uj} to such graphons via a few examples.

Consider the family of unbounded kernels
\begin{align}\label{eqn:4}
K(x, y)= (x\vee y)^{-a}\, , \ \ \ \ 0<a<1/3\, .
\end{align}
Note that the kernel $(x\vee y)^{-a}$ is unbounded and in $L^3([0, 1]^2)$ iff $0<a<2/3$.
However, Theorem~\ref{thm:graphon-Ui-Uj} applies only when $0<a<1/3$.
To see that Theorem~\ref{thm:graphon-Ui-Uj} can be applied when $K$ is as in \eqref{eqn:4}, it is enough to verify that $K$ satisfies Condition~\ref{ass:W-Ui-Uj}.
To this end, choose $b\in\big( \big(3(1-a)\big)^{-1},\, 1/2\big)$, and for each $\eps\in(0, 1)$, set
\[
\cA_{\eps}=[0, \eps^b]\times [0, \eps^b]\ \ \text{ and }\ \ 
\cB_{\eps}=\emptyset\, .
\]
Then Condition~\ref{ass:W-Ui-Uj} is satisfied with these choices of $\cA_{\eps}, \cB_{\eps}$, and any $p_1>8$, $0<\delta_0< 1/3-4/(3p_1)$ and $\varpi_0\in [0\, ,\, 2\delta_0+1/3)$.
Similar calculations will show that Condition~\ref{ass:W-Ui-Uj} is also satisfied when working with the kernels
\begin{align}\label{eqn:6}
(x+y)^{-a}\, ,\ \ 0<a<1/3\, ,
\end{align}
and consequently, Theorem~\ref{thm:graphon-Ui-Uj} is applicable.

Now, let us look at two examples where $K$ is not bounded away from zero.
Consider the family of kernels
\begin{align}\label{eqn:45}
K(x, y)=(x\wedge y)^a\, ,\ \ a>0\, .
\end{align}
The kernel $x\wedge y$ is particularly important to us as it is related to the RGIV model;
see the proof of Theorem~\ref{thm:rgiv-scaling-limit} stated below.
First, consider the case $a>1/6$.
To see that the kernels in \eqref{eqn:45} satisfy Condition~\ref{ass:W-Ui-Uj} when $a>1/6$, we choose 
\[
\delta_0\in \Big(\frac{2a}{3(2a+1)}\,\, ,\ \frac{1}{3}\Big)\, , \ \ \
\varpi_0\in \Big(\big(1-\delta_0/a\big)\vee 0\,\, ,\ 2\delta_0 +1/3\Big),
\] 
and for each 
$\eps\in (0, 1)$, set
$
\cA_{\eps}=\emptyset$ 
and
$\cB_{\eps}=\big[0,\, \eps^{\delta_0/a}\big]
$.
Then Condition~\ref{ass:W-Ui-Uj} is satisfied with these choices of $\cA_{\eps}, \cB_{\eps}, \delta_0, \varpi_0$, and any $p_1>8$ such that $\delta_0< 1/3-4/(3p_1)$.
When $a\in (0, 1/6]$, Condition~\ref{ass:W-Ui-Uj}~\eqref{item:1} is not satisfied.
However, as mentioned right below Condition~\ref{ass:W-Ui-Uj}, for the kernels in \eqref{eqn:45}, we can directly verify that \eqref{eqn:21} holds if $H$ is as in the statement of Theorem~\ref{thm:graphon-Ui-Uj}.
Indeed, the relation 
\[
n^{-1}\sum_{i\in [n]}\big( ( V_i^{\sss (n)} )^a - i^a/n^a \big)^2
=
o_P(n^{-2/3})
\]
and some straightforward calculation will show that convergence in \eqref{eqn:21} holds in distribution whenever $H$ is as in the statement of Theorem~\ref{thm:graphon-Ui-Uj}.
Thus, the conclusion of Theorem~\ref{thm:graphon-Ui-Uj} remains true when $K$ is as in \eqref{eqn:45}.

Next, we consider the kernels in \eqref{eqn:452}.
Fix any $\delta_0\in \big( a/(2+a)\, ,\, 1/3\big)$.
Then Condition~\ref{ass:W-Ui-Uj} is satisfied with this choice of $\delta_0$, 
$\cA_{\eps}=\emptyset, \cB_{\eps}=\emptyset$,
and any  $\varpi_0\in [0,\, 2\delta_0 +1/3)$ and 
any $p_1>8$ such that $\delta_0< 1/3-4/(3p_1)$.
Similar calculations show that Theorem~\ref{thm:graphon-Ui-Uj} is applicable to the kernels in \eqref{eqn:453}.
Further, the analogue of Conjecture~\ref{conj:x-max-y} and the discussion around Conjecture~\ref{conj:x-max-y} are valid here as well.

We now look at a few examples of kernels that are both unbounded and not bounded away from zero.
Consider
\begin{align}\label{eqn:8}
K(x, y)
=
\frac{1}{(x\vee y)^a}-1\, ,\ \ \ 0<a<1/3\, .
\end{align}
Kernels of this form are in $L^3([0, 1]^2)$ when $0<a<2/3$.
As was the case in \eqref{eqn:4} and \eqref{eqn:6}, Theorem~\ref{thm:graphon-Ui-Uj} applies to such kernels when $0<a<1/3$.
To verify Condition~\ref{ass:W-Ui-Uj} for the kernels in \eqref{eqn:8}, choose 
$b\in\big( \big(3(1-a)\big)^{-1},\, 1/2\big)$,
$\delta_0\in (2/9, 1/3)$, 
$\varpi_0\in (1-\delta_0,\, 2\delta_0 +1/3)$, and let
\[
\cA_{\eps}=[0, \eps^b]\times [0, \eps^b]
\ \ \text{ and }\ \ 
\cB_{\eps}=\big[1-2\eps^{\delta_0}/a\, ,\, 1\big]\, .
\]
Then Condition~\ref{ass:W-Ui-Uj} is satisfied with these choices of $\cA_{\eps}, \cB_{\eps}, \delta_0, \varpi_0$, and any $p_1>8$ such that $\delta_0< 1/3-4/(3p_1)$.
Similarly, Theorem~\ref{thm:graphon-Ui-Uj} is applicable when
\[
K(x, y)= 
\Big(\frac{1}{x\vee y} - 1\Big)^a
\ \ \text{ or }\ \ \
\Big(\frac{2}{x+y} - 1\Big)^a , 
\ \ \ 0<a<1/3\, .
\]

\subsection{Special case: the edge weights in \eqref{eqn:56}}
Let us now turn to the edge weights in \eqref{eqn:56}, and formulate sufficient conditions for Theorem~\ref{thm:scaling-limit-graphon} to be applicable.

\begin{ass}\label{ass:W-l2-approx}
Assume that $W$ is an $L^3$ graphon satisfying the following:
\begin{enumeratea}
\item\label{item:1-l2}
We have, 
\begin{align}\label{eqn:41}
\sum_{1\leq j<i\leq n}\
\int_{\sss Q_{ij}^{\sss (n)}}
\Big(W-\dashint_{\sss Q_{ij}^{\sss (n)}} W\Big)^2
	=
o(n^{-2/3})\, .
\end{align}

\item\label{item:2-l2}
Condition~\ref{ass:W-Ui-Uj}~\eqref{item:2} holds.

\item\label{item:3-l2}
Let $p_1$ be as in \eqref{item:2-l2}. 
Then there exist $\delta_0\in (0,\, 1/3-4/(3p_1))$ and $\varpi_0\in [0\, ,\, 2\delta_0+1/3)$ such that the following hold:
For each $n$ there exists $B_n\subseteq [n]$ such that
\begin{gather}
|B_n|
=
o(n^{\varpi_0}) \, ,\ \ \
\medint\int_{ [0, 1]\times \cB_{1/n}} W^2(x, y)dx dy 
= 
o(n^{-2/3})\, ,\ \text{ and}  
\label{eqn:34-l2}
\\
\mu_{\sss \square}\big(
\big\{
(x, y) \in \big([0, 1]\setminus\cB_{1/n}\big)^2\, :\, 
W(x, y)\leq n^{-\delta_0}
\big\}
\big)
=o(n^{\delta_0-1})\ \ \text{ as }n\to \infty\, ,
\label{eqn:35-l2}
\end{gather}
where $\cB_{1/n}:=\bigcup_{i\in B_n}\big((i-1)/n\, ,\, i/n\big)$.
\end{enumeratea}
\end{ass}

If $W$ has weak partial derivatives $D_x W$ and $D_y W$ on the region $\{0<y<x<1\}$, then by the Poincar\'e inequality for functions in Sobolev spaces \cite[Display~(7.45)]{gilbarg-trudinger-pde-book}, for any $1\leq j<i\leq n$,
\begin{align}\label{eqn:68}
\int_{\sss Q_{ij}^{\sss (n)}}
\Big(W-\dashint_{\sss Q_{ij}^{\sss (n)}} W\Big)^2
\leq
Cn^{-2}\int_{\sss Q_{ij}^{\sss (n)}}\Big[\big(D_x W\big)^2 + \big(D_y W\big)^2\Big]\, ,
\end{align}
where $C$ is a universal constant.
Further, we have the trivial bound
$
\int_{\sss Q_{ij}^{\sss (n)}}
\big(W-\dashint_{\sss Q_{ij}^{\sss (n)}} W\big)^2
\leq
\int_{\sss Q_{ij}^{\sss (n)}} W^2
$,
which can be combined with \eqref{eqn:68} to get a sufficient condition for \eqref{eqn:41} to be satisfied.
Finally, if $W$ is H\"{o}lder$(\gamma)$ on $[0, 1]^2$ with $\gamma>1/3$, then $W$ satisfies \eqref{eqn:41}.
We collect these observations in the next lemma.

\begin{lem}\label{lem:21}
An $L^3$ graphon $W$ satisfies \eqref{eqn:41} if either
\begin{enumeratea}
\item\label{item:6}
$W$ is H\"{o}lder$(\gamma)$ on $[0, 1]^2$ with $\gamma>1/3$, or if
\item\label{item:7}
$W$ has weak partial derivatives $D_x W$ and $D_y W$ on the region $\{0<y<x<1\}$ and
\begin{align}\label{eqn:40}
\sum_{1\leq j<i\leq n}
\min\Big\{
\medint\int_{\sss Q_{ij}^{\sss (n)}} W^2\, ,\
\frac{1}{n^2}\medint\int_{\sss Q_{ij}^{\sss (n)}}\Big[\big(D_x W\big)^2 + \big(D_y W\big)^2\Big]
\Big\}
=
o(n^{-2/3})\, .
\end{align}
\end{enumeratea}
\end{lem}

\begin{thm}\label{thm:l-2-approx}
Consider $\graphonn_n$ constructed using edge weights of the form \eqref{eqn:56}.
If $W$ is a critical $L^3$ graphon satisfying Condition~\ref{ass:W-l2-approx} and either
\begin{enumeratea}
\item\label{item:4a}
$H=\lambda W$ for some $\lambda\in\bR$, or
\item\label{item:5a}
$H^+\in L^3([0, 1]^2)$ and 
$H^-\in L^q([0, 1]^2)$ where $q=q(\delta_0):=3(1-\delta_0)/(1-3\delta_0)$,
\end{enumeratea}
then the convergences in \eqref{eqn:65} and \eqref{eqn:66} hold with the same constants as in \eqref{eqn:def-alpha-chi-zeta}.

If further $\inf_{x, y\in [0, 1]}W(x, y)>0$, then the convergences in \eqref{eqn:65} and \eqref{eqn:66} hold with the same constants as in \eqref{eqn:def-alpha-chi-zeta} for any $H\in L^3([0, 1]^2)$.
\end{thm}

The following result is immediate from Theorem~\ref{thm:l-2-approx}:

\begin{cor}
Suppose $W$ is a critical $L^3$ graphon that satisfies $\inf_{x, y\in [0, 1]} W(x, y)>0$.
Assume further that at least one of the following two conditions hold:
\begin{enumeratea}
\item 
$W$ is H\"{o}lder$(\gamma)$ on $[0, 1]^2$ with $\gamma>1/3$.
\item 
$W$ satisfies Condition~\ref{ass:W-Ui-Uj}~\eqref{item:2} and the condition in Lemma~\ref{lem:21}~\eqref{item:7}.
\end{enumeratea}
Consider $\graphonn_n$ constructed using edge weights of the form \eqref{eqn:56}, where $H\in L^3([0, 1]^2)$.
Then the convergences in \eqref{eqn:65} and \eqref{eqn:66} hold with the same constants as in \eqref{eqn:def-alpha-chi-zeta}.
\end{cor}

For some specific examples, note that Theorem~\ref{thm:l-2-approx} is applicable to $W$ of the form \eqref{eqn:37}, where $K$ is any of the kernels in \eqref{eqn:5}, \eqref{eqn:30}, \eqref{eqn:55}, \eqref{eqn:451}, \eqref{eqn:452}, \eqref{eqn:453}, and \eqref{eqn:30-a}
(here, we can use Lemma~\ref{lem:21}~\eqref{item:7} to verify that \eqref{eqn:41} is satisfied).

\subsection{The metric scaling limit of Aldous and Pittel's RGIV process}
We now state the result on the metric scaling limit of the RGIV process inside the critical window; recall the notation from Section~\ref{sec:def-rgiv}.

\begin{thm}\label{thm:rgiv-scaling-limit}
Fix $\lambda\in\bR$.
For $i\geq 1$, let $\cC_i\big(\rgiv_n(\lambda)\big)$ be the $i$-th largest component in $\rgiv_n(\lambda)$, and view it as a metric measure space using the counting measure. 
Then
\begin{equation}\label{eqn:64}
	\Big(
	\scl\Big( 
	\frac{a_1^{\sss\text{rgiv}}}{n^{1/3}}\ ,\
	\frac{a_2^{\sss\text{rgiv}}}{n^{2/3}} \Big) 
	\cC_i\big(\rgiv_n(\lambda)\big)\, ,
	\ i \geq 1 \Big) 
	\weakc 
	\vCrit\big(a_3^{\sss\text{rgiv}}\cdot\lambda\big)
	\ \ \text{ as }\ \ n\to\infty \, ,
\end{equation}
where
\[
a_1^{\sss\text{rgiv}}
=
\frac{2^{8/3}}{3^{2/3}\cdot\pi}\, ,\ \ \
a_2^{\sss\text{rgiv}}
=
\frac{2^{1/3}}{3^{1/3}}\, ,\ \ \ \text{ and }\ \ \
a_3^{\sss\text{rgiv}}
=
\frac{3^{2/3}}{2^{2/3}}\, .
\]
\end{thm}

\section{Proofs: Universality}
\label{sec:proofs-universality}
In this section we prove Theorem \ref{thm:gen-2}.
Sections \ref{sec:def-size-bias}--\ref{sec:proof-tilt-p-trees} collect some preliminary results related to the model.
Section \ref{sec:scaling-rank-one-connected} studies the scaling limit of a related model of a {\it connected} random graph. 
Building on this result, Section \ref{sec:proof-aldous-gen-2} completes the proof.

\subsection{Size-biased random order}\label{sec:def-size-bias}
Given a collection of positive weights $\vx=(x_i, i\in [n])$, a \emph{size-biased random order} of $[n]$ is a random permutation $(v(1), v(2), \ldots, v(n))$ with $\pr(v(1) =k) \propto x_k$ for $k\in [n]$, and for $2\leq i\leq n$ and $k\in [n]\setminus \set{v(1), \ldots, v(i-1)}$,
\begin{align}
\pr\big(v(i) =k|v(1), v(2), \ldots, v(i-1) \big) \propto x_k. \label{eqn:size-bias}
\end{align}
One construction of such an ordering is as follows:  
Generate independent exponentials $\xi_i \sim \text{Exp}(x_i)$, $i\in [n]$, and arrange them in increasing order as $\xi_{v(1)} < \xi_{v(2)} < \cdots < \xi_{v(n)}$.  
Then $(v(1), \ldots, v(n))$ would be a size-biased permutation of $[n]$ using the weight sequence $\vx$. In  \cite[Section 3.1]{aldous1997brownian} Aldous used this to construct
$\cG(\vx,q)$ simultaneously with a breadth-first exploration of the graph such that the vertices in the exploration appeared in a size-biased random order. 
We describe this construction next. 
Let $\xi_{i,j} \sim \text{Exp}(qx_j)$,
$i\neq j \in [n]$, be independent random variables. 
The exploration process initializes by selecting a vertex $v(1)$ with 
$\pr(v(1) =k) \propto x_k$ for $k\in [n]$.
The neighbors (sometimes referred to as children) of $v(1)$ are the vertices
$\set{i: \xi_{v(1),i} \leq x_{v(1)}}$. Write $c(1)$ for the number of children of
$v(1)$, labeled as $v(2), v(3), \ldots, v(c(1)+1)$ in
increasing order of the $\xi_{v(1), v(i)}$ values. Now move to $v(2)$ and obtain
the children of $v(2)$ as the vertices 
$i\notin\{ v(1),\ldots, v(c(1)+1)\}$ 
such that $ \xi_{v(2), i}\leq x_{v(2)}$. 
Label them as $v(c(1)+2), \ldots, v(c(1)+c(2)+1)$ in increasing order of their $\xi_{v(2),i}$ values. 
Proceed recursively until the component of $v(1)$ has been explored. 
Then select a new vertex amongst the unexplored vertices with probability proportional to their weights and proceed until all vertices have been explored.

\subsection{Partition into connected components}
\label{sec:part-conn-comp}
{We mainly follow \cite{SBSSXW14}.}
Recall that $(\cC_{i}, i \geq 1)$ denoted the ranked components of $\cG(\vx,q)$, and let 
$\cV_i := V(\cC_i)$. 
Thus, $(\cV_i, i\geq 1)$ is a random partition of the vertex set $[n]$. 
For $\cV\subseteq [n]$ write $\bG_{\cV}^{\con}$ {for} the space of simple connected graphs with vertex set $\cV$. 
Given $a > 0$ and a pmf $\vp = (p_v, v \in \cV)$,
define the probability distribution $\pr_{\con}(\; \cdot \; ; \vp, a, \cV)$ on $\bG_{\cV}^{\con}$ via
\begin{equation}
\label{eqn:pr-con-vp-a-cV-def}
\pr_{\con}(G; \vp, a, \cV) \propto \prod_{(u,v)\in E(G)} (1-e^{-a p_u p_v}) \prod_{(u,v)\notin E(G)} e^{-a p_u p_v}\, .
\end{equation}
For $i \geq 1$, define
\begin{equation}
\label{eqn:vp-i-a-i-def}
\vp_i =(p_{i,v}, v\in \cV_i) := \big( {x_v}/{\sum_{v \in \cV_i}x_v },\ v \in \cV_i\big)
\ \text{ and }\
a_i:= q  \cdot \big( \sum_{v \in \cV_i}x_v \big)^2.
\end{equation}

\begin{prop}[{\cite[Proposition 6.1]{SBSSXW14}}]
	\label{prop:generate-nr-given-partition}
	Let $N$ be the number of components in $\cG(\vx,q)$. Then for any $G_i \in  \bG_{\cV_i}^{\con}$, $1 \leq i \leq N$,
	\begin{equation*}
	\pr\left(\cC_i = G_i, \;\forall 1\leq i\leq N \mid  N, (\cV_i)_{1 \leq i \leq N} \right) = \prod_{1\leq i \leq  N} \pr_{\con}( G_i; \vp_i, a_i, \cV_i).
	\end{equation*}
\end{prop}

\subsection{Tilted $\vp$-trees and scaling limits of components }
\label{sec:proof-tilt-p-trees}
We recall here a result from \cite{SBSSXW14} on scaling limits of graphs with distribution $\pr_{\con} = \pr_{\con} (\; \cdot \; ; \vp, a, [m])$ under regularity assumptions on $a = a^{\sss(m)}> 0$ and the pmfs $\vp= \vp^{\sss(m)}=(p_i^{\sss(m)}, 1\leq i \leq m)$. 
We suppress dependence on $m$ to simplify notation. 
Let $\sigma(\vp) := (\sum_{i\in [m]}p_i^2)^{1/2}$, $p_{\max} := \max_{i \in [m]} p_i$, and $p_{\min} := \min_{ i \in [m]} p_i$. Suppose $\cG^{\vp}\sim \pr_{\con}$.  
Recall the notation $\tilde \ve^{\theta}$ from Section \ref{sec:cont-limit-descp}.
\begin{thm}[{\cite[Theorem 7.3]{SBSSXW14}}]
	\label{thm:sbssxw-conn-res}
	Assume that there exist $\bar \gamma, r_0, \eta_0 \in (0,\infty)$ such that
	\begin{equation}
	\label{eqn:ass-connected-1}
	\sigma(\vp) \to 0, \; \; \frac{p_{\max}}{[\sigma(\vp)]^{3/2+\eta_0}} \to 0, \;\;  \frac{[\sigma(\vp)]^{r_0}}{p_{\min}} \to 0, \;\;
	\text{ and }\, \, 	  a \sigma(\vp) \to \bar \gamma,
	\end{equation}
	as $m\to \infty$. 
	View $\cG^{\vp}$ as a metric measure space by endowing it with the graph distance and the probability measure that assigns mass $p_i$ to the vertex $i$.
	Then
	\begin{equation*}
	\sigma(\vp) \cdot \cG^{\vp} = \scl\left(\sigma(\vp), 1\right)  \cG^{\vp} \weakc \cG( 2 \tilde \ve^{\bar \gamma}, \bar \gamma \tilde \ve^{\bar \gamma}, \cP),
	\end{equation*}
	where the metric measure space $\cG( 2 \tilde \ve^{\bar \gamma}, \bar \gamma \tilde \ve^{\bar \gamma}, \cP)$ is as defined in Section \ref{sec:cont-limit-descp}.
\end{thm}
We recall from \cite{SBSSXW14} an ingredient in the proof of this theorem that will be required in this paper. 
Write $\bT_m^{\ord}$ for the space of ordered rooted trees with vertex set $[m]$. 
Here, ``ordered'' means that the children of every vertex is arranged from ``left'' to ``right.''
(In other words, these are rooted labeled trees with a plane embedding.)
Let $\vt\in \bT_m^{\ord}$ with root $\phi$.
For a vertex $v\neq\phi$, let $(\phi=v_0, v_1, \ldots, v_k=v)$ be the path in $\vt$ from $\phi$ to $v$.
Let $R(v; \vt)$ be the set of all vertices $u$ such that for some $0\leq i\leq k-1$, $u$ is a child of $v_i$ and $u$ appears on the right side of $v_{i+1}$ in the plane embedding.
Define the collection of \emph{permitted edges} to be
\[
\sP(\vt):= \set{\set{v, u}: v\in [m]\setminus\phi, u\in R(v; \vt)}\, . 
\]
For $\vt\in
\bT_m^{\ord}$ and $v\in [m]$, write $d_v(\vt)$ for the number of
children of $v$ in $\vt$. 
Define a probability measure $\pr_{\ord}(\cdot)=\pr_{\ord}(\; \cdot \; ; \vp)$ {on $\bT_m^{\ord}$ via}
\begin{equation}
\label{eqn:ordered-p-tree-def}
\pr_{\ord}(\vt) = \prod_{v\in [m]} {p_v^{d_v(\vt)}}/{(d_v(\vt)) !}, \qquad \vt \in \bT_m^{\ord}.
\end{equation}
A random tree with distribution $\pr_{\ord}$ is a $\vp$-tree (\cite{aldous2004exploration,pitman2001random}) together with a plane embedding.
Define $L : \bT_m^{\ord} \to \bR_+$ by
\begin{equation}
\label{eqn:ltpi-def}
\displaystyle L(\vt):= \Big(\prod_{\set{i,j}\in E(\vt)} \Big[\frac{\exp(a p_i p_j)- 1}{ap_ip_j} \Big]\Big) \exp\big(\sum_{\set{i,j} \in \sP(\vt)} a p_i p_j\big).
\end{equation}
Define the \emph{tilted} $\vp$-tree distribution $\tilde{\pr}_{\ord}(\cdot) = \tilde \pr_{\ord}(\; \cdot \; ; \vp, a)$ via 
\begin{equation}
\label{eqn:tilt-p-tree-dist}
\frac{d\tilde{\pr}_{\ord}}{d\pr_{\ord}}(\vt) = \frac{L(\vt)}{\E_{\ord}[L] }\, , \qquad  \vt \in \bT_m^{\ord},
\end{equation}
where $\E_{\ord}$ is the expectation operator with respect to $\pr_{\ord}$.
\begin{prop}[{\cite[Proposition 7.4]{SBSSXW14}}]\label{prop:gp-const-tilt-surp}
Fix a pmf $\vp$ on $[m]$ and $a> 0$. Then $\cG^{\vp}\sim \pr_{\con}$ can be constructed as follows: Generate $\cT^{\vp}\sim \tilde \pr_{\ord}$. Conditional on $\cT^{\vp}$, add each permitted edge $\set{u,v}\in\sP(\cT^{\vp})$ independently with probability $1-\exp(-ap_up_v)$, and then forget about the root of $\cT^{\vp}$ and the plane embedding of the tree $\cT^{\vp}$.
\end{prop}

\subsection{Connected random graphs with blobs}\label{sec:scaling-rank-one-connected}
For $m \geq 1$, consider $a = a^{\sss(m)} > 0$, and pmfs $\vp = \vp^{\sss(m)} = (p_i, i\in [m])$. 
Suppose we are given the following:

\begin{enumeratea}
\item 
{\bf Blob level superstructure:}  Let $\cG^{\vp}\sim \pr_{\con}(\; \cdot \; ;\vp,a,[m]) $.
\item 
{\bf Blobs:} Let $\vM = \vM^{\sss(m)} = \big\{(M_i, d_i, \mu_i) :  i\in [m]\big\}$ be a family of compact metric measure spaces where $\mu_i$, $i\in [m]$, are probability measures. 
Recall $u_i$ from \eqref{eqn:uik-def}.
\item 
{\bf Junction points:} $\vX = (X_{i,j},\ i,j \in [m])$ be a family of independent random variables (that is also independent of $\cG^{\vp}$) such that for each $i,j\in [m]$, $X_{i,j}\in M_i$ has law $\mu_i$.
\end{enumeratea}
Recall the operation $\Gamma$ from Section \ref{sec:inter-blob-distance}.
Define
\begin{equation}\label{eqn:g-barp-Am-def}
\bar \cG^{\vp} := \Gamma(\cG^{\vp}, \vp, \vM, \vX), \;\; 
A_m := \sum_{i \in [m]} p_i u_{i},\, \,
\text{ and }\, \,
{\diam_{\max}} := \max_{i \in [m]} \diam(M_i).
\end{equation}

\begin{thm}\label{thm:augment-metric-space-connected}
	Assume that \eqref{eqn:ass-connected-1} holds and that for some $\eta_0 \in (0,\infty)$,
	\begin{equation}
	\label{eqn:ass-connected-2}
	\lim_{m \to \infty} {[\sigma(\vp)]^{1/2-\eta_0} {\diam_{\max}}}/{(A_m+1)} = 0.
	\end{equation}
	Then
	$\frac{\sigma(\vp)}{A_m + 1} \cdot \bar \cG^{\vp} \weakc \cG(2 \tilde \ve^{\bar \gamma}, \bar \gamma \tilde \ve^{\bar \gamma}, \cP), \mbox{ as } m \to \infty.$
\end{thm}

\noindent \textbf{Proof: }  
We construct $ \bar \cG^{\vp}$ in the following way:
Consider a coupling $(\cG^{\vp}, \cT^{\vp}) \in \bG_{[m]}^{\con}\times \bT_m^{\ord}$ as in Proposition \ref{prop:gp-const-tilt-surp} with $\cG^{\vp} \sim \pr_{\con}$ and $\cT^{\vp} \sim \tilde \pr_{\ord}$.  
Generate the random junction points $\vX$ independent of $(\cG^{\vp}, \cT^{\vp})$.
Set $ \bar \cG^{\vp} := \Gamma(\cG^{\vp}, \vp, \vM, \vX)$.
Further, let $\bar \cT^{\vp} := \Gamma(\cT^{\vp}, \vp, \vM, \vX)$.
By Theorem \ref{thm:sbssxw-conn-res} we have 
$\sigma(\vp)\cdot \cG^{\vp} \weakc \cG(2 \tilde \ve^{\bar \gamma}, \bar \gamma \tilde \ve^{\bar \gamma}, \cP)$. Thus we only need to prove
\begin{equation}
\label{eqn:1481}
d_{\GHP}\Big( \sigma(\vp)\cdot \cG^{\vp}, \frac{\sigma(\vp)}{A_m + 1}\cdot \bar \cG^{\vp}  \Big) \weakc 0, \mbox{ as } m \to \infty.
\end{equation}
It would be convenient to work with slight variants of the original spaces $\cG^{\vp}$ and $\bar \cG^{\vp}$. 
Write $\cG_*^{\vp}$ (resp. $\bar \cG_*^{\vp}$) for the metric measure space obtained from $\cT^{\vp}$ (resp. $\bar \cT^{\vp}$) by identifying $i$ and $j$ (resp. $X_{i,j}$ and $X_{j,i}$) for each edge $\{i,j\} \in E(\cG^{\vp}) \setminus E(\cT^{\vp})$, instead of placing an edge of length one between them. 
Recall that $\spls(\cG^{\vp})$ denotes the number of surplus edges in $\cG^{\vp}$. Then one can check that
\begin{equation}
\label{eqn:1600}
d_{\GHP}( \cG_*^{\vp}, \cG^{\vp}) \leq \spls(\cG^{\vp}),\ \text{ and }\
d_{\GHP}( \bar \cG_*^{\vp}, \bar \cG^{\vp}) \leq \spls(\cG^{\vp}).
\end{equation}
View $\cT^{\vp}$ as a metric measure space by endowing it with the tree distance and the probability measure that assigns mass $p_i$ to the vertex $i$.
Let $\cR_m\subseteq \cT^{\vp}\times \bar \cT^{\vp}$ be given by $\cR_m := \{ (i,x): i \in [m], x \in M_i\}$.
Let $\nu_m$ be the measure on $\cT^{\vp}\times \bar \cT^{\vp}$ given by
$\nu_m( \{i\} \times A) := p_i \mu_i(A \cap M_i)$ for $i\in [m]$, $A \subseteq \bar \cT^{\vp}$.
Write 
\[
\dis(\cR_m)=\dis\Big(\cR_m; \sigma(\vp)\cT^{\vp}, \frac{\sigma(\vp)}{A_m+1}\bar \cT^{\vp}\Big)\ \ 
\text{ and }\ \ 
\dsc(\nu_m)= \dsc\Big(\nu_m; \sigma(\vp)\cT^{\vp},  \frac{\sigma(\vp)}{A_m+1}\bar \cT^{\vp}\Big).
\] 
By \cite[Lemma 4.2]{addario2013scaling}, we have
\begin{equation}
\label{eqn:1503}
d_{\GHP}\Big( \sigma(\vp)\cG_*^{\vp},  \frac{\sigma(\vp)}{A_m+1} \bar \cG_*^{\vp}\Big)
\leq
(\spls(\cG^{\vp})+1) \max\Big\{ \frac{\dis(\cR_m)}{2}, \dsc(\nu_m), \nu_m(\cR_m^c) \Big\}\, .
\end{equation}
It is easy to check that for all $m$, $\dsc(\nu_m)=0$ and $\nu_m(\cR_m^c)=0$. Therefore by \eqref{eqn:1503} and \eqref{eqn:1600}, to obtain \eqref{eqn:1481} we only need to show that
\begin{align}\label{eqn:1145-1}
(\spls(\cG^{\vp}), m\geq 1) \mbox{ is tight, and }
\dis(\cR_m) \weakc 0 \mbox{ as } m \to \infty. 
\end{align}
We will prove \eqref{eqn:1145-1} in Lemma \ref{lem:1542} and Lemma \ref{lem:1567} below, which will complete the proof of Theorem \ref{thm:augment-metric-space-connected}. 
\qed

\begin{lem}
	\label{lem:1542}
	The sequence $\big(\spls(\cG^{\vp});\, m\geq 1\big)$ is tight.
\end{lem}
\noindent \textbf{Proof:}
By \cite[Corollary 7.13]{SBSSXW14}, there exists a constant $K_1 > 0$ such that for $m$ large,
\begin{equation}
\label{eqn:1725}
\E_{\ord}[L^2] \leq K_1.
\end{equation}
Since $\cG^{\vp}$ is obtained by adding each permitted edge $\set{u,v}$ of $\cT^{\vp}$ independently with probability $1-\exp(-ap_u p_v)$,
$\E [ \spls (\cG^{\vp})] \leq \E\big[ \sum_{\{i,j\} \in \sP(\cT^{\vp})} a p_i p_j\big]$.
Using \eqref{eqn:tilt-p-tree-dist} and \eqref{eqn:1725} together with the relations $L(\vt) \geq 1$ and
$ \sum_{\{i,j\} \in \sP(\vt)} a p_i p_j \leq L(\vt)$,
we have, for all large $m$,
$	\E [ \spls (\cG^{\vp})] \leq \E[L(\cT^{\vp})] =\E_{\ord}[L^2]/\E_{\ord}[L]\leq K_1.$
Thus $\sup_{m \geq 1}\E[\spls(\cG^{\vp})] < \infty$, which implies the desired tightness. \qed

\begin{lem}\label{lem:1567}
As $m \to \infty$, $\dis(\cR_m) \weakc 0$.
\end{lem}
\noindent \textbf{Proof:} 
Let $L(\cdot)$ be as in \eqref{eqn:ltpi-def}. 
Suppose $(\Omega,\cF,\bP)$ is a probability space where $\cT^{\vp}$ and $\bar \cT^{\vp}$ are defined.
Let $\bP'$ be the probability measure on $\Omega$ given by
\begin{equation}
\label{eqn:1537}
\frac{d \pr}{ d \pr'}(\omega) = \frac{L(\cT^{\vp}(\omega))}{\E_{\ord}[L]} \mbox{ for } \omega \in \Omega.
\end{equation}
By \eqref{eqn:tilt-p-tree-dist}, $\cT^{\vp} \sim \pr_{\ord}$ under $\pr'$.
Suppose we show that for any $\eps>0$,
\begin{align}\label{eqn:9}
\pr'(\dis(\cR_m) > \eps)\to 0\ \ \text{ as }\ \ m\to\infty\, .
\end{align}
Then 
$
\big(\pr(\dis(\cR_m) > \eps)\big)^2
\leq
\pr'(\dis(\cR_m) > \eps)\cdot\bE_{\ord}[L^2]
\to 0
$
as $m\to\infty$, 
where the first step uses \eqref{eqn:1537}, the fact that $L(\vt)\geq 1$, and the Cauchy-Schwarz inequality, and the second step uses \eqref{eqn:1725} and \eqref{eqn:9}.
So we just have to prove \eqref{eqn:9}.

We introduce some additional random variables: Let $J$ and $Y_1,\ldots, Y_m$ be respectively $[m]$ and $M_1,\ldots, M_m$ valued independent random variables that are also independent of $\cT^{\vp}$ such that $J \sim \vp$, and $Y_i \sim \mu_i$ for $i \in [m]$. We extend the probability space to incorporate these random variables and still write $\pr'$ for the underlying probability measure.

For $ x \in \bar \cT^{\vp}$, write $i(x)$ for the $i \in [m]$ with $x \in M_{i}$. Write $d_{\cT}$ and $d_{\bar \cT}$ for the metrics on $\cT^{\vp}$ and $ \bar \cT^{\vp}$ respectively.   
Then
\begin{equation}
\label{eqn:1157}
\dis(\cR_m) = \sup_{x,y \in \bar \cT^{\vp}}\Big\{ \Big|\sigma(\vp) d_{\cT}(i(x),i(y)) -\frac{\sigma(\vp)}{A_m+1} d_{\bar \cT}(x,y)\Big|\Big\}.
\end{equation}
We first show for any fixed $i_0 \in [m]$ and $x_0 \in M_{i_0}$,
\begin{equation}\label{eqn:1142}
\dis(\cR_m) 
\leq 
4 \sup_{y \in \bar \cT^{\vp}}
\Big\{ 
\Big|\sigma(\vp) d_{\cT}(i_{0},i(y)) - \frac{\sigma(\vp)}{A_m+1} d_{\bar \cT}(x_0,y)\Big|
\Big\} 
+ 
\frac{2\sigma(\vp)}{A_m+1} {\diam_{\max}} \, .
\end{equation}
Indeed, for any two points $x,y \in \bar \cT^{\vp}$, there are two unique paths between $i_0$ and $i(x)$ and between $i_0$ and $i(y)$ in $\cT^{\vp}$. Write $(i_0, \ldots , i_{k})$ for the longest common path shared by these two paths. 
Let $i_* = i_{k}$ and $x_* = X_{i_k, i_{k-1}}$. Since $\cT^{\vp}$ is a tree,
\begin{equation}
\label{eqn:1146}
d_{\cT}(i(x), i(y)) = d_{\cT}(i_0, i(x)) + d_{\cT}(i_0, i(y)) - 2d_{\cT}(i_0, i_*).
\end{equation}
By a similar observation for $\bar \cT^{\vp}$, we have
\begin{equation}
\label{eqn:1151}
d_{\bar \cT}(x , y) \leq d_{\bar \cT}(x_0, x) + d_{\bar \cT}(x_0, y) - 2d_{\bar \cT}(x_0, x_*) \le	d_{\bar \cT}(x , y) + 2 {\diam_{\max}} .
\end{equation}
Hence, \eqref{eqn:1142} follows using \eqref{eqn:1146} and \eqref{eqn:1151} in \eqref{eqn:1157}.

Next, replace $(i_0, x_0)$ by $(I, Y_I)$, where $I\in [m]$ is the root of $\cT^{\vp}$, and replace every $y \in M_j$ in \eqref{eqn:1142} with $Y_j \in M_j$, which incurs an error of at most ${4\sigma(\vp) {\diam_{\max}}}/{(A_m+1)}$. Thus,
\begin{equation*}
\dis(\cR_m) \leq 4 \sup_{j \in [m]} \Big\{ \Big|\sigma(\vp) d_{\cT}(I,j) - \frac{\sigma(\vp)}{A_m+1} d_{\bar \cT}(Y_I,Y_j)\Big|\Big\} + \frac{6\sigma(\vp) {\diam_{\max}}}{A_m+1}.
\end{equation*}
Using \eqref{eqn:ass-connected-2} \ch{and $\sigma(\vp)\to 0$} we can find $m_0$ such that ${6\sigma(\vp) {\diam_{\max}}}/{(A_m+1)} < \eps/5$ for $m > m_0$. Then for $m>m_0$,
\begin{align}
\pr'( \dis(\cR_m) > \eps )
\leq& \pr' \Big( \sup_{j \in [m]} { \Big|\sigma(\vp) d_{\cT}(I,j) - \frac{\sigma(\vp)}{A_m+1} d_{\bar \cT}(Y_I,Y_j)\Big|} > \frac{\eps}{5} \Big) \nonumber \\
\leq& \frac{1}{p_{\min}}\sum_{j \in [m]} p_j \pr'  \Big( \Big|\sigma(\vp) d_{\cT}(I,j) - \frac{\sigma(\vp)}{A_m+1} d_{\bar \cT}(Y_I,Y_j)\Big| > \frac{\eps}{5} \Big) \nonumber \\
=& \frac{1}{p_{\min}}\pr' \Big( \Big|\sigma(\vp) d_{\cT}(I,J) - \frac{\sigma(\vp)}{A_m+1} d_{\bar \cT}(Y_I,Y_J)\Big| > \frac{\eps}{5} \Big). \label{eqn:1173}
\end{align}
Write $R^* := d_{\cT}(I, J)+1$,  and let  $(I_0=I, I_1, ..., I_{R^*-1} = J)$ be the path between $I$ and $J$ in $\cT^{\vp}$. 
Define $\xi_0^* := d_I(Y_I, X_{I,I_1})$, $\xi_{R^*-1}^* := d_J(X_{J,I_{R^*-2}}, Y_J)$, and 
$\xi_i^* := d_{I_i}(X_{I_i,I_{i-1}}, X_{I_i, I_{i+1}})$ 
for $1\leq i \leq R^*-2$.
In this notation $d_{\bar \cT}(Y_I,Y_J) = \sum_{i=0}^{R^*-1} \xi_i^* + (R^* - 1)$ and $d_{\cT}(I,J) = R^*-1$. Thus,
\ch{\begin{equation}
	\label{eqn:1590}
	\sigma(\vp)d_{\cT}(I,J) - \frac{\sigma(\vp)}{A_m+1} d_{\bar \cT}(Y_I,Y_J)  = \frac{\sigma(\vp)}{A_m+1} \Big( \Big[\sum_{i=0}^{R^*-1} (A_m - \xi_i^*)\Big] - A_m \Big).
	\end{equation}}
The summation above involves the path in a $\vp$-tree connecting the root and a random vertex selected according to the distribution $\vp$. This admits the following alternate construction: 
Let $\vJ:=(J_i, i\geq 0)$ be a sequence of iid $[m]$-valued random variables with law $\vp$.  For each $j \in [m]$, let $\mvxi^{\sss(j)}:=(\xi_i^{\sss(j)}, i\geq 0)$ be an \chh{iid} sequence with $\xi_1^{\sss(j)}\equald d_j(X_{j,1}, X_{j,2})$ such that $\mvxi^{\sss(1)},\ldots, \mvxi^{\sss(m)}$ and $\vJ$ are jointly independent. Let $(\Omega'', \cF'', \pr'')$ be the probability space on which the random variables $\vJ$ and $\big(\mvxi^{\sss(j)}; j\in [m] \big)$ are defined. 
Let $R$ denote the first repeat time of the sequence $\vJ$, i.e.,
$R:= \inf\, \big\{ k \geq 1: {J_k=J_i \text{ for some } i<k} \big\}$.
By \cite[Corollary 3]{camarri2000limit}, 
\begin{equation*}
(I_0, ..., I_{R^*-1}; R^*)_{\pr'} \stackrel{d}{=} (J_0, ..., J_{R-1}; R)_{\pr''}.
\end{equation*}
Owing to the independence structure,  \ch{we further} have
\begin{equation*}
(I_0, ..., I_{R^*-1}; R^*; \xi_0^*, ...,\xi_{R^*-1}^*)_{\pr'} \stackrel{d}{=} (J_0, ..., J_{R-1}; R; \xi_0^{\sss(J_0)}, ..., \xi_{R-1}^{\sss(J_{R-1})} )_{\pr''}.
\end{equation*}
For $i\geq 0$ write $\Delta_i := \xi_i^{\sss(J_i)} - A_m$. Using \eqref{eqn:1590} and the fact that $\sigma(\vp) \to 0$ \chsen{as $m\to\infty$, we have}
\begin{equation}
\label{eqn:1223}
\pr'\Big(\Big |\sigma(\vp) d_{\cT}(I,J) - \frac{\sigma(\vp)}{A_m+1} d_{\bar \cT}(Y_I,Y_J)\Big| > \frac{\eps}{5} \Big) \leq \pr''\Big( \frac{\sigma(\vp)}{A_m+1} \Big| \sum_{i=0}^{R-1} \Delta_i \Big| > \frac{\eps}{6} \Big)
\end{equation}
\chsen{for all large $m$.}
Note that $\big(\xi_i^{\sss (J_i)}; i\geq {0}\big)$ is a collection of iid random variables with mean
$\sum_{k\in [m]} p_k \E[\xi_0^{\sss(k)}] = A_m$.
Thus $(\sum_{i=0}^k \Delta_i,\ k\geq 0)$ is a martingale with respect to \ch{its} natural filtration.
Then, for any $t > 0$,
\begin{equation}
\label{eqn:1199}
\pr''\Big(\frac{\sigma(\vp)}{A_m+1} \Big| \sum_{i=0}^{R-1} \Delta_i \Big| > \frac{\eps}{6}\Big) 
\leq 
\pr''(R \geq t) + \pr''\Big(\sup_{0\leq k \leq t-1}  \Big|  \sum_{i=0}^{k} \Delta_i \Big| > \frac{\eps (A_m+1)}{6 \sigma(\vp)}\Big).
\end{equation}
For the first term, by \ch{\cite[Lemma 10.7]{SBSSXW14}} (see also 
\cite[Displays (26) and (29)]{camarri2000limit}),
\begin{equation}
\label{eqn:1794}
\pr''\left(R \geq t\right)  \leq 2 \exp \big( -t^2 \sigma^2(\vp)/24\big)\ \mbox{ for }\ t \in (0, 1/p_{\max}).
\end{equation}
The second term on the right side of \eqref{eqn:1199} can be bounded by using \chh{Markov's} inequality and the Burkholder-Davis-Gundy (BDG) inequality \cite{burkholder-davis-gundy}.
For fixed $r \geq 1$, we have
\begin{align}
\pr''\Big(\sup_{0\leq k \leq t-1}  \big|  \sum_{i=0}^{k} \Delta_i \big| > \frac{\eps (A_m+1)}{6 \sigma(\vp)}\Big)
\leq
\Big(\frac{6 \sigma(\vp)}{\eps (A_m+1)}\Big)^{2r}  \E''\Big[ \sup_{0\leq k \leq t-1}  \big|  \sum_{i=0}^{k} \Delta_i \big|^{2r} \Big]& \nonumber\\
\hskip10pt
\leq \Big(\frac{6 \sigma(\vp)}{\eps (A_m+1)}\Big)^{2r} K_2(r) \E''\Big[ \big(\sum_{i=0}^{t-1} \Delta_i^2\big)^{r} \Big]
\leq  \Big(\frac{6 \sigma(\vp)}{\eps (A_m+1)}\Big)^{2r} K_2(r)\cdot t^r {\diam_{\max}^{2r}} \, , \label{eqn:1240}&
\end{align}
where $\bE''$ denotes expectation with respect to $\bP''$, the constant $K_2(r)>0$ comes from the BDG inequality, and the last step uses $|\Delta_i| \leq {\diam_{\max}}$. 
Define 
\chsen{$\alpha_m := \sqrt{-24 r \log \sigma(\vp)}$}
and $t_m := \alpha_m/\sigma(\vp)$. 
Using \eqref{eqn:1199}, \eqref{eqn:1794}, and \eqref{eqn:1240} with $t = t_m $, we have for sufficiently large $m$,
\begin{align}
\pr''\Big( \frac{\sigma(\vp)}{A_m+1} \big| \sum_{i=0}^{R-1} \Delta_i \big| > \frac{\eps}{6} \Big)
\leq& 2 \exp \Big( -\frac{t_m^2 \sigma^2(\vp)}{24}\Big) + K_3(r,\eps)\cdot \frac{\alpha_m^r \sigma^r(\vp) {\diam_{\max}}^{2r}}{(A_m+1)^{2r}} \nonumber \\
\leq& 2 [\sigma(\vp)]^r + K_3(r,\eps) \alpha_m^r [\sigma(\vp)]^{2r\eta_0}\, , \label{eqn:1244}
\end{align}
where $K_3(r,\eps) = 6^{2r} K_2(r)/\eps^{2r}$, and the last line uses  \eqref{eqn:ass-connected-2}. 
(It is easy to check that $t_m p_{\max}\to 0$ as $m \to\infty$. 
Thus, we can apply \eqref{eqn:1794} for $m$ large, and \eqref{eqn:1244} is valid for large $m$.) 
By \eqref{eqn:1244}, \eqref{eqn:1223}, and \eqref{eqn:1173}, taking 
$r =r_{\ast}:= \max\big\{r_0 \, , \lceil (r_0+1)/(2\eta_0) \rceil \big\}$
\chsen{(so that 
$\sigma(\vp)^{r_{\ast}} \leq \sigma(\vp)^{r_0}$ and 
$\sigma(\vp)^{2 \eta_0 r_{\ast}} \leq \sigma(\vp)^{r_0 + 1}$)}, 
we get
	\begin{equation*}
	\pr' \big( \dis(\cR_m) > \eps \big)
	\leq 
	\Big(2 [\sigma(\vp)]^{r_0} + K_3(r_{\ast},\eps) \alpha_m^{r_{\ast}} [\sigma(\vp)]^{r_0+1}\Big)/p_{\min} \, .
	\end{equation*}
By \eqref{eqn:ass-connected-1}, the above expression goes to zero as $m \to \infty$.
This completes the proof of \eqref{eqn:9}. 
\qed

\subsection{Proof of Theorem \ref{thm:gen-2}}
\label{sec:proof-aldous-gen-2}
The previous section deals with asymptotics for a single connected component. This section leverages the above results to study maximal components of $\bar\cG(\vx,q,\vM)$ and prove Theorem \ref{thm:gen-2}.
Recall that $\vx$ is the weight sequence, $\vM$ is the collection of blobs, $\vX$ is the collection of random junction points, and $u_1,\ldots, u_n$ are the average distances.
Also recall that $\sigma_r = \sum_{i \in [n]} x_i^r$, $r \geq  1$, and $\tau = \sum_{i \in [n]} x_i^2 u_i$.  We start with the following auxiliary result.

\begin{prop}\label{prop:moments-convergence-each-component}
Under Conditions \ref{ass:aldous-basic-assumption} and \ref{ass:gen-2}, 
for each $i\geq 1$,
\begin{equation*}
\frac{\sum_{v\in \cC_i} x_v^2}{\sum_{v\in \cC_i} x_v } \cdot \frac{\sigma_2}{\sigma_3} 
\weakc 1,
\ \ \text{ and }\ \ 
\frac{\sum_{v\in \cC_i}  x_v  u_{v}}{\sum_{v\in \cC_i}  x_v} \cdot \frac{\sigma_2}{\tau}  
\weakc 1\, , 
\ \ \text{ as }\ \  n \to \infty.
\end{equation*}
\end{prop}
We will make use of the following lemma in the proof.
\begin{lem}[\cite{SBSSXW14}, Lemma 8.2]\label{lem:size-biased-partial-sum}
	For $n\geq 1$, fix $\mvx := (x_i > 0,\ i \in [n])$ and $\mva  := (a_i \geq 0,\ i \in [n])$ such that $c_n:= \sum_{i \in [n]} x_i a_i/{\sum_{i\in [n]} x_i} > 0$.  Let $(v(i), i \in [n])$ be a size-biased random \ch{ordering} of $[n]$ using the weights $\mvx$.
	Let $x_{\max}:= \max_{i \in [n]} x_i$ and $a_{\max} := \max_{i \in [n]}a_i$. 
	Let $\ell= \ell(n)  \in [n]$ be such that
	\begin{equation}
	\label{eqn:ass-size-biased-partial-sum}
	{\ell x_{\max}}/{\sigma_1} \to 0\, ,\ \text{ and }\ \   a_{\max}/( \ell c_n) \to 0 \ \ \text{ as }\ \ n\to\infty.
	\end{equation}
	Then
	\begin{equation*}
	\sup_{k \leq \ell }\Big| \frac{\sum_{i = 1}^k a_{v(i)}}{\ell c_n} - \frac{k}{\ell} \Big| \weakc 0\ \ 
	\text{ as }\ \ n\to\infty.
	\end{equation*}
\end{lem}
We are now ready to prove Proposition \ref{prop:moments-convergence-each-component}.

\noindent\textbf{Proof of Proposition \ref{prop:moments-convergence-each-component}: }
We only work with $i=1$ to keep the notation simple.
Recall the breadth-first construction of $\cG(\vx,q)$ from Section \ref{sec:def-size-bias}, where vertices are explored in a size-biased order $(v(i), 1\leq i\leq n)$ using the weight sequence $\vx$. 
The following properties of this exploration were derived in \cite{aldous1997brownian}:
\begin{inparaenuma}
	\item There exist random variables $m_{L}, m_{R} \in [n]$ such that $V(\cC_1)=\set{ v(i) : m_L +1 \leq i \leq m_R} $.
	\item Under \chh{Condition} \ref{ass:aldous-basic-assumption}, $\sum_{i=1}^{m_{R}} x_{v(i)}$ is tight.
	\item By Theorem~\ref{thm:aldous-review},  as $n \to \infty$, $\sum_{i=m_L+1}^{m_R} x_{v(i)} \weakc \chh{\gamma_1(\lambda)}$.
\end{inparaenuma}
In this notation, the claim in Proposition~\ref{prop:moments-convergence-each-component}, for $i=1$, is equivalent to
\begin{equation}
\label{eqn:prop67-equiv}
\frac{\sum_{i=m_L+1}^{m_R} x_{v(i)}^2}{\sum_{i=m_L+1}^{m_R} x_{v(i)}}\cdot \frac{\sigma_2}{\sigma_3} \weakc 1, 
\ \ \text{ and }\ \
\frac{\sum_{i=m_L+1}^{m_R} x_{v(i)}  u_{v(i)}}{\sum_{i=m_L+1}^{m_R} x_{v(i)}}\cdot \frac{\sigma_2}{\tau} \weakc 1.
\end{equation}
Thus, it is enough to prove \eqref{eqn:prop67-equiv}.
We first show that
\begin{equation} \label{eqn:1429}
\frac{\sum_{i =m_{L}+1}^{m_R}  x_{v(i)} }{m_R - m_L} \cdot \frac{\sigma_1}{\sigma_2} \weakc 1, \mbox{ as } n \to \infty.
\end{equation}
Fix $\eta >0$. Since $\big(\sum_{i=1}^{m_{R}} x_{v(i)};\, n\geq 1\big)$ is a tight sequence, there exists $T >0$ such that for all $n\geq 1 $,
$	\pr\left(\sum_{i=1}^{m_{R}} x_{v(i)} \geq T\right) < \eta$.
Let $m_0 := \sigma_1/\sigma_2$ and apply Lemma \ref{lem:size-biased-partial-sum} with $\ell = 2T m_0$ and $a_i = x_i$ (thus  $c_n = \sigma_2/\sigma_1$). Condition~\eqref{eqn:ass-size-biased-partial-sum} is equivalent to $ x_{\max}/\sigma_2 \to 0$, which follows from \chsen{Condition}~\ref{ass:aldous-basic-assumption}. By Lemma \ref{lem:size-biased-partial-sum}, there exists $N_\eta > 0$ such that for all $n > N_\eta$,
\begin{equation}\label{eqn:1438}
\pr\Big( \sup_{k \leq 2Tm_0}\big|\sum_{i=1}^k x_{v(i)} - \frac{k}{m_0}\big| > \eta\Big) < \eta.
\end{equation}
On the event
$\big\{\sup_{k \leq 2Tm_0}\left|\sum_{i=1}^k x_{v(i)} - \frac{k}{m_0}\right| \leq \eta\big\} 
\cap 
\big\{\sum_{i=1}^{m_{R}} x_{v(i)} \leq T\big\}$,
we must have
$m_L < m_R < 2T m_0$ (assuming $\eta < T$), and hence $|\sum_{i=m_L+1}^{m_R} x_{v(i)} - {(m_R - m_L)}/{m_0}| < 2 \eta$. 
Since $\eta$ can be arbitrarily small, using \eqref{eqn:1438} we have
\begin{equation}\label{eqn:47}
\Big|\sum_{i=m_L+1}^{m_R} x_{v(i)} - \frac{m_R - m_L}{m_0}  \Big|\weakc 0 \mbox{ as } n \to \infty.
\end{equation}
By property (c) of the exploration and \eqref{eqn:47},  
${(m_R - m_L)}/{m_0} \weakc \gamma_1(\lambda)$,
which coupled with \eqref{eqn:47} yields \eqref{eqn:1429}.

Now, under Conditions \ref{ass:aldous-basic-assumption} and \ref{ass:gen-2},
\begin{equation}
\label{eqn:we-want-moment}
\frac{\sigma_2 x_{\max}^2}{\sigma_3} \to 0\, ,\ \ 
\text{ and }\ \  \
\frac{\sigma_2 x_{\max} {\diam_{\max}}}{\tau} \to 0\, ,\ \  \mbox{ as }\ n \to \infty.
\end{equation}
Thus, repeating the above argument with respectively $ x_{v(i)}^2$  and $ x_{v(i)} u_{v(i)}$ in place of $x_{v(i)}$ and using Lemma \ref{lem:size-biased-partial-sum}, we get
\begin{equation}\label{eqn:1447}
\frac{\sum_{i =m_{L+1}}^{m_R}  x_{v(i)}^2 }{m_R - m_L} \cdot \frac{\sigma_1}{\sigma_3} \weakc 1\, ,\ \ 
\text{ and }\ \  \
\frac{\sum_{i =m_{L+1}}^{m_R}  x_{v(i)}  u_{v(i)} }{m_R - m_L} \cdot \frac{\sigma_1}{\tau} \weakc 1.
\end{equation}
Combining \eqref{eqn:1429} and \eqref{eqn:1447} completes the proof of \eqref{eqn:prop67-equiv}, and thus of Proposition \ref{prop:moments-convergence-each-component}. \qed\\

\noindent \textbf{Proof of Theorem \ref{thm:gen-2}:}
Since we work with the product GHP topology, it is enough to prove the assertion for the maximal component $\bar \cC_1$, i.e.,
\begin{equation}
\label{eqn:conv-first-component}
\scl\Big(\frac{\sigma_2^2}{\sigma_2 + \tau}, 1 \Big)\bar \cC_1 \weakc \cG(2 \tilde \ve_{\gamma_1}, \tilde \ve_{\gamma_1}, \cP_1 ), \mbox{ as } n \to \infty.
\end{equation}
In \eqref{eqn:conv-first-component} and the proof below, to simplify notation, we suppress dependence on $\lambda$ and write $\gamma_1$ for $\gamma_1(\lambda)$. 
By Theorem \ref{thm:aldous-review} and Proposition \ref{prop:moments-convergence-each-component},  without loss of generality, we \ch{may} consider \chh{a} probability space on which the following \chh{convergences} hold almost surely: 
\begin{equation}
\mass(\cC_1) \convas \gamma_1,   \ \
\frac{\sum_{v\in \cC_1} x_v^2}{\mass(\cC_1)} \cdot \frac{\sigma_2}{\sigma_3} \convas 1, \ \ 
\text{ and }\ \ 
\frac{\sum_{v\in \cC_1}  x_v  u_{v}}{\mass(\cC_1)} \cdot \frac{\sigma_2}{\tau}  \convas 1. \label{eqn:1616-1109}
\end{equation}
Consider the following construction of $\bar \cC_1$: 
Define $\vp := \big( {x_v}/{\mass(\cC_1)}, v \in \cC_1\big)$ and 
$a := q \cdot [\mass(\cC_1)]^2$.
Conditioned on $\cV_1 := V(\cC_1)$, let $\cC_1'$ be a $\bG^{\con}_{\cV_1}$-valued random variable with distribution {$\bP_{\con}(\cdot; \vp, a, \cV_1)$}.
By Proposition \ref{prop:generate-nr-given-partition}, 
\begin{align}\label{eqn:94}
\bar \cC_1 {\equald} \scl(1, \mass(\cC_1)) \Gamma(\cC_1', \vp, \vM, \vX)\, .
\end{align} 
In order to apply Theorem \ref{thm:augment-metric-space-connected}, we will verify the statements in \eqref{eqn:ass-connected-1} and \eqref{eqn:ass-connected-2} with $\bar \gamma = \gamma_1^{3/2}$. Note that $p_{\max} \leq {x_{\max}}/{\mass(\cC_1)}$, $p_{\min} \geq {x_{\min}}/{\mass(\cC_1)}$, and by \chh{\eqref{eqn:1616-1109}},
\begin{equation}
\label{eqn:1866}
\sigma(\vp) = \frac{\sqrt{\sum_{v \in \cC_1} x_v^2}}{\mass(\cC_1)} \sim \frac{\sigma_2}{\gamma_1^{1/2}}
\, ,
\ \ \text{ and }\ \
A_m = \frac{\sum_{v\in \cC_1}  x_v  u_{v}}{\mass(\cC_1)} \sim \frac{\tau}{\sigma_2}\, ,
\end{equation}
where the simplification uses the relation $\sigma_3 \sim \sigma_2^3$ (by Condition \ref{ass:aldous-basic-assumption}). 
In addition, {using} $q \sim 1/\sigma_2$ from Condition \ref{ass:aldous-basic-assumption}, $a \cdot \sigma(\vp) \sim q\gamma_1^2 \cdot \frac{\sigma_2}{\gamma_1^{1/2}} \sim {\gamma_1^{3/2}}$. The remaining conditions in \eqref{eqn:ass-connected-1} and \eqref{eqn:ass-connected-2} follow directly from Conditions \ref{ass:aldous-basic-assumption} and \ref{ass:gen-2}, and the details are omitted.

Applying Theorem~\ref{thm:augment-metric-space-connected} and using \eqref{eqn:94}, we see that
\begin{equation*}
\scl\Big( \frac{\sigma_2^2}{(\sigma_2 + \tau)\gamma_1^{1/2}}, \frac{1}{\gamma_1} \Big) \bar \cC_1
\weakc \cG(2 \tilde \ve^{\gamma_1^{3/2}}, \gamma_1^{3/2} \tilde \ve^{\gamma_1^{3/2}}, \cP_1).
\end{equation*}
Note that for any excursions $h$ and $g$ and Poisson point process $\cP$, for $\alpha, \beta >0$, we have
$	\scl(\alpha, \beta) \cG(h,g,\cP) \stackrel{d}{=} \cG\big( \alpha h( \cdot/\beta), \frac{1}{\beta}g(\cdot/\beta), \cP\big)$.
Thus,
\begin{equation*}
\scl(\gamma_1^{1/2}, \gamma_1) \scl\Big( \frac{\sigma_2^2}{(\sigma_2 + \tau)\gamma_1^{1/2}}, \frac{1}{\gamma_1} \Big) \bar \cC_1  \weakc \cG\Big(2 \gamma_1^{1/2}\tilde \ve^{\gamma_1^{3/2}}(\cdot/ \gamma_1), \gamma_1^{1/2} \tilde \ve^{\gamma_1^{3/2}}(\cdot/\gamma_1), \cP_1\Big).
\end{equation*}
By \cite[Display (4.9)]{SBSSXW14}, we have $\gamma_1^{1/2} \tilde \ve^{\gamma_1^{3/2}}(\cdot/\gamma_1) \stackrel{d}{=} \tilde \ve_{\gamma_1}(\cdot)$. 
Combining this observation with the last display completes the proof of \eqref{eqn:conv-first-component}. 
\qed


\section{Proofs: Associated random graph models}\label{sec:graphon-proofs}
In this section, we will prove the results stated in Section~\ref{sec:res-models}.
We start in Section~\ref{sec:g1} with the proof of Lemma~\ref{lem:1}.
In Section~\ref{sec:g2} we state a result (Proposition~\ref{prop:graphon-intermediate}) that claims that the conclusions in Theorems~\ref{thm:scaling-limit-graphon} and \ref{thm:component-limit-graphon} hold under a stronger assumption.
Here, we also introduce the notation to be used in the sequel.
From Section~\ref{sec:g3} through Section~\ref{sec:g9}, we prove various estimates working towards the proof of Proposition~\ref{prop:graphon-intermediate}.
Section~\ref{sec:g11} collects all these results to complete the proof of Proposition~\ref{prop:graphon-intermediate} via an application of Theorem~\ref{thm:gen-2}.
We then use  Proposition~\ref{prop:graphon-intermediate} in Section~\ref{sec:g12} to prove Theorems~\ref{thm:scaling-limit-graphon} and \ref{thm:component-limit-graphon}.
Finally, the proofs of Theorems~\ref{thm:graphon-Ui-Uj}, \ref{thm:l-2-approx}, and \ref{thm:rgiv-scaling-limit} are given in Sections~\ref{sec:g13} and \ref{sec:g14}.

\subsection{Proof of Lemma~\ref{lem:1}}\label{sec:g1}
We start with an elementary result.
\begin{lem}\label{lem:2}
Suppose $K\in L^3([0, 1]^2)$. 
Then 
$
\|T_K\|_{p,\, q}
\leq
\big(\int_{[0, 1]^2}|K|^3\big)^{1/3}
$
for any $p, q\in[3/2,\, 3]$.
\end{lem}
\noindent{\bf Proof:}
For any $f\in L^{3/2}([0, 1])$,
\begin{align*}
|T_Kf(x)|^3
\leq
\Big(\medint\int_0^1 |K(x, y)f(y)|\, dy\Big)^3
\leq
\Big(\medint\int_0^1 |K(x, y)|^3\, dy\Big)
\Big(\medint\int_0^1 |f|^{3/2}\Big)^2\, ,
\end{align*}
where the last step uses H\"{o}lder's inequality.
Hence,
$\|T_K f\|_3\leq \big(\int_{[0, 1]^2}|K|^3\big)^{1/3}\|f\|_{3/2}\, $, and the claim follows for $p=3/2$ and $q=3$.
The proof is completed upon noting that $\|T_K\|_{p,\, q}\leq \|T_K\|_{3/2,\, 3}$ for any $p, q\in[3/2,\, 3]$.
\qed

\medskip

\noindent{\bf Proof of Lemma~\ref{lem:1}:}
Let us first show that $T_W$ is strongly band irreducible, i.e., there do not exist Borel $A_1, A_2\subseteq [0, 1]$ with $\mu_{\circ}(A_1)>0$ and $\mu_{\circ}(A_2)>0$ such that $W(x, y)=0$ for 
$\mu_{\sss \square}$~a.e. $(x, y)\in A_1\times A_2$.
Assume on the contrary that such $A_1$ and $A_2$ exist.
Then
\begin{align}\label{eqn:89}
n^{1/3}
\langle \ind_{A_2}, T_{W_n}\ind_{A_1}\rangle
=
n^{1/3}
\langle \ind_{A_2}, T_{W_n}\ind_{A_1}\rangle
-
n^{1/3}
\langle \ind_{A_2}, T_W\ind_{A_1}\rangle
\to
\langle \ind_{A_2}, T_H\ind_{A_1}\rangle
\end{align}
as $n\to\infty$, where the last step follows from \eqref{eqn:21}.
By Condition~\ref{ass:graphon-3}, 
$\mu_{\sss \square}\big(\{(x, y) : W_n(x, y)\leq n^{-\delta_0}\}\big)=o(1)$.
Hence, 
\[
n^{1/3}
\langle \ind_{A_2}, T_{W_n}\ind_{A_1}\rangle
\geq n^{1/3}\cdot n^{-\delta_0}\cdot \mu_{\circ}(A_1)\mu_{\circ}(A_2)\cdot \big(1-o(1)\big)\, ,
\] 
which contradicts \eqref{eqn:89} as $\delta_0<1/3$.
Thus, $T_W$ is strongly band irreducible.
In particular, this implies that  $T_W$ is band irreducible, i.e., there does not exist a Borel $A\subseteq [0, 1]$ with $\mu_{\circ}(A)\in (0, 1)$ such that $W(x, y)=0$ for $\mu_{\sss \square}$~a.e. $(x, y)\in A\times A^c$.

Now, $W\in L^3([0,1]^2)\subseteq L^2([0, 1]^2)$ implies that $T_W$ is a compact operator on $L^2([0, 1])$.
Since $W\geq 0$ and $T_W$ is band irreducible, by \cite[Theorem~43.8]{zaanen-book-operator-theory},
the spectral radius of $T_W$ (which equals $+1$ under Condition~\ref{ass:graphon-1}) is an eigenvalue of $T_W$.
Further, the corresponding eigenspace is one dimensional, and is spanned by an eigenfunction $\psi$ satisfying $\psi(x)>0$ for $\mu_{\circ}$~a.e. $x$.
We can further assume that $\psi$ is normalized so that $\|\psi\|_2=1$.
Finally, 
$\|\psi\|_3
=
\|T_W \psi\|_3
\leq 
\|T_W\|_{2, 3}
<\infty$,
where the last step follows from Lemma~\ref{lem:2}.
This proves Lemma~\ref{lem:1}~\eqref{item:21} and \eqref{item:22}.

The claim in Lemma~\ref{lem:1}~\eqref{item:23} follows from \cite[Theorem~44.10]{zaanen-book-operator-theory}.
In our rather simple setting, this result follows from an easy argument, which we include here.
Note that $T_W^2=T_K$, where $K(x, y)=\int_{z=0}^1 W(x, z)W(z, y) dz$.
If $T_W^2$ were not strongly band irreducible, then we could find Borel $A_1, A_2\subseteq [0, 1]$ with positive measures such that
\[
0 
= 
\medint\int_{x\in A_1} \medint\int_{y\in A_2} K(x, y) dx dy
=
\medint\int_{z=0}^1 
\bigg( \medint\int_{x\in A_1} W(x, z) dx\bigg) 
\bigg(\medint\int_{y\in A_2} W(y, z) dy\bigg) dz\, .
\]
Hence, we can find a Borel $A_3\subseteq [0, 1]$ such that 
$W(x, z)=0$ for $\mu_{\sss \square}$~a.e. $(x, z)\in A_1\times A_3$, and 
$W(y, z)=0$ for $\mu_{\sss \square}$~a.e. $(y, z)\in A_1\times A_3^c$.
This contradicts the strong band irreducibility property of $T_W$, since at least one of the sets $A_3$ and $A_3^c$ has positive measure.
Thus, $T_W^2=T_K$ is a compact, strongly band irreducible operator on $L^2([0, 1])$ with nonnegative kernel $K$.
Since $\| T_W^2\|_{2, 2}=1$ and $T_W^2\psi =\psi$, by \cite[Theorem~43.8]{zaanen-book-operator-theory}, the eigenspace of $T_W^2$ associated with the eigenvalue $1$ is spanned by $\psi$.
If there exists $\phi\in L^2([0, 1])$ with $\|\phi\|_2=1$ and $T_W\phi = -\phi$, then $T_W^2\phi=\phi$, and consequently, either $\phi=\psi$ or $\phi= -\psi$ which leads to a contradiction.
Hence, $-1$ is not an eigenvalue of $T_W$, and this completes the proof of Lemma~\ref{lem:1}~\eqref{item:23}.
\qed

\subsection{An intermediate result and related notation}\label{sec:g2}
As a first step, we will prove Theorem~\ref{thm:scaling-limit-graphon} when the following stronger assumption replaces Condition~\ref{ass:graphon-3}.

\begin{ass}\label{ass:graphon-3a}
Let $\theta_0$ be as in Condition~\ref{ass:graphon-2}.
Then there exist $\delta_0\in (0\, ,\, 1/3-2\theta_0)$ and $\varpi_0\in [0\, ,\, 2\delta_0+1/3)$ such that the following hold:
	For each $n\geq 1$, there exists $B_n\subseteq [n]$ satisfying
	\begin{gather}
		|B_n|
		=
		o(n^{\varpi_0}) \, , \ \ \ \
	     \beta_{ij }^{\sss(n)}	= \beta_{ji}^{\sss(n)}	 =  0
		\ \ \text{ if }\ \ (i, j)\in B_n\times [n]\, ,
		\label{eqn:31a}\\
		\beta_{ii }^{\sss(n)} = 0 \ \ \text{ if }\ \
		i\in [n]\setminus B_n\ \, ,\ \text{ and }\ \  
		\beta_{ij }^{\sss(n)}>n^{-\delta_0}
		\ \ \text{ if }\ \
		i, j\in [n]\setminus B_n\text{ with }i\neq j .
		\label{eqn:33a}
	\end{gather}
\end{ass}

\begin{prop}\label{prop:graphon-intermediate}
Under Conditions~\ref{ass:graphon-1}, \ref{ass:graphon-2}, and \ref{ass:graphon-3a}, the convergences in \eqref{eqn:65} and \eqref{eqn:66} hold with the same constants as in \eqref{eqn:def-alpha-chi-zeta}.
\end{prop}

The next few sections are devoted to the proof of Proposition~\ref{prop:graphon-intermediate}.
For the reader's convenience, we collect some notation here that will be used frequently in the proof.
We define 
\begin{align*}
\kappa^{\sss(n)}_{ij}
:=
\left\{
\begin{array}{l}
e^{-n}\, ,\ \ \text{ if } i=j \in [n]\setminus B_n\, ,\\[3pt]
\beta_{ij }^{\sss(n)}-n^{-\delta_0}\, ,\ \ \text{ if } i, j\in [n]\setminus B_n\text{ with }i\neq j\, ,\\[3pt]
0\, ,\ \ \text{ otherwise},
\end{array}
\right.
\end{align*}
and
\begin{align}\label{eqn:88}
\kappa_n(x, y)
:=
\left\{
\begin{array}{l}
\kappa_{ij }^{\sss(n)}\, ,\ \  \text{ if } (x, y)\in Q_{ij}^{\sss(n)},\ \ 1\leq i,j\leq n,\\[5pt]
0,\ \ \text{ if } (x, y)\in [0, 1]^2\setminus \bigcup_{i, j\in [n]}Q_{ij}^{\sss(n)}\, .
\end{array}
\right.
\end{align}
Let $\cG_n^{\bullet}$ denote the random graph on the vertex set $[n]$ where the edges $\{i, j\}$ appear independently with respective probabilities $1\wedge(\kappa^{\sss(n)}_{ij}/n)$, $1\leq i<j\leq n$.
For a finite graph $G$ and $k\geq 1$, define 
\[
s_k(G) := \sum_{\cC\subseteq G}\ \frac{|\cC|^k}{|G|}
\ \  \text{ and }\ \
\cD(G) := \frac{1}{|G|}
\sum_{\cC\subseteq G}
\sum_{i\in V(\cC)} 
\sum_{j\in V(\cC)} d_{\cC}(i,j) \, ,
\]
where $\sum_{\cC\subseteq G}$ denotes the sum over all connected components of $G$.
The function $s_k(G)$ is called the $k$-th susceptibility function.
We will write
\begin{align}\label{eqn:898}
D_{\max}^{\bullet} := \diam(\cG_n^{\bullet})\, ,\ \ 
\cD^{\bullet}:= \cD(\cG_n^{\bullet})\, ,\ \ \text{ and }\ \
s_k^{\bullet} := s_k(\cG_n^{\bullet})\, , \ k\geq 1\, .
\end{align}
We write $\cC_i^{\bullet}$ for the $i$-th largest component in $\cG_n^{\bullet}$, and let $\cC^{\bullet}(j)$ denote the  component in $\cG_n^{\bullet}$ containing the vertex $j\in [n]$.
For $\ell\in\bZ_{>0}$ and $i, j\in [n]$ with $i\neq j$, we write $i\connectsell j$ for the event that there is a self-avoiding path with $\ell$ edges in $\cG_n^{\bullet}$ whose endpoints are $i$ and $j$.
We will write $i\connects j$ for the event $i\in \cC^{\bullet}(j)$
and $i\leftrightarrow j$ for the event that there is an edge between $i$ and $j$ in $\cG_n^{\bullet}$. 
(Thus, $\big\{ i\leftrightarrow j \big\} = \big\{ i\stackrel{1}{\leftrightsquigarrow} j \big\}$.)

It follows from \eqref{eqn:20} that 
$\max_{i, j\in [n]}\beta^{\sss(n)}_{ij}
=
O(n^{2/3})$.
In fact, we will prove later (see \eqref{eqn:22-b}) that 
$\max_{i, j\in [n]}\beta^{\sss(n)}_{ij}
=
o(n^{2/3})$.
Consequently,
\begin{align}\label{eqn:88-a}
\max_{i, j\in [n]}\kappa^{\sss(n)}_{ij}
\leq
e^{-n} + \max_{i, j\in [n]}\beta^{\sss(n)}_{ij}
=
o(n^{2/3})\, .
\end{align}
Hence, we will simply write $(\kappa^{\sss(n)}_{ij}/n)$ and $(\beta^{\sss(n)}_{ij}/n)$ for the edge probabilities in the sequel (i.e., we will omit the `$\wedge 1$').
Further, from now on, we will drop the superscript `$(n)$' and write $\beta_{ij}$, $\kappa_{ij}$, and $Q_{ij}$.
We let
\[
T_n=T_{\kappa_n}
\ \ \ \text{ and }\ \ \ 
\|T_n\|=\|T_{\kappa_n}\|_{2,2}\, .
\]
We make note of a very useful bound here: for any $p, q\in [3/2, 3]$ and for all $n$,
\begin{align}\label{eqn:T-n-norm-bound}
\| T_n \|_{p, q}
\leq
\big(\medint{\int}_{[0,1]^2} \kappa_n^3\big)^{1/3}
\leq
\big( e^{-3n} + \medint{\int}_{[0, 1]^2} W_n^3\big)^{1/3}
\leq 
C\, ,
\end{align}
where the first step uses Lemma~\ref{lem:2}, and the last step follows from \eqref{eqn:20}.

Note that if Condition~\ref{ass:graphon-3a} (or Condition~\ref{ass:graphon-3}) is satisfied for some choice of $\big\{\delta_0, \varpi_0, (B_n, n\geq 1)\big\}$ with $\delta_0\in (0\, ,\, 1/3-2\theta_0)$, $\varpi_0\in [0\, ,\, 2\delta_0+1/3)$, and $B_n\subseteq [n]$, $n\geq 1$, then the condition is also satisfied for the choice $\big\{\delta_1, \varpi_0, (B_n, n\geq 1)\big\}$ for each
$\delta_1\in (\delta_0, 1/3-2\theta_0)$.
Since $\theta_0<1/12$ by assumption and consequently $1/3-2\theta_0>1/6$, we will, from now on, assume that
\begin{align}\label{eqn:17}
\delta_0\in (1/6, 1/3-2\theta_0)\, .
\end{align}

There is an obvious correspondence between the non-zero eigenvalues of the (compact and self-adjoint) operator $T_n$ (resp. their associated eigenfunctions) and the non-zero eigenvalues of the symmetric matrix $((\kappa_{ij}))_{n\times n}$ (resp. their associated eigenvectors).
Under Condition~\ref{ass:graphon-3a}, this matrix has a submatrix of dimension $(n-|B_n|)\times (n-|B_n|)$ with positive entries and all its other entries are zeros.
Thus, by the Perron-Frobenius theorem, $\|T_n\|=\|T_n\|_{2,2}$ is a simple eigenvalue of $T_n$, and the corresponding eigenspace is spanned by a nonnegative eigenfunction.
Further, we can write
\begin{align}\label{eqn:53}
\kappa_n(x, y)
=
\|T_n\|\psi_n(x)\psi_n(y)
+
\sum_{j=2}^n\lambda_{n;j}\psi_{n;j}(x)\psi_{n;j}(y)\, ,\ \ x, y\in[0, 1]\, ,
\end{align}
where 
(a) $\lambda_{n;j}$, $2\leq j\leq n$, are eigenvalues of $T_n$ (some of which may be zeros) satisfying 
$|\lambda_{n;n}| \leq \cdots \leq |\lambda_{n;2}| < \|T_n\|$, 
(b) $\fF:=\big\{\psi_n, \psi_{n;2}, \ldots, \psi_{n;n}\big\}$ is an orthonormal collection of functions in $L^2([0,1])$ such that each $f\in\fF$ is constant on $\big( (i-1)/n ,\, i/n \big)$ for $i\in [n]$ and $f(i/n)=0$ for $i=0, 1,\ldots, n$, and
(c) $\psi_n(x)\geq 0$ for $x\in [0, 1]$ with $\mu_{\circ}\big(\{x\, :\, \psi_n(x)=0\}\big)=|B_n|/n$.

Our aim is to use Theorem~\ref{thm:gen-2} to prove Proposition~\ref{prop:graphon-intermediate}.
In the present setting, the functionals appearing in Conditions~\ref{ass:aldous-basic-assumption} and \ref{ass:gen-2} will be rescaled versions of $s_2^{\bullet}$, $s_3^{\bullet}$, $\cD^{\bullet}$, $D_{\max}^{\bullet}$, and $|\cC_1^{\bullet}|$, and consequently, we need to derive asymptotics for these random variables.
Our strategy is to prove the relevant asymptotics for a closely related branching process and then transfer the results to the above random variables by approximation.
We now define this branching process and introduce some associated functions.

Consider the multitype branching process with type space $[n]$ where every vertex of type $i\in [n]$ has $\text{Bernoulli}(\kappa_{ij}/n)$ many type $j$ children for each $j\in [n]$.
For $i\in[n]$ and $\ell\geq 0$, let $\fG_{\ell}^{\sss\triangle}(i)=\fG_{n;\ell}^{\sss\triangle}(i)$ denote the collection of vertices in generation $\ell$ of such a branching process started from a single vertex of type $i$.
(Thus, $\fG_0^{\sss\triangle}(i)$ consists of a single vertex of type $i$.)
Let $\fX_n^{\sss\triangle}(i):=\sum_{\ell\geq 0}|\fG_{\ell}^{\sss\triangle}(i)|$ be the total progeny, and let
$\fH_n^{\sss\triangle}(i):=\max\big\{\ell\geq 0\, :\, \fG_{\ell}^{\sss\triangle}(i)\neq \emptyset\big\}$ denote the height.
Further, define
\begin{align}\label{eqn:52}
\fZ_n^{\sss\triangle}(i):=\sum_{\ell\geq 1}\ell\cdot |\fG_{\ell}^{\sss\triangle}(i)| \, ,\ \text{ and }\
g_{n;k}^{\sss\triangle}(i):=\bE\big[\fX_n^{\sss\triangle}(i)^k\big]\, ,\ \ k\geq 1\, .
\end{align}
We will simply write $g_n^{\sss\triangle}(i)$ instead of $g_{n; 1}^{\sss\triangle}(i)$.
For $x\in [0, 1]$, we define
\begin{align}\label{eqn:52a}
g_n(x)
:=
\left\{
\begin{array}{l}
g_n^{\sss\triangle}(\lceil nx\rceil)\, , \text{ if } x\in [0, 1]\setminus\{0, 1/n, 2/n, \ldots, 1\} \, ,\\[2pt]
0,\text{ if } x\in \{0, 1/n, 2/n, \ldots, 1\}\, .
\end{array}
\right.
\end{align}
Similarly, for $x\in [0, 1]$, $|\fG_{\ell}(x)|, \fX_n(x), \fH_n(x), \fZ_n(x), g_{n;k}(x)$ etc. are defined in an analogous manner.
Note that the functions $g_{n; k}$ can be viewed as elements of $L^p([0, 1])$ for any $p\in [1, \infty)\cup\{\infty\}$.
Further, if $v$ is uniformly distributed on $[n]$, then
$\bE\big[g_{n;k}^{\sss\triangle}(v)\big]=\int_0^1 g_{n;k}(x) dx$.

We define $\mvone_{\circ} = \mvone_{\circ}^{\sss (n)}: [0, 1]\to\bR$ by
$\mvone_{\circ}(x) := \ind_{A}(x)$ where $A = [0, 1]\setminus\{0,\, 1/n, 2/n, \ldots, 1\}$.
We let $\mvone_{\sss\square} = \mvone_{\sss\square}^{\sss (n)}: [0, 1]^2\to\bR$ be the function 
$\mvone_{\sss\square}(x, y) = \mvone_{\circ}(x) \mvone_{\circ}(y)$.

\subsection{Asymptotics for the eigenfunction $\psi_n$}\label{sec:g3}
Throughout this section, we work under Conditions~\ref{ass:graphon-1}, \ref{ass:graphon-2}, and \ref{ass:graphon-3a}.
\begin{lem}\label{lem:3}
We have, as $n\to\infty$,
\begin{gather}
\|\psi_n - \psi\|_3 = o(1)\, ,\ \ \text{ and}\label{eqn:22}\\
\|\psi_n - \psi\|_2 = O(n^{-\delta_0})\, .\label{eqn:22-a}
\end{gather}
\end{lem}

The following result will be needed in the proof.
\begin{lem}\label{lem:4}
We have, 
$
\|T_{W_n} - T_W\|_{3,3} \to 0
$
as $n\to\infty$.
\end{lem}

\noindent{\bf Proof:}
We will first show that
\begin{align}\label{eqn:22-b}
	\max_{i, j\in [n]}\, \beta_{ij} = o(n^{2/3})\ \ \text{ as }\ \ n\to\infty.
\end{align}
If \eqref{eqn:22-b} fails then there exist $\eps_0>0$, $n_1<n_2<\cdots$, and $i_k, j_k\in [n_k]$ such that $\beta_{i_k j_k}\geq \eps_0 n_k^{2/3}$ for all $k\geq 1$.
For $k\geq 1$, let $f_k:=\sqrt{n_k}\cdot\ind_{J_k}$, where $J_k=\big((j_k -1)/n_k,\, j_k/n_k\big)$.
Then
\begin{align*}
|(T_H f_k)(x)|
=
\sqrt{n_k}\cdot \big|\medint\int_{y\in J_k} H(x, y) dy\big|
\leq
\Big(\medint\int_{y\in J_k} H^2(x, y) dy\Big)^{1/2} \, ,
\end{align*}
and consequently,
\begin{align}\label{eqn:32}
\| T_H f_k\|_2
\leq
\Big(\medint\int_{[0, 1]\times J_k} H^2(x, y) dx dy\Big)^{1/2}
=o(1) \ \ \text{ as }\ \ k\to\infty\, .
\end{align}
Let $I_k=\big((i_k -1)/n_k,\, i_k/n_k\big)$.
Then
\begin{align}\label{eqn:32-a}
\medint\int_{x\in I_k}\big(T_W f_k (x)\big)^2 dx
&=
n_k \medint\int_{x\in I_k}\Big( \medint\int_{y\in J_k} W(x, y) dy\Big)^2 dx
\notag\\
&\leq 
\medint\int_{I_k\times J_k} W^2
\leq
\Big(\medint\int_{Q_{i_k j_k}} W^3\Big)^{2/3}
\cdot n_k^{-2/3}
=
o(n_k^{-2/3})\, ,
\end{align}
as $k\to\infty$.
Hence,
\begin{align}\label{eqn:27}
\| T_{W_{n_k}}f_k - T_Wf_k \|_2
&\geq
\Big[
\medint\int_{x\in I_k}
\Big(
T_{W_{n_k}}f_k(x)
-
T_Wf_k(x)
\Big)^2 dx
\Big]^{1/2}
\notag\\
&\geq
\Big[
\medint\int_{x\in I_k}
\Big(
T_{W_{n_k}}f_k(x)
\Big)^2 dx
\Big]^{1/2}
-
\Big[
\medint\int_{x\in I_k}
\Big(
T_Wf_k(x)
\Big)^2 dx
\Big]^{1/2}
\notag\\
&=
n_k^{-1}\beta_{i_k j_k}
-
\Big[
\medint\int_{x\in I_k}
\Big(
T_Wf_k(x)
\Big)^2 dx
\Big]^{1/2}
\geq 
\eps_0 n_k^{-1/3}\big(1-o(1)\big)\, ,
\end{align}
where the second step follows from the Minkowski inequality and the last step uses \eqref{eqn:32-a}.
Combining \eqref{eqn:27} with \eqref{eqn:32}, we see that
\[
\liminf_{k\to\infty}\,
\big\| 
n_k^{1/3}\big( T_{W_{n_k}}f_k - T_Wf_k \big)- T_H f_k 
\big\|_2
\geq \eps_0\, ,
\]
which contradicts \eqref{eqn:21}.
This proves \eqref{eqn:22-b}.

Fix $\eps>0$ and define
$W_n^{\sss(\leq)}:= W\ind_{\{W\leq\eps n^{2/3}\}}$ and $W_n^{\sss (>)}:=W-W_n^{\sss(\leq)}$.
Then
\begin{gather}
\|T_{\sss W_n^{\sss(>)}}\|_{3,3}
\leq
\Big(\medint\int_{W>\eps n^{2/3}} W^3\Big)^{1/3}
\, ,\ \ \text{ and }\label{eqn:4a}\\[2pt]
\|T_{\sss W_n^{\sss(>)}}\|_{2,2}
\leq
\Big(\medint\int_{W>\eps n^{2/3}} W^2\Big)^{1/2}
\leq
\frac{1}{\eps^{1/2} n^{1/3}}
\Big(\medint\int_{W>\eps n^{2/3}} W^3\Big)^{1/2}\, .\label{eqn:4b}
\end{gather}
Using \eqref{eqn:22-b}, we can choose $n_*\geq 1$ such that $\max_{i, j\in [n]}\beta_{ij}\leq \eps n^{2/3}$ for all $n\geq n_*$.
Let $f\in L^3([0, 1])$ with $\| f\|_3=1$.
Then for $n\geq n_*$ and $x\in [0, 1]$,
\begin{align}\label{eqn:58-a}
\big| 
T_{W_n}f(x)-T_{W_n^{\sss(\leq)}}f(x)
\big|
\leq
\medint\int_0^1 \big|W_n(x, y)-W_n^{\sss(\leq)}(x, y)\big|\cdot \big|f(y)\big|\, dy
\leq
2\eps n^{2/3} \|f\|_1
\leq
2\eps n^{2/3} ,
\end{align}
where the last step uses the fact $\|f\|_1\leq \|f\|_3=1$.
Hence, for $n\geq n_*$,
\begin{align}\label{eqn:311}
&\|T_{W_n}f-T_{\sss W_n^{\sss(\leq)}}f\|_3^3
\leq
2\eps n^{2/3} \|T_{W_n}f-T_{\sss W_n^{\sss(\leq)}}f\|_2^2 
\leq
2\eps n^{2/3} \|T_{W_n}-T_{\sss W_n^{\sss(\leq)}}\|_{\sss 2, 2}^2 \cdot\|f\|_2^2\notag\\
&\hskip50pt
\leq
2\eps n^{2/3} \|T_{W_n}-T_{\sss W_n^{\sss(\leq)}}\|_{\sss 2, 2}^2
\leq
4\eps n^{2/3} 
\Big(\|T_{W_n}-T_W\|_{\sss 2, 2}^2 + \|T_{\sss W_n^{\sss(>)}}\|_{\sss 2,2}^2\Big) \notag\\
&\hskip100pt
\leq
4\eps n^{2/3} \|T_{W_n}-T_W\|_{\sss 2, 2}^2 + 4\medint\int_{W>\eps n^{2/3}} W^3\, ,
\end{align}
where the first step uses \eqref{eqn:58-a}, the third step uses the fact $\|f\|_2\leq \|f\|_3 = 1$, and the last step follows from \eqref{eqn:4b}.
Since
$
\|T_{W_n}f-T_Wf\|_3
\leq
\|T_{W_n}f-T_{\sss W_n^{\sss(\leq)}}f\|_3 + \| T_{\sss W_n^{\sss(>)}}\|_{3,3} \,
$, it follows from \eqref{eqn:4a}, \eqref{eqn:311}, and \eqref{eqn:21} that for all large $n$,
\[
\|T_{W_n}-T_W\|_{3,3}
\leq
\Big(
4\eps\|T_H\|_{2,2}^2 + \eps + 4\medint\int_{W>\eps n^{2/3}} W^3
\Big)^{1/3}
+
\Big(\medint\int_{W>\eps n^{2/3}} W^3\Big)^{1/3}\, ,
\]
which yields the desired claim.
\qed

\medskip

\noindent{\bf Proof of Lemma~\ref{lem:3}:}
Note that for $p=2, 3$,
\begin{align}\label{eqn:101}
\|T_n-T_{W_n}\|_{p, p}
\leq
\Big[ \int_{[0, 1]^2}
\big( \kappa_n - W_n \big)^p
\Big]^{1/p}
\lesssim
n^{-\delta_0}\, ,
\end{align}
where the last step uses the definition of $\kappa_n$ from \eqref{eqn:88}.
Hence,
\begin{gather}
\| T_n -T_W\|_{3,3} =o(1) \ \ \text{ and} \label{eqn:102}\\
\big|1-\| T_n\| \big|
\leq
\| T_n -T_W\|_{2,2} =O(n^{-\delta_0}) \label{eqn:103}
\end{gather}
as $n\to\infty$, where \eqref{eqn:102} uses \eqref{eqn:101} and Lemma~\ref{lem:4}, and  \eqref{eqn:103} follows from \eqref{eqn:101} and \eqref{eqn:21}.
Now,
\[
\| T_n \|\cdot \|\psi_n\|_3
=
\|T_n\psi_n\|_3
\leq
\|T_n\|_{2,3}
\leq
C\, ,
\]
where the last step uses \eqref{eqn:T-n-norm-bound}.   
It thus follows from \eqref{eqn:103} that $\big(\psi_n\, ,\, n\geq 1\big)$ is a bounded sequence in $L^3([0, 1])$.
Consider the sequence $\big(T_W\psi_n\, ,\, n\geq 1\big)$.
Since $W\in L^3([0, 1]^2)$, 
\[
\medint{\int_{x=0}^1}\big(
\medint{\int_{y=0}^1} W(x, y)^{3/2} dy
\big)^2 
\, dx
\leq
\medint{\int}_{[0, 1]^2} W^3
<\infty\, .
\]
Thus, by \cite[Chapter~21, Theorem~41.6]{zaanen-book-operator-theory}, $T_W$ is a compact operator on $L^3([0, 1])$, and consequently, for any subsequence 
$\big(T_W\psi_{n_k}\, ,\, k\geq 1\big)$, there exists $\phi\in L^3([0, 1])$ and a further subsequence 
$\big(T_W\psi_{n_{k_{\ell}}} \, ,\, \ell\geq 1\big)$ such that
\begin{align}\label{eqn:666}
\|T_W\psi_{n_{k_{\ell}}} - \phi\|_3 \to 0\ \ \text{ as }\ \ \ell\to\infty\, .
\end{align}
Writing 
\begin{align}\label{eqn:76}
T_W\psi_n
=
\psi_n + (T_W-T_n)\psi_n + \big(\|T_n\| - 1\big)\psi_n\, ,
\end{align}
and using \eqref{eqn:102}, \eqref{eqn:103} in conjunction with \eqref{eqn:666} yields 
\begin{align}\label{eqn:667}
\|\psi_{n_{k_{\ell}}} - \phi\|_3 \to 0\ \ \text{ as }\ \ \ell\to\infty\, .
\end{align}
This in turn implies that $\| T_W\psi_{n_{k_{\ell}}} - T_W\phi\|_3 \to 0$, which combined with \eqref{eqn:666} yields 
\begin{align}\label{eqn:668}
T_W\phi=\phi\, .
\end{align}
Now, \eqref{eqn:667} shows that $\phi(x)\geq 0$ for $\mu_{\circ}$ a.e. $x$ and $\|\phi\|_2=1$.
Hence, using \eqref{eqn:668} and Lemma~\ref{lem:1}, we conclude that $\phi=\psi$.
Thus, for any subsequence we can extract a further subsequence along which $\|\psi_n-\psi\|_3\to 0$.
This proves \eqref{eqn:22}.

Since $T_W$ is a compact and self-adjoint operator on $L^2([0, 1])$, we can write
$
T_W
=
\langle \psi,\cdot\,\rangle \psi
+
\sum_{j\geq 2}\lambda_j \langle \phi_j,\cdot\,\rangle \phi_j$,
where $\lambda_j, j\geq 2$, are the non-zero eigenvalues (if any) of $T_W$ other than $1$ counted according to their multiplicity, and $\{\psi\}\cup\{\phi_j\, :\, j\geq 2\}$ is an orthonormal collection of elements of $L^2([0, 1])$.
Let $\cP: L^2([0, 1])\to\text{kernel}(T_W)$ denote the orthogonal projection operator onto the kernel of $T_W$.
Then
\begin{align}\label{eqn:232}
\psi_n
=
\langle \psi, \psi_n\rangle \psi
+
\sum_{j\geq 2}\langle \phi_j, \psi_n\rangle \phi_j
+
\cP\psi_n\, ,\ \text{ and }\ 
T_W\psi_n
=
\langle \psi, \psi_n\rangle \psi
+
\sum_{j\geq 2}\lambda_j \langle \phi_j, \psi_n\rangle \phi_j\, .
\end{align}
Hence,
\begin{align}
&Cn^{-2\delta_0}
\geq
\big\|(T_W-T_n)\psi_n + \big(\|T_n\| - 1\big)\psi_n\big\|_2^2
=
\sum_{j\geq 2}\big(1-\lambda_j\big)^2 \langle \phi_j, \psi_n\rangle^2
+
\|\cP\psi_n\|_2^2
\notag
\\
&\hskip80pt
\geq
\Delta_W^2\Big(
\sum_{j\geq 2}\langle \phi_j, \psi_n\rangle^2
+
\|\cP\psi_n\|_2^2
\Big)
=
\Delta_W^2\big(1- \langle \psi, \psi_n\rangle^2\big)\, ,\label{eqn:5a}
\end{align}
where the first step uses \eqref{eqn:103}, the second step uses  \eqref{eqn:232}, and the third step follows from Lemma~\ref{lem:1}~\eqref{item:23}. 
We know from \eqref{eqn:22} that $\|\psi_n-\psi\|_2=o(1)$, which combined with \eqref{eqn:5a} yields the lower bound 
$\langle \psi, \psi_n\rangle \geq 1-C n^{-2\delta_0}$.
Consequently,
\[
\|\psi-\psi_n\|_2^2
=
\|\psi\|_2^2+\|\psi_n\|_2^2-2\langle \psi, \psi_n\rangle 
=
2-2\langle \psi, \psi_n\rangle 
=
O(n^{-2\delta_0})\, ,
\]
as claimed in \eqref{eqn:22-a}.
\qed

\subsection{Asymptotics for $\| T_n \|$ and $\| g_n\|_1$}\label{sec:g4}
Throughout this section, we will work under Conditions~\ref{ass:graphon-1}, \ref{ass:graphon-2}, and \ref{ass:graphon-3a}.
Recall the constant $\zeta$ from \eqref{eqn:def-alpha-chi-zeta} and the function $g_n=g_{n; 1}$ from \eqref{eqn:52}.
Recall also from \eqref{eqn:17} that we are working with $\delta_0>1/6$.
The aim of this section is to prove the following result:

\begin{lem}\label{lem:20}
We have
\[
\lim_{n\to\infty}
\frac{\big(\int_0^1 g_n(x)dx\big)-n^{\delta_0}}{n^{2\delta_0-1/3}}
=
\zeta\, .
\]
\end{lem}

We will make use of the following result:

\begin{lem}\label{lem:99}
Suppose $A_1, A_2$, and $A_3$ are compact, self-adjoint operators on $L^2([0, 1])$.
Further, assume that $\|A_1\|_{\sss 2,2}$ is an eigenvalue of $A_1$ with multiplicity one, and 
$-\|A_1\|_{\sss 2,2}$ is not an eigenvalue of $A_1$.
Then
\begin{align}\label{eqn:86}
\lim_{(x, y)\to (0, 0)} \ 
\frac{1}{y}\cdot
\Big(
\|A_1 + xA_2 + yA_3\|_{\sss 2,2}
-
\|A_1 + xA_2\|_{\sss 2,2}
\Big)
=
\big\langle \phi_{A_1} ,\, A_3\phi_{A_1} \big\rangle\, ,
\end{align}
where $\phi_{A_1}$ is an eigenfunction of $A_1$ corresponding to the eigenvalue $\|A_1\|_{2,2}$ with $\|\phi_{A_1}\|_2=1$.
\end{lem}

\noindent{\bf Proof:}
We claim that there exist an open $U\subseteq\bR^2$ containing the origin,
 and infinitely differentiable functions $\Lambda: U\to\bR$ and $\Phi: U\to L^2([0, 1])$ (the latter differentiability to be interpreted in the Fr\'echet sense) such that 
$\Lambda(0, 0)=\|A_1\|_{\sss 2,2}\, $,  $\Phi(0, 0)=\phi_{A_1}$, and $\|\Phi(x, y)\|_2=1$ for $(x, y)\in U$.
Further, for each $(x, y)\in U$, $\Lambda(x, y)$ is an eigenvalue of $(A_1 + xA_2 + yA_3)$, and $\Phi(x, y)$ is an associated eigenfunction.

Let us first complete the proof assuming the above claim.
Choose $\eps>0$ small so that $(-\eps, \eps)^2\subseteq U$.
For $(x, y)\in (-\eps, \eps)^2$,
\begin{align}
	&
	\big\| 
	\big(A_1 + xA_2 + yA_3 \big)
	- 
	\Lambda(x, y)\big\langle \Phi(x, y)\, , \, \cdot\, \big\rangle \Phi(x, y)
	\big\|_{\sss 2,2}
	\notag\\
	&\hskip30pt
	\leq
	\eps\|A_2\|_{\sss 2,2}
	+
	\eps\|A_3\|_{\sss 2,2}
	+
	\big\| 
	A_1 
	-
	\|A_1\|_{\sss 2,2} \big\langle \phi_{A_1}\, , \, \cdot\, \big\rangle \phi_{A_1}
	\big\|_{\sss 2,2}
	\notag\\
	&\hskip60pt
	+
	\big\| 
	\|A_1\|_{\sss 2,2} \big\langle \phi_{A_1}\, , \, \cdot\, \big\rangle \phi_{A_1}
	-
	\Lambda(x, y)\big\langle \Phi(x, y)\, , \, \cdot\, \big\rangle \Phi(x, y)
	\big\|_{\sss 2,2}\, .
	\label{eqn:900}
\end{align}
By our assumption,
$\big\| 
A_1 
-
\|A_1\|_{\sss 2,2} \big\langle \phi_{A_1}\, , \, \cdot\, \big\rangle \phi_{A_1}
\big\|_{\sss 2,2}
<
\|A_1\|_{\sss 2,2} 
=
\Lambda(0, 0)
$, 
and consequently, the right side of \eqref{eqn:900} is strictly smaller \chsen{than} $\Lambda(x, y)$ (by choosing a sufficiently small $\eps$).
Hence, for each $(x, y)\in (-\eps, \eps)^2$, 
$\Lambda(x, y)=\| A_1 + xA_2 + yA_3\|_{\sss 2,2}$ and further, it is an eigenvalue of $(A_1 + xA_2 + yA_3)$ of multiplicity one and the corresponding (one-dimensional) eigenspace is spanned by $\Phi(x, y)$.
Thus, the left side of \eqref{eqn:86} is simply 
$\frac{\partial}{\partial y}\Lambda (0, 0)$,
so we just need to show that
\begin{align}\label{eqn:600}
\frac{\partial}{\partial y}\Lambda (0, 0)
=
\big\langle \phi_{A_1} ,\, A_3\phi_{A_1} \big\rangle\, .
\end{align}
To this end, for $y\in (-\eps, \eps)$, write $\Phi_2(y)=\Phi(0, y)$, $D_y\Phi_2$ for the Fr\'echet derivative of $\Phi_2$, and identify the linear operator 
$D_y\Phi_2 (y): \bR\to L^2([0, 1])$ with $\big[D_y\Phi_2(y)\big](1)\in L^2([0, 1])$ .
Then differentiating both sides of the relation 
$\langle \Phi_2(y), \Phi_2(y)\rangle=1$, $y\in (-\eps, \eps)$, we get
\begin{align}\label{eqn:87}
\big\langle \Phi_2(y),\, D_y \Phi_2(y)  \big\rangle =0\ \ \text{ for all }\ \ y\in (-\eps, \eps)\, .
\end{align}
Hence,
\begin{align}
\frac{1}{y}\Big(
\big\langle \Phi_2(y)\, ,\, A_1 \Phi_2( y)\big\rangle
-
\Lambda(0, 0)
\Big)
\stackrel{\sss y\to 0}{\longrightarrow}
2\Lambda(0, 0)\cdot \big\langle \Phi_2(0),\, D_y \Phi_2(0) \big\rangle=0\, ,
\label{eqn:300}
\end{align}
where the last equality follows from \eqref{eqn:87}.
Now, $\big\langle \Phi_2(y)\, ,\, (A_1+yA_3)\Phi_2( y)\big\rangle=\Lambda(0, y)$. 
Consequently,
\begin{align*}
y\cdot\big\langle \Phi_2(y),\, A_3 \Phi_2(y)\big\rangle
&=
\Lambda(0, y)-\big\langle \Phi_2(y)\, ,\, A_1 \Phi_2( y)\big\rangle \\
&
=
\Big(\Lambda(0, y)-\Lambda(0, 0)\Big)
+
\Big(\Lambda(0, 0)
-
\big\langle \Phi_2(y)\, ,\, A_1 \Phi_2( y)\big\rangle
\Big)\, .
\end{align*}
Now \eqref{eqn:600} follows upon dividing by $y$ and letting $y\to 0$ in the last equation and using \eqref{eqn:300}.

It remains to prove the existence of the functions $\Lambda$ and $\Phi$.
This is well-known in the setting of one parameter perturbations of matrices, and a proof using the implicit function theorem can be found in the note \cite{kazdan}. 
A similar approach can be taken in the present setting; we only outline the steps.
Consider the function $F:\bR\times\bR\times L^2([0, 1]) \times\bR\to L^2([0, 1])\times\bR$ given by
\[
\big(x, y, \phi, \lambda\big)
\mapsto
\Big(
\big(A_1 + xA_2 + yA_3\big)\phi  - \lambda\phi
\, ,\ 
\|\phi\|_2^2
\Big)
\, .
\]
Then $F\big(0, 0, \phi_{A_1} , \|A_1\|_{\sss 2,2}\big)=(\mvzero, 1)$, where $\mvzero$ is the element of $L^2([0, 1])$ that is identically zero.
Further, $F$ is infinitely differentiable.
A direct calculation shows that the Fr\'echet derivative $D_{\phi, \lambda}F$, at the point  
$\big(0, 0, \phi_{A_1} , \|A_1\|_{\sss 2,2}\big)$, is the following linear transformation on 
$L^2([0, 1]) \times\bR$:
\[
L^2([0, 1]) \times\bR
\ni
(\phi, \lambda)
\mapsto
\big(
A_1\phi- \|A_1\|_{\sss 2,2}\phi - \lambda\phi_{A_1}\, ,\,
2\langle \phi, \phi_{A_1}\rangle
\big)\, .
\]
It can be easily shown that this is a one-to-one and onto linear transformation (and consequently its inverse is a bounded linear operator).
Thus, the implicit function theorem in Banach spaces \cite[Chapter I, Theorem~5.9]{lang-differential-geometry} guarantees the existence of an open $U\ni (0, 0)$, and infinitely differentiable functions 
$\Lambda: U\to\bR$ and $\Phi: U\to L^2([0, 1])$ having the desired properties.
This completes the proof.
\qed

\medskip

\noindent{\bf Proof of Lemma~\ref{lem:20}:}
Recall the notation $\mvone_{\circ}$ and $\mvone_{\sss\square}$ introduced at the end of Section~\ref{sec:g2}.
Let 
\[
\kappa_n^-:= W+n^{-1/3}\cdot H-n^{-\delta_0}\cdot\mvone_{\sss\square}
\ \ \ \text{ and }\ \ \
\kappa_n^{=}:= W_n-n^{-\delta_0}\cdot\mvone_{\sss\square}\, .
\]
Then
\begin{align}\label{eqn:760}
&\big|
\|T_n\| - \|T_{\kappa_n^-}\|_{\sss 2,2}
\big|^2
\leq 
2
\|T_n - T_{\kappa_n^{=}}\|_{\sss 2,2}^2
+
2
\|T_{\kappa_n^{=}} - T_{\kappa_n^-}\|_{\sss 2,2}^2
\notag\\
&\hskip40pt
\leq
2\int_{[0, 1]^2} \big(\kappa_n - \kappa_n^{=}\big)^2
+
o(n^{-2/3})
\leq
\frac{4}{n^{2\delta_0}}\cdot\frac{|B_n|+1}{n}
+
4e^{-2n}
+
o(n^{-2/3})
\notag\\[1pt]
&\hskip80pt
=
o(n^{-2\delta_0-1+\varpi_0})
+
o(n^{-2/3})
=
o(n^{-2/3})\, ,
\end{align}
where the second step uses \eqref{eqn:21} and the last two steps use Condition~\ref{ass:graphon-3a}.
Now, as seen in the proof of Lemma~\ref{lem:99},
the function $x\mapsto \| T_W - xT_{\mvone_{\sss\square}} \|_{\sss 2,2}$ is infinitely differentiable in a neighborhood of the origin.
Thus, using Lemma~\ref{lem:99} and the fact $\|T_W\|_{\sss 2,2}=1$, we get
\begin{align}\label{eqn:654}
1 - \| T_W - n^{-\delta_0}T_{\mvone_{\sss\square}} \|_{\sss 2,2}
=
n^{-\delta_0}\big\langle \psi,\, T_{\mvone_{\sss\square}}\psi\big\rangle
+O(n^{-2\delta_0})
=
n^{-\delta_0}\big(
\medint\int_0^1 \psi
\big)^2
+O(n^{-2\delta_0})\, ,
\end{align}
Further,
\begin{align}\label{eqn:655}
&n^{\delta_0}\big(1 - \| T_n\|\big)
=
n^{\delta_0}\big(1 - \|T_{\kappa_n^-}\|_{\sss 2,2}\big)
+
o(n^{\delta_0-1/3})
\notag\\
&\hskip30pt
=
n^{\delta_0}\big(1 - \| T_W - n^{-\delta_0}T_{\mvone_{\sss\square}} \|_{\sss 2,2}\big)
+
O(n^{\delta_0-1/3})
\to
\big(
\medint\int_0^1 \psi
\big)^2
\end{align}
as $n\to\infty$, where the first step uses \eqref{eqn:760}, the second step uses the deifinition of $\kappa_n^-\,$, and the last step uses \eqref{eqn:654} and the fact $\delta_0<1/3$.
Next, recall \eqref{eqn:53}, and note that
\[
|\lambda_{n;2}|
=
\big\|
T_n - \|T_n\| \langle \psi_n\, ,\, \cdot\,\rangle \psi_n
\big\|_{\sss 2,2}
\
\stackrel{\sss n\to\infty}{\longrightarrow}
\
\big\|
T_W - \langle \psi\, ,\, \cdot\,\rangle \psi
\big\|_{\sss 2,2}
\leq 
1-\Delta_W\, ,
\]
where the penultimate step uses \eqref{eqn:22-a} and \eqref{eqn:103}, and the last step uses Lemma~\ref{lem:1}~\eqref{item:23}. 
In view of the last display and \eqref{eqn:655}, we can assume without loss of generality that
\begin{align}\label{eqn:765}
1/2\leq \|T_n\|<1
\ \ \ \text{ and }\ \ \
|\lambda_{n;2}|
\leq
1-(\Delta_W/2)\ \ \ \text{ for all }\ \ \ n\, .
\end{align}

Recall the definition of $|\fG_{\ell}(\cdot)|$ from around \eqref{eqn:52a}. 
It is easy to see that for any $\ell\geq 1$,
\begin{align}\label{eqn:105}
\bE\big(
|\fG_{\ell}(x)|
\big)
=
T_n^{\ell}\mvone_{\circ}(x)\, ,\ x\in [0, 1]\, ,
\ \text{ and }\ 
T_n^{\ell}
=
\|T_n\|^{\ell} \langle\psi_n\, ,\, \cdot\, \rangle\psi_n
+
\sum_{j=2}^n \lambda_{n;j}^{\ell} \langle\psi_{n;j}\, ,\, \cdot\, \rangle \psi_{n;j} \, .
\end{align}
It thus follows that
\begin{align}\label{eqn:540}
g_n
=
\frac{1}{1-\|T_n\|} \big\langle\psi_n\, ,\, \mvone_{\circ} \big\rangle \psi_n
+
\sum_{j=2}^{n}
\frac{1}{1-\lambda_{n;j}} \big\langle\psi_{n; j}\, ,\, \mvone_{\circ} \big\rangle \psi_{n;j} \, ,
\end{align}
and consequently
\begin{align*}
\medint{\int}_0^1 g_n(x) dx
=
\frac{1}{1-\|T_n\|} \big\langle\psi_n\, ,\, \mvone_{\circ} \big\rangle^2
+
\sum_{j=2}^{n}
\frac{1}{1-\lambda_{n;j}} \big\langle\psi_{n; j}\, ,\, \mvone_{\circ} \big\rangle^2
=
\frac{1}{1-\|T_n\|} \big(\medint{\int_0^1} \psi_n \big)^2
+
O(1)\, ,
\end{align*}
where we have used \eqref{eqn:765} and the fact that
$
\sum_{j=2}^n \big\langle\psi_{n; j}\, ,\, \mvone_{\circ} \big\rangle^2
\leq 1$.
Hence,
\begin{align}\label{eqn:11}
\frac{\big(\int_0^1 g_n\big)-n^{\delta_0}}{n^{2\delta_0-1/3}}
=
\frac{\big(\int_0^1\psi_n\big)^2 - n^{\delta_0}\big(1-\|T_n\|\big)}{(1-\|T_n\|) \cdot n^{2\delta_0-1/3}}
+o(1)\, .
\end{align}
It follows from \eqref{eqn:22-a} that
\begin{align}\label{eqn:881}
\medint{\int_0^1}\psi_n
=
\|\psi_n\|_1
=
\|\psi\|_1 + O(n^{-\delta_0})\, .
\end{align}
Hence,
\begin{align}
&\frac{\big(\int_0^1\psi_n\big)^2 - n^{\delta_0}\big(1-\|T_n\|\big)}{(1-\|T_n\|) \cdot n^{2\delta_0-1/3}}
\sim
\frac{\big(\int_0^1\psi_n\big)^2 
- 
n^{\delta_0}\big(1-\|T_n\|\big)}{\big(\int_0^1\psi\big)^2\cdot n^{\delta_0-1/3}}
=
\frac{\big(\int_0^1\psi\big)^2 
- 
n^{\delta_0}\big(1-\|T_n\|\big)}{\big(\int_0^1\psi\big)^2\cdot n^{\delta_0-1/3}}
+o(1) 
\notag\\
&\hskip40pt
=
\frac{
n^{\delta_0}\big(1 - \| T_W - n^{-\delta_0}T_{\mvone_{\sss\square}} \|_{\sss 2,2}\big)
- 
n^{\delta_0}\big(1-\|T_n\|\big)}{\big(\int_0^1\psi\big)^2\cdot n^{\delta_0-1/3}}
+o(1)
\notag\\
&\hskip80pt
=
\frac{
n^{\delta_0}\Big(
\|T_{\kappa_n^-}\|_{\sss 2,2} 
- 
\| T_W - n^{-\delta_0}T_{\mvone_{\sss\square}} \|_{\sss 2,2}
\Big)
}{\big(\int_0^1\psi\big)^2\cdot n^{\delta_0-1/3}
}
+o(1)
=
\frac{
\big\langle \psi,\, T_H\psi \big\rangle
}{
(\int_0^1\psi)^2
}
+o(1)\, , \label{eqn:569}
\end{align}
where the first step uses \eqref{eqn:655},
the second step uses \eqref{eqn:881},
the third step uses \eqref{eqn:654},
the fourth step uses \eqref{eqn:760},
and the last step uses Lemma~\ref{lem:99};
in addition, we have repeatedly used the fact that $\delta_0>1/6$.
Combining \eqref{eqn:11} and \eqref{eqn:569} completes the proof.
\qed

\subsection{Norm estimates for the functions $g_{n; k}\, $}\label{sec:g10}
Throughout this section, we work under Conditions~\ref{ass:graphon-1}, \ref{ass:graphon-2}, and \ref{ass:graphon-3a}.
Recall the functions $g_{n; k}$ from \eqref{eqn:52}.
Recall also that we are working under the assumptions of \eqref{eqn:765} and \eqref{eqn:17}.
The following lemma is the main result of this section:

\begin{lem}\label{lem:12}
The following hold:
\begin{enumeratea}
\item\label{item:11}
$\|g_n\|_3 = O(n^{\delta_0})$;

\medskip

\item\label{item:12}
$\|g_{n; 2}\|_3 = O(n^{3\delta_0})$
and further,
\begin{align}\label{eqn:555}
\frac{1}{n^{3\delta_0}}\cdot \medint\int_0^1 g_{n;2} 
\to 
\frac{\int_0^1\psi^3 }{ ( \int_0^1 \psi )^3}\ \ \ \text{ as }\ \ \ n\to\infty\, ;
\end{align}

\item\label{item:13}
$\|g_{n; k}\|_{3/2} = O\big(n^{(2k-1)\delta_0}\big)$ for $k=3, 4, 5$.
\end{enumeratea}
\end{lem}

We will need the next two lemmas in the proof.
Recall the notation used in  \eqref{eqn:53}, and define
\begin{align}\label{eqn:999}
\widetilde h_n(x, y)
:=
\sum_{j=2}^n
\Big(
\frac{\lambda_{n;j}}{1-\lambda_{n;j}} 
\Big)
\psi_{n;j}(x)\psi_{n;j}(y)\, ,\ \ x, y\in[0, 1]\, .
\end{align}

\begin{lem}\label{lem:12a}
We have, 
$\ \sup_{n\geq 1}\medint\int_{[0, 1]^2} |\widetilde h_n|^3<\infty$.
\end{lem}

\begin{lem}\label{lem:12b}
\begin{inparaenuma}
\item\label{item:37}
For any $p\in [3/2, 3]$, $(I-T_n)$ is an invertible operator on $L^p([0, 1])$, where $I$ denotes the identity operator.
Further, there exists $C>0$ such that
$\|(I-T_n)^{-1}\|_{p, p}\leq C n^{\delta_0}$ for all $n$.

\medskip

\noindent\item\label{item:38}
We have,
$\| T_n \|_{3, \infty} = o(n^{\theta_0})$.
\end{inparaenuma}
\end{lem}

\noindent{\bf Proof of Lemma~\ref{lem:12a}:}
Define
\begin{align}\label{eqn:755}
h_n(x, y)
:=
\sum_{j=2}^n
\lambda_{n;j}
\psi_{n;j}(x)\psi_{n;j}(y)
=
\kappa_n(x, y) - \|T_n\| \psi_n(x)\psi_n(y)
\, , \ \ \ 
x, y\in[0, 1] .
\end{align}
Then note that
\begin{gather}
\medint{\int_{y=0}^{1}}h_n^2(x, y) dy\,
\lesssim\,
\medint{\int_{y=0}^{1}}\kappa_n^2(x, y) dy
+
\psi_n^2(x)
\leq
\medint{\int_{y=0}^{1}} W_n^2(x, y) dy
+
\psi_n^2(x)\, ,\label{eqn:70}
\end{gather}
and
\begin{align}\label{eqn:71}
\medint{\int_{y=0}^{1}}\widetilde h_n^2(x, y) dy\,
&=
\sum_{j=2}^n
\Big(
\frac{\lambda_{n;j}}{1-\lambda_{n;j}} 
\Big)^2
\psi_{n;j}^2(x)
\lesssim
\sum_{j=2}^n
\lambda_{n;j}^2 \psi_{n;j}^2(x)
\notag\\
&=
\medint{\int_{y=0}^{1}} h_n^2(x, y) dy
\leq
\medint{\int_{y=0}^{1}} W_n^2(x, y) dy
+
\psi_n^2(x)\, ,
\end{align}
where the second step uses  \eqref{eqn:765}, and the last step uses \eqref{eqn:70}.
Now,
\begin{align}\label{eqn:871}
\widetilde h_n(x, y)
=
h_n(x, y)
+
\medint\int_{z=0}^1 h_n(x, z)\, \widetilde h_n(z, y) dz
=
h_n(x, y)
+
\medint\int_{z=0}^1 \widetilde h_n(x, z)\, h_n(z, y) dz  \, ,
\end{align}
and consequently,
\begin{align}\label{eqn:541}
|\widetilde h_n(x, y)|
\leq
e^{-n}
+
W_n(x, y) + \psi_n(x)\psi_n(y)
+
\Big( \medint{\int_{z=0}^1} h_n^2(x, z)dz \Big)^{1/2}
\cdot
\Big( \medint{\int_{z=0}^1} \widetilde h_n^2(z, y)dz \Big)^{1/2} .
\end{align}
We thus have
\begin{align}
\medint\int_{[0, 1]^2} |\widetilde h_n|^3
&\lesssim
e^{-3n}
+
\medint\int_{[0, 1]^2} W_n^3
+
\big(\medint\int_0^1 \psi_n^3\big)^2
\notag\\
&
+
\left[
\int_{x=0}^1 \Big( \psi_n^2(x) + \medint{\int_{z=0}^1} W_n^2(x, z)dz\Big)^{3/2}dx
\right]
\cdot
\left[
\int_{y=0}^1 \Big( \psi_n^2(y) + \medint{\int_{z=0}^1} W_n^2(y, z)dz \Big)^{3/2} dy
\right]
\notag\\
&\hskip40pt
\lesssim
e^{-3n}
+
\medint\int_{[0, 1]^2} W_n^3
+
\big(\medint\int_0^1 \psi_n^3\big)^2
+
\Big[
\medint\int_0^1 \psi_n^3 
+
\medint{\int}_{[0, 1]^2} W_n^3
\Big]^2\, , \label{eqn:542}
\end{align}
where the first step uses 
\eqref{eqn:541}, \eqref{eqn:70}, and \eqref{eqn:71}
\chsen{together with the fact that
\begin{align}\label{eqn:0101}
| a_1+\ldots+a_k |^3 \leq k^2\cdot \big( |a_1|^3+\ldots+ |a_k|^3\big)
\end{align}
for any $a_1,\ldots, a_k\in\bR$.
}
By \eqref{eqn:22}, $(\|\psi_n\|_3\, ,\, n\geq 1)$ is a bounded sequence, and hence the desired result follows from \eqref{eqn:542} and \eqref{eqn:20}.
\qed

\medskip

\noindent{\bf Proof of Lemma~\ref{lem:12b}~\eqref{item:37}:}
Fix $p\in [3/2,\, 3]$, and regard all operators as acting on $L^p([0, 1])$ in this proof.
It follows from \eqref{eqn:755} and \eqref{eqn:871} that
\[
T_n
=
T_{h_n} + \|T_n\| \langle\psi_n,\, \cdot\ \rangle \psi_n
\ \ \ \text{ and }\ \ \ 
T_{\widetilde h_n} 
= 
T_{h_n} + T_{h_n}  T_{\widetilde h_n} 
=
T_{h_n} + T_{\widetilde h_n} T_{h_n}  \, .
\]
From this, it can be seen easily that $(I-T_n)$ is an invertible operator on $L^p([0, 1])$, and
\begin{align}\label{eqn:937}
(I-T_n)^{-1}
=
I+T_{\widetilde h_n}
+
\frac{\| T_n \|}{1-\|T_n\|} \langle\psi_n,\, \cdot\ \rangle \psi_n\, .
\end{align}
Let $p':=p/(p-1)$.
Since $p, p'\in [3/2, 3]$, it follows from \eqref{eqn:22} that
\begin{align}\label{eqn:231}
 \|\psi_n\|_p = O(1) \ \ \text{ and }\ \  \|\psi_n\|_{p'} = O(1)\, .
\end{align}
It thus follows from \eqref{eqn:937} that
\begin{align*}
\| (I-T_n)^{-1}\|_{p, p}
\leq
1 + 
\big(\medint\int_{[0, 1]^2} |\widetilde h_n|^3\big)^{1/3}
+
Cn^{\delta_0}  \|\psi_n\|_p\cdot \|\psi_n\|_{p'}
\leq
C'n^{\delta_0}\, ,
\end{align*}
where the first step uses Lemma~\ref{lem:2} and \eqref{eqn:655},
and the second step uses Lemma~\ref{lem:12a} and \eqref{eqn:231}.
This yields the claim in Lemma~\ref{lem:12b}~\eqref{item:37}.
\qed

\medskip

\noindent{\bf Proof of Lemma~\ref{lem:12b}~\eqref{item:38}:}
For any $f\in L^3([0, 1])$ and $x\in [0, 1]$,
\[
|T_n f(x)|
\leq
\medint{\int_0^1}\kappa_n(x, y)|f(y)| dy
\leq
\| f\|_3 \cdot \max_{x\in [0, 1]}\big(\medint{\int}_0^1 W_n^{3/2}(x, y) dy\big)^{2/3}\, .
\]
Now the desired result follows from Condition~\ref{ass:graphon-2}.
\qed

\medskip

\noindent{\bf Proof of Lemma~\ref{lem:12}~\eqref{item:11}:}
Let $\eps_{ij}$, $i, j\in [n]$, and $\fX_n'(j)$, $j\in [n]$, be independent random variables such that $\eps_{ij}\sim\text{Bernoulli}(\kappa_{ij}/n)$ and $\fX_n'(j) \equald \fX_n^{\sss\triangle}(j)$.
Then for each $i\in [n]$,
\begin{align}\label{eqn:737}
\fX_n^{\sss\triangle}(i) - 1
\equald
\sum_{j\in [n]}\eps_{ij}\cdot \fX_n'(j)\, .
\end{align}
Taking expectation on both sides yields
$
g_n(x) = \mvone_{\circ}(x) + T_n g_n (x)
$,
$x\in [0, 1]$.
Hence,
\begin{align}\label{eqn:123}
g_n = (I-T_n)^{-1}\mvone_{\circ}\, ,
\end{align}
and consequently,
$\|g_n\|_3 \leq \| (I-T_n)^{-1} \|_{\sss 3, 3} \cdot \| \, \mvone_{\circ}\|_3 \lesssim n^{\delta_0}$, 
where we have used Lemma~\ref{lem:12b}~\eqref{item:37}.
This completes the proof.
\qed

\medskip

\noindent{\bf Proof of Lemma~\ref{lem:12}~\eqref{item:12}:}
Squaring both sides of \eqref{eqn:737} and taking expectation yields
\begin{align*}
g_{n; 2}^{\sss\triangle}(i) - 2 g_n^{\sss\triangle}(i) + 1
=
\sum_{j\in [n]} \frac{\kappa_{ij}}{n} g_{n; 2}^{\sss\triangle}(j)
+
\Big(\sum_{j\in [n]} \frac{\kappa_{ij}}{n} g_n^{\sss\triangle}(j) \Big)^2
-
\sum_{j\in [n]} \frac{\kappa_{ij}^2}{n^2} \big(g_n^{\sss\triangle}(j)\big)^2\, ,
\end{align*}
for $i\in [n]$.
Hence, for $x\in [0, 1]$,
\begin{align}
g_{n; 2}(x)
&=
T_n g_{n; 2}(x) + 
\big(T_n g_n (x)\big)^2 + 2 g_n(x) - \mvone_{\circ}(x)
- 
\frac{1}{n}\medint{\int}_{y=0}^1 \kappa_n^2(x, y)g_n^2(y) dy
\label{eqn:659}\\
&\leq 
T_n g_{n; 2}(x) +
\big(T_n g_n (x)\big)^2 + 2 g_n(x)\, .\label{eqn:660}
\end{align}
Since $\| T_n\|<1$, $(I-T_n)^{-1}$ viewed as an operator on $L^2([0, 1])$ admits the expansion
\[
(I-T_n)^{-1}
=
I+T_n +T_n^2+\ldots\, .
\]
Thus, if $f_1\leq f_2$ are two functions in $L^2([0, 1])$, then 
$T_n^{\ell}f_1\leq T_n^{\ell}f_2$ for each $\ell\geq 1$ and consequently, 
$(I-T_n)^{-1}f_1\leq (I-T_n)^{-1}f_2$.
Hence, \eqref{eqn:660} implies that
\begin{align}\label{eqn:661}
g_{n; 2}
\leq
2(I-T_n)^{-1} g_n + (I-T_n)^{-1} \big[ ( T_n g_n )^2 \big]\, ,
\end{align}
which in turn yields
\begin{align}\label{eqn:662}
\| g_{n; 2} \|_{3/2}
&
\, \lesssim\,
\| (I-T_n)^{-1} \|_{\sss 3/2, 3/2} \cdot \| g_n \|_{3/2} 
+ 
\| (I-T_n)^{-1} \|_{\sss 3/2, 3/2} \cdot
\Big(
\medint\int_0^1 (T_n g_n)^3
\Big)^{2/3}
\notag\\
&
\, \lesssim\,
n^{2\delta_0}
+
n^{\delta_0}\cdot \| T_n g_n\|_3^2
\, \lesssim\,
n^{2\delta_0}
+
n^{\delta_0}\cdot \| T_n\|_{\sss 3, 3}^2 \| g_n\|_3^2
\, \lesssim\,
n^{3\delta_0}\, ,
\end{align}
where the second step uses Lemma~\ref{lem:12b}~\eqref{item:37} and Lemma~\ref{lem:12}~\eqref{item:11}, and the last step follows from \eqref{eqn:T-n-norm-bound} and Lemma~\ref{lem:12}~\eqref{item:11}.
Going back to \eqref{eqn:660}, we see that
\begin{align*}
&\| g_{n; 2}\|_3
\lesssim
\| T_n\|_{\sss 3/2, 3}\cdot \| g_{n; 2} \|_{3/2}
+
\| g_n\|_3
+
\big( \medint{\int}_0^1 (T_n g_n)^6 \big)^{1/3}\\
&\hskip50pt
\lesssim\,
n^{3\delta_0} + n^{\delta_0} + \| T_n g_n\|_{\infty} \cdot \| T_n g_n\|_3
\\
&\hskip90pt
\lesssim\,
n^{3\delta_0} + n^{\theta_0} \|g_n\|_3\cdot \| T_n\|_{3, 3} \cdot \| g_n\|_3
\lesssim\,
n^{3\delta_0} + n^{2\delta_0 + 1/12} 
= 
O(n^{3\delta_0})\, ,
\end{align*}
where the second step uses \eqref{eqn:T-n-norm-bound}, \eqref{eqn:662}, Lemma~\ref{lem:12}~\eqref{item:11},
the third step uses Lemma~\ref{lem:12b}~\eqref{item:38},
the fourth step uses \eqref{eqn:T-n-norm-bound}, Lemma~\ref{lem:12}~\eqref{item:11}, and the relation $\theta_0<1/12$ from Condition~\ref{ass:graphon-2}, and the last step uses the fact that $\delta_0>1/6$ from \eqref{eqn:17}.
This proves the first claim made in Lemma~\ref{lem:12}~\eqref{item:12}.

We now turn to the proof of \eqref{eqn:555}.
To this end, we examine the terms appearing on the right side of \eqref{eqn:659}.
It follows from \eqref{eqn:540} and \eqref{eqn:999} that
\begin{align}\label{eqn:998}
T_n g_n
=
\frac{\|T_n\|}{1-\|T_n\|} \big\langle\psi_n\, ,\, \mvone_{\circ} \big\rangle \psi_n
+
T_{\widetilde h_n}\mvone_{\circ}\, \, .
\end{align}
Hence,
\begin{gather}
(T_n g_n)^2
=
\frac{\|T_n\|^2}{(1-\|T_n\|)^2} \big(\medint\int_0^1\psi_n \big)^2\cdot \psi_n^2
+
a_n \, ,  \label{eqn:997} 
\end{gather}
where
\begin{align}\label{eqn:996}
\|a_n\|_{3/2}
\, \lesssim\,
n^{\delta_0} \big\| \psi_n\big\|_3 \cdot \big\| T_{\widetilde h_n}\mvone_{\circ} \big\|_3
+
\big\| \big( T_{\widetilde h_n}\mvone_{\circ} \big)^2 \big\|_{3/2}
\,\lesssim\,
n^{\delta_0}\, ;
\end{align}
here, we have used \eqref{eqn:655}, \eqref{eqn:22}, Lemma~\ref{lem:2}, and Lemma~\ref{lem:12a}.
Now, using \eqref{eqn:937}, we see that
\begin{align}\label{eqn:471}
(I-T_n)^{-1}\big[\psi_n^2\big]
=
\frac{\| T_n \|}{1-\|T_n\|} \big(\medint{\int}_0^1 \psi_n^3\big)\cdot\psi_n
+ 
b_n\, ,
\end{align}
where
\begin{align}\label{eqn:472}
\|b_n\|_{3/2}
\leq
\big\| I+T_{\widetilde h_n} \big\|_{\sss 3/2, 3/2} \cdot \big\| \psi_n^2\big\|_{3/2} 
\leq 
C\, .
\end{align}
Using Lemma~\ref{lem:12b}~\eqref{item:37} and \eqref{eqn:996}, we get
\begin{align}\label{eqn:473}
\big\|(I-T_n)^{-1}a_n \big\|_{3/2}
\, \lesssim\,
n^{2\delta_0}\, .
\end{align}
Combining \eqref{eqn:997},  \eqref{eqn:471}, \eqref{eqn:472}, and \eqref{eqn:473} 
with \eqref{eqn:22} and \eqref{eqn:655} yields
\begin{align}\label{eqn:474}
n^{-3\delta_0} 
\medint{\int}_0^1 (I-T_n)^{-1}\big[ (T_n g_n)^2\big]
\to
\big(\medint\int_0^1 \psi^3\big) \big(\medint\int_0^1 \psi \big)^{-3}  \ \ \text{ as }\ \ n\to\infty\, .
\end{align}
Next,
\begin{align}\label{eqn:475}
\bigg| \medint{\int}_0^1 (I-T_n)^{-1}\big[ 2g_n-1 \big] \bigg|
\leq
\big\|(I-T_n)^{-1}\big\|_{\sss 3/2, 3/2} \cdot \| 2g_n - 1\|_{3/2}
\, \lesssim\,
n^{2\delta_0}\, ,
\end{align}
where we have used  Lemma~\ref{lem:12}~\eqref{item:11} and Lemma~\ref{lem:12b}~\eqref{item:37}.
Finally,
\begin{align}\label{eqn:476}
\bigg\|
\frac{1}{n}\medint{\int}_{y=0}^1 \kappa_n^2(\, \cdot\, , y)g_n^2(y) dy
\bigg\|_{3/2}
\lesssim\,
\frac{1}{n^{1/3}} \big\| T_n g_n^2 \big\|_{3/2}
\lesssim\,
\frac{1}{n^{1/3}} \big\| g_n^2 \big\|_{3/2}
\lesssim
n^{2\delta_0-1/3}\, ,
\end{align}
where the first step uses \eqref{eqn:88-a}, the second step uses \eqref{eqn:T-n-norm-bound}, and the last step uses Lemma~\ref{lem:12}~\eqref{item:11}.
Combining \eqref{eqn:659} with \eqref{eqn:474}, \eqref{eqn:475}, and \eqref{eqn:476} yields the convergence in \eqref{eqn:555}.
\qed

\medskip

\noindent{\bf Proof of Lemma~\ref{lem:12}~\eqref{item:13}:}
The results of Lemma~\ref{lem:12}~\eqref{item:11} and \eqref{item:12} give the asymptotic behavior for $\| g_{n; k} \|_3$, $k=1,2$.
Now the asymptotics for $\| g_{n; k} \|_{3/2}$, $k=3, 4, 5$, can be obtained inductively.
We will only give the details for $k=5$, i.e., we will assume that
\begin{align}\label{eqn:477}
\| g_{n; k} \|_3 = O\big( n^{(2k-1)\delta_0} \big)\, ,\  k = 1, 2,
\ \ \text{ and }\ \ 
\| g_{n; k} \|_{3/2} = O\big( n^{(2k-1)\delta_0} \big)\, ,\ k=3, 4,
\end{align}
and from this we will deduce the desired result for $\| g_{n; 5} \|_{3/2}$.

Raising both sides of \eqref{eqn:737} to the power $5$ and taking expectation, and then using arguments similar to the one leading to \eqref{eqn:661}, we get
\begin{align}\label{eqn:434}
g_{n; 5}
\, \lesssim\,
(I-T_n)^{-1}
\Big[
&
(T_n g_{n; 3}) (T_n g_{n; 2}) + (T_n g_n)^5 + (T_n g_n)^3(T_n g_{n; 2}) 
\notag\\
&
+ (T_n g_{n;4}) (T_n g_n) + (T_n g_{n; 2})^2 (T_n g_n) + (T_n g_n)^2 (T_n g_{n;3}) + g_{n; 4}
\Big]
\, .
\end{align}
We will estimate the contribution coming from the first three terms on the right side of \eqref{eqn:434}; the other terms can be handled similarly.
We start with
\begin{align*}
\big\| (I-T_n)^{-1}
\big[
(T_n g_{n; 3}) (T_n g_{n; 2}) 
\big] \, \big\|_{3/2}\, 
&
\lesssim
n^{\delta_0}
\big\| 
(T_n g_{n; 3}) (T_n g_{n; 2}) 
\big\|_{3/2}
\leq
n^{\delta_0}
\| T_n g_{n; 3} \|_3 \cdot \| T_n g_{n; 2} \|_3
\\
&
\leq
n^{\delta_0}
\| T_n \|_{\sss 3/2, 3}^2 \cdot \| g_{n; 3} \|_{3/2} \cdot \| g_{n; 2} \|_{3/2}
\leq
C n^{9\delta_0} \, ,
\end{align*}
where the first step uses Lemma~\ref{lem:12b}~\eqref{item:37}, and the last step uses \eqref{eqn:477} and \eqref{eqn:T-n-norm-bound}.
Next,
\begin{align*}
\big\| (I-T_n)^{-1}
\big[ (T_n g_n)^5 \big] \, \big\|_{3/2}\, 
&
\lesssim
n^{\delta_0}	
\bigg( \medint\int_0^1 (T_n g_n)^{15/2}\bigg)^{2/3}
\leq
n^{\delta_0}	
\bigg(
\| T_n g_n\|_{\infty}^{9/2}
\medint\int_0^1 (T_n g_n)^3
\bigg)^{2/3}
\\
&
\lesssim
n^{\delta_0}	
\bigg(
n^{\theta_0} \| g_n\|_{3}
\bigg)^3
\big\| T_n g_n \big\|_3^2
\, \lesssim
n^{6\delta_0 + 3\theta_0}
= 
o(n^{9\delta_0}) \, ,
\end{align*}
where the third step uses Lemma~\ref{lem:12b}~\eqref{item:38}, and the last step uses the fact $\theta_0<1/12<1/6<\delta_0$.
Similarly,
\begin{align*}
&
\big\| (I-T_n)^{-1}
\big[ (T_n g_n)^3 (T_n g_{n; 2}) \big] \, \big\|_{3/2}\, 
\lesssim
n^{\delta_0}	
\big\| (T_n g_n)^3 (T_n g_{n; 2}) \big\|_{3/2}
\leq
n^{\delta_0}	
\big\| (T_n g_n)^3 \big\|_3
\cdot
\big\| T_n g_{n; 2}\big\|_3
\\[2pt]
&\hskip100pt
\leq
n^{\delta_0}\big\| T_n g_n \big\|_{\infty}^2\cdot \big\| T_n g_n \big\|_3\cdot \big\| T_n g_{n; 2}\big\|_3
\lesssim
n^{7\delta_0 + 2\theta_0}
=
o(n^{9\delta_0})\, .
\end{align*}
As mentioned before, the other terms on the right side of \eqref{eqn:434} can be bounded in a similar fashion. 
The claimed result follows upon combining all these bounds.
\qed

\subsection{The second susceptibility $s_2^{\bullet}$}\label{sec:g5}
Throughout this section, we work under Conditions~\ref{ass:graphon-1}, \ref{ass:graphon-2}, and \ref{ass:graphon-3a}.
Recall the constant $\zeta$ from \eqref{eqn:def-alpha-chi-zeta} and the susceptibility $s_2^{\bullet}$ from \eqref{eqn:898}.
Recall also that we are working under the assumptions of \eqref{eqn:765} and \eqref{eqn:17}.
The aim of this section is to prove the following result:

\begin{lem}\label{lem:20a}
	We have,
	$
	n^{1/3-2\delta_0} 
	\big(
	s_2^{\bullet}-n^{\delta_0} 
	\big) 
	\weakc
	\zeta\, 
	$,
	as $n\to\infty$.
\end{lem}

Lemma~\ref{lem:20a} follows upon combining Lemma~\ref{lem:20} with the next lemma and using the fact $\delta_0<1/3$.

\begin{lem}\label{lem:54}
The following hold:
\begin{enumeratea}
\item\label{item:54a}
$0\leq ( \int_0^1 g_n )- \bE[s_2^{\bullet}] \lesssim n^{4\delta_0 -1}$, and
\smallskip
\item\label{item:54b}
$\var(s_2^{\bullet}) \lesssim n^{5\delta_0 - 1}$.
\end{enumeratea}
\end{lem}

The rest of this section is devoted to the proof of Lemma~\ref{lem:54}.
For 
$x, y\in [0, 1]$ and $\ell\in\bZ_{>0}$, define
\begin{gather}
\rho_{n; \ell}(x, y)
:=
\| T_n\|^{\ell} \psi_n(x)\psi_n(y)
+
\sum_{j=2}^n \lambda_{n;j}^{\ell}\, \psi_{n; j}(x)\psi_{n; j}(y)\, , \ \ \text{ and} \label{eqn:323}\\
\rho_n(x, y)
:=
\frac{\| T_n \|}{1 - \| T_n \| } \psi_n(x)\psi_n(y)
+
\widetilde{h}_n(x, y)\, , \label{eqn:324}
\end{gather}
where $\widetilde{h}_n$ is as given in \eqref{eqn:999}.
Recall the notation $\connects$ and $\connectsell$ introduced below \eqref{eqn:898}.
We will need the following simple bounds on the connection probabilities in $\cG_n^{\bullet}$.

\begin{lem}\label{lem:55}
For $1\leq i\neq j\leq n$,
\[
\pr\big(i \connectsell j\big) \leq \frac{1}{n}\rho_{n; \ell}(x, y)
\ \ \ \text{ and }\ \ \
\pr\big(i \connects j\big) \leq \frac{1}{n}\rho_n(x, y)
\]
for any $x\in \big((i-1)/n\, ,\ i/n\big)$ and $y\in \big((j-1)/n\, ,\ j/n\big)$.
\end{lem}

\noindent{\bf Proof:}
Clearly,
\begin{align*}
&\pr\big( i \connectsell j \big)
\leq
\sum_{ k_1,\ldots, k_{\ell -1} }
\frac{1}{n^{\ell}} \cdot \kappa_{i k_1} \cdot \kappa_{k_1 k_2}\ldots \kappa_{k_{\ell -1} j} 
\\
&\hskip20pt
=
\frac{1}{n}\medint\int_{[0, 1]^{\ell-1}} \kappa_n(x, z_1)\kappa_n(z_1, z_2)\ldots\kappa_n(z_{\ell-1}, y)
\ dz_1 dz_2\ldots dz_{\ell -1}
=
\frac{1}{n}
\rho_{n; \ell}(x, y)\, .
\end{align*}
The bound on $\pr\big(i \connects j\big)$ follows if we sum over all $\ell\geq 1$.
\qed

\medskip

\noindent{\bf Proof of Lemma~\ref{lem:54}~\eqref{item:54a}:}
Consider $i\in [n]$, and explore its component in $\cG_n^{\bullet}$, namely $\cC^{\bullet}(i)$, in a breadth-first fashion to build its breadth-first spanning tree (BFS tree).
The exploration starts at $i$ which is designated as the root of the BFS tree.
In each step of this exploration, one looks for edges between the current vertex being explored, say $j$, and all vertices that have not been discovered up to that point.
The vertices found are added as children of $j$ in the BFS tree and they are ordered from left to right using their labels.
The vertices in each generation of the BFS tree are explored from left to right, and the exploration process moves to the next generation when all vertices in a generation have been explored.
There is an obvious way of coupling this BFS tree with a multitype branching process started from a vertex of type $i$ as described above \eqref{eqn:52}:
For each $j$ in the BFS tree and each $k$ discovered before the exploration of the vertex $j$ (this includes the vertices $j$ and $i$), we give $j$ an extra child of type $k$ with probability $\kappa_{jk}/n$, and we do this independently for each such pair $(j, k)$ of vertices.
Then we run independent branching processes starting from each new vertex added to the BFS tree.
This will yield a multitype branching process tree that contains the BFS tree as a subtree.
In particular, 
\begin{align}\label{eqn:coupling}
\diam\big( \cC^{\bullet}(i) \big)
\leq
2\fH_n^{\sss\triangle}(i)
\ \ \text{ and }\ \
|\cC^{\bullet}(i)| \leq \fX_n^{\sss\triangle}(i)  \ \ \text{ in this coupling.}
\end{align}

Let $v$ be uniformly distributed on $[n]$.
Note that $\bE[s_2^{\bullet}] = \bE\big( |\cC^{\bullet}(v) | \big)$ and 
$\bE\big[\fX_n^{\sss\triangle}(v)\big] = \int_0^1 g_n$.
It then follows from the above coupling that
\begin{align}\label{eqn:72}
	0
	\leq   
	\medint\int_0^1 g_n - \bE[s_2^{\bullet}] 
	&
	\leq 
	\frac{1}{n}\left(
	\sum_{1\leq i\neq j\leq n} 
	\pr\big(i \connects j\big)\cdot
	\Big(
	\frac{\kappa_{ij}}{n}g_n^{\sss\triangle}(i)
	+
	\frac{\kappa_{jj}}{n}g_n^{\sss\triangle}(j)
	\Big)
	\right.
	\notag\\
	&\hskip30pt
	+
	\left.
	\sum_{i, j, k\text{ (p.d.)}} 
	\pr\bigg(i\connects j\, ,\ i\connects k\, ,\ k\text{ before }j\bigg)
	\frac{\kappa_{jk}}{n}g_n^{\sss\triangle}(k)
	\right)\, ,
\end{align}
where `$i, j, k$ (p.d.)' stands for the sum over all $i, j, k\in [n]$ that are pairwise distinct, and `$k$ before $j$' represents the event that $k$ is discovered before $j$ is explored in the breadth-first exploration.

\begin{figure}[t]
\centering
\includegraphics[trim=1cm 16.3cm 1.7cm 5cm, clip=true, angle=0, scale=.82]{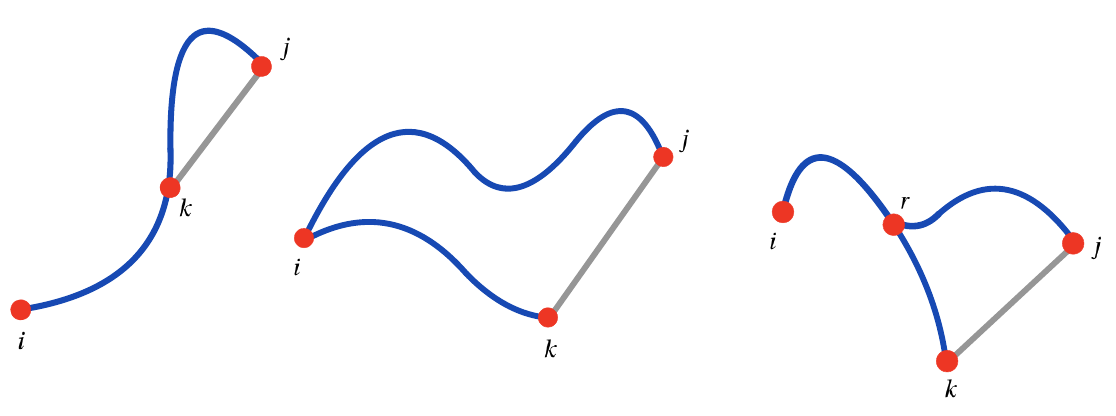}
\caption{
The blue curves represent paths of positive length in $\cG_n^{\bullet}$. 
The gray lines indicate that $j$ gets an extra child of type $k$ with probability $\kappa_{jk}/n$.}
\label{fig:1}
\end{figure}

We first bound the second term appearing on the right side of \eqref{eqn:72}.
\chsen{Recall the notation `$\circ$' signifying disjoint occurrence of events from Section~\ref{sec:bk-inequality}, and} 
note that
\begin{align}
	\big\{ i\connects j\, ,\ i\connects k\, ,\ k\text{ before }j \big\}
	&\subseteq
	\bigg( \bigcup_{r\neq i, j, k} 
	\big\{ i\connects r\big\}\circ 
	\big\{ r\connects k\big\}\circ 
	\big\{ r\connects j\big\}\bigg)
	\notag
	\\
	&\hskip10pt
	\bigcup\
	\bigg( \big\{ i\connects k\big\} \circ \big\{ k\connects j \big\} \bigg)
	\ \bigcup\
	\bigg( \big\{ i\connects k\big\} \circ \big\{ i\connects j \big\} \bigg)\, ;
	\label{eqn:278}
\end{align}
these events are depicted in Figure~\ref{fig:1}.
\chsen{Now, the right side of \eqref{eqn:278} involves only connection events pertaining to the random graph $\cG_n^{\bullet}$.
Thus, by Lemma~\ref{lem:bk-inequality},}
\begin{align*}
\pr\big(i\connects j\, ,\ i\connects k\, ,\ k\text{ before }j \big)
&
\leq	
\pr\big(i\connects k\big) \pr\big(k\connects j\big)
+
\pr\big(i\connects k\big) \pr\big(i\connects j\big)
\notag\\
&\hskip40pt
+
\sum_{r\neq i, j, k} \pr\big(i\connects r\big) \pr\big(r\connects j\big) \pr\big(r\connects k\big)\, .
\end{align*}
Combining this bound with Lemma~\ref{lem:55}, we see that
\begin{align}\label{eqn:112}
\sum_{i, j, k\text{ (p.d.)}} 
\pr\bigg(i\connects j\, ,\ i\connects k\, ,\ k\text{ before }j\bigg)
\frac{\kappa_{jk}}{n}g_n^{\sss\triangle}(k)
\leq
F_1 + F_2 +F_3\, ,
\end{align}
where
\begin{gather*}
F_1 : =
\medint\int_{[0, 1]^3} \rho_n(x, z) \rho_n(y, z) \kappa_n(y, z) g_n(z) \ dx dy dz
\, ,\\
F_2 : =
\medint\int_{[0, 1]^3} \rho_n(x, z) \rho_n(x, y) \kappa_n(y, z) g_n(z) \ dx dy dz
\, , \ \text{ and}\\
F_3 : =
\medint\int_{[0, 1]^4} \rho_n(x, w) \rho_n(w, y) \rho_n(w, z) \kappa_n(y, z) g_n(z) \ dx dy dz dw\, .
\end{gather*}
In order to bound $F_1$, recall from \eqref{eqn:123} that $g_n = \mvone_{\circ} + T_n g_n$.
Hence,
\begin{align}\label{eqn:124}
\| g_n \|_{\infty}
\leq 
1+\| T_n \|_{3, \infty}\cdot \| g_n \|_3
\, \lesssim\, 
n^{\theta_0 + \delta_0}\, ,
\end{align}
where the last step uses Lemma~\ref{lem:12}~\eqref{item:11} and Lemma~\ref{lem:12b}~\eqref{item:38}.
Further, from the expression for $\rho_n$ given in \eqref{eqn:324}, we see that
\begin{align}\label{eqn:125}
\medint{\int}_{[0, 1]^3}\rho_n^3
\, \lesssim\,
n^{3\delta_0} \big(\medint{\int}_0^1\psi_n^3\big)^2 + \medint{\int}_{[0, 1]^3} |\widetilde h_n|^3
\lesssim\,
n^{3\delta_0}\, ,
\end{align}
where the first step uses \eqref{eqn:655}, 
and the second step uses Lemma~\ref{lem:12a} and \eqref{eqn:22}.
Hence,
\begin{align}\label{eqn:126}
F_1
\lesssim
n^{\theta_0 + \delta_0}
\medint\int_{[0, 1]^3} \rho_n(x, z) \rho_n(y, z) \kappa_n(y, z) \ dx dy dz
\, \lesssim\,
n^{\delta_0 + \theta_0}\cdot n^{2\delta_0}
=
o(n^{4\delta_0})\, ,
\end{align}
where the first step uses \eqref{eqn:124}, the second step uses H\"{o}lder's inequality in conjuction with \eqref{eqn:125} and the fact $\int_{[0, 1]^2}\kappa_n^3 = O(1)$, and the last step uses the relation $\theta_0<1/12<1/6<\delta_0$.

Turning to $F_3$, recall the function $h_n$ from \eqref{eqn:755} and note that
\begin{align}\label{eqn:756}
\medint{\int}_{[0, 1]^2} |h_n|^3
\lesssim
\medint{\int}_{[0, 1]^2} \kappa_n^3
+
\big(\medint{\int}_0^1\psi_n^3\big)^2
=
O(1)\, ,
\end{align}
\chsen{where the first step uses \eqref{eqn:0101} with $k=2$.}
Now,
\begin{align*}
\widetilde\rho_n(z, w)
&
:=
\medint\int_{y=0}^1 \kappa_n(z, y)\rho_n(w, y)\ dy
=
\frac{\| T_n \|^2}{1 - \| T_n \| } \psi_n(w)\psi_n(z)
+
\sum_{j=2}^n
\frac{\lambda_{n;j}^2}{1 - \lambda_{n;j} } \psi_{n; j}(w)\psi_{n;j}(z)
\\
&\hskip50pt
=
\frac{\| T_n \|^2}{1 - \| T_n \| } \psi_n(w)\psi_n(z) 
+ \widetilde h_n(w, z)  - h_n(w, z)\, ,
\end{align*}
and consequently,
\begin{align}\label{eqn:757}
\medint{\int}_{[0, 1]^2} \widetilde\rho_n^3
\lesssim
n^{3\delta_0} \big(\medint{\int}_0^1\psi_n^3\big)^2 
+ 
\medint{\int}_{[0, 1]^3} |\widetilde h_n|^3
+
\medint{\int}_{[0, 1]^3} | h_n|^3
\lesssim
n^{3\delta_0}\, .
\end{align}
Hence,
\begin{align}\label{eqn:758}
F_3
&
=
\medint{\int}_{[0, 1]^3}
\big(\rho_n(x, w) g_n(z) \big)\cdot \widetilde\rho_n(w, z)\cdot \rho_n(w, z) \ dx dz dw
\notag\\
&
\leq
\Big(\medint{\int}_{[0, 1]^2}\rho_n^3 \cdot \medint\int_0^1 g_n^3\Big)^{1/3}\cdot
\Big(\medint{\int}_{[0, 1]^2}\widetilde\rho_n^3\Big)^{1/3}\cdot
\Big(\medint{\int}_{[0, 1]^2}\rho_n^3\Big)^{1/3}
\lesssim
n^{4\delta_0}\, ,
\end{align}
where we have used \eqref{eqn:125}, \eqref{eqn:757}, and Lemma~\ref{lem:12}~\eqref{item:11}.
Similarly,
\begin{align}\label{eqn:759}
F_2
=
\medint{\int}_{[0, 1]^2}
\rho_n(x, z)  \cdot \widetilde\rho_n(x, z)\cdot g_n(z) \ dx dz
\,\lesssim\,
n^{3\delta_0} \, .
\end{align}
Combining \eqref{eqn:112}, \eqref{eqn:126}, \eqref{eqn:758}, and \eqref{eqn:759} yields
\begin{align}\label{eqn:112-a}
\sum_{i, j, k\text{ (p.d.)}} 
\pr\bigg(i\connects j\, ,\ i\connects k\, ,\ k\text{ before }j\bigg)
\frac{\kappa_{jk}}{n}g_n^{\sss\triangle}(k)
\lesssim
n^{4\delta_0}\, .
\end{align}

Using similar arguments and the fact that $\kappa_{jj}=e^{-n}$, we can show that
\[
\sum_{1\leq i\neq j\leq n} 
\pr\big(i \connects j\big)\cdot
\Big(
\frac{\kappa_{ij}}{n}g_n^{\sss\triangle}(i)
+
\frac{\kappa_{jj}}{n}g_n^{\sss\triangle}(j)
\Big)
\lesssim
n^{2\delta_0}\, ,
\]
which combined with \eqref{eqn:112-a} and \eqref{eqn:72} completes the proof.
\qed

\medskip

\noindent{\bf Proof of Lemma~\ref{lem:54}~\eqref{item:54b}:}
Let $u$ and $v$ be uniformly distributed on $[n]$ such that $u, v$, and $\cG_n^{\bullet}$ are independent.
Note that for any $k\geq 1$,
\begin{align}\label{eqn:768}
s_{k+1}^{\bullet} 
= 
\E\big( |\cC^{\bullet}(v)|^k\, \big|\, \cG_n^{\bullet}\big)
\ \ \text{ and }\ \
\bE\big[\big(s_{k+1}^{\bullet}\big)^2\big]
=
 \E\big( |\cC^{\bullet}(v)|^k\cdot  |\cC^{\bullet}(u)|^k\big)\, .
\end{align}
For any nonempty $A\subseteq [n]$, write $\cG^{A-}$ for the graph induced by $\cG_n^{\bullet}$ on the vertex set $[n]\setminus A$. 
Then
\begin{align*}
\E\Big( |\cC^{\bullet}(u)|^k \big| V\big(\cC^{\bullet}(v)\big)=A \Big)
= 
\frac{|\cC^{\bullet}(v)|^{k+1}}{n} 
+ 
\frac{1}{n} \E\Big(\sum_{\substack{\cC\text{ compo-}\\ \text{nent of }\cG^{A-}}} |\cC|^{k+1} \Big)
\leq 
\frac{1}{n}|\cC^{\bullet}(v)|^{k+1} 
+ 
\E(s_{k+1}^{\bullet})\, .
\end{align*}
Multiplying both sides by $|\cC^{\bullet}(v)|^k$, taking expectation, and using \eqref{eqn:768}, we get
$
\bE\big[\big(s_{k+1}^{\bullet}\big)^2\big]
\leq
n^{-1}\bE\big( |\cC^{\bullet}(v)|^{2k+1} \big)
+
\big(
\bE s_{k+1}^{\bullet}
\big)^2
$.
This leads to the following useful bound:
\begin{align}\label{eqn:769}
\var\big( s_{k+1}^{\bullet} \big)
\leq 
n^{-1}\bE\big( |\cC^{\bullet}(v)|^{2k+1} \big)
\leq
n^{-1}\medint{\int}_0^1 g_{n; 2k+1} 
\leq
n^{-1} \cdot\| g_{n; 2k+1} \|_{3/2}\, ,
\end{align}
where the second step follows from \eqref{eqn:coupling}.
Now an application of Lemma~\ref{lem:12}~\eqref{item:13} shows that
$\var\big( s_2^{\bullet} \big)
\leq
n^{-1} \cdot\| g_{n; 3} \|_{3/2}
\lesssim
n^{5\delta_0 -1}
$,
which completes the proof.
\qed

\subsection{Tail bounds for $|\cC_1^{\bullet}|$ and $D_{\max}^{\bullet}$}\label{sec:g7}
Throughout this section, we work under Conditions~\ref{ass:graphon-1}, \ref{ass:graphon-2}, and \ref{ass:graphon-3a}.
We will further assume that \eqref{eqn:765} and \eqref{eqn:17} hold.
The following lemma is the main result of this section:

\begin{lem}\label{lem:43}
For each $r\geq 1$, there exist $C_1, C_2>0$ such that for all $n$,
\begin{gather}
\pr\big(D_{\max}^{\bullet} \geq C_1 n^{\delta_0}\log n\big)\leq C_2 n^{-r}\ \ \ \ \text{ and}
\label{eqn:981}\\
\pr\big(|\cC_1^{\bullet}| \geq C_1 n^{2\delta_0 + \theta_0}\log n\big)\leq C_2 n^{-r}\, .
\label{eqn:982}
\end{gather}
\end{lem}

\noindent{\bf Proof:}
It follows from \eqref{eqn:coupling} that for any $\ell\geq 1$,
\begin{align}\label{eqn:215}
\pr\big(D_{\max}^{\bullet} \geq 2\ell\big)
&\leq
\sum_{i\in [n]}\pr\big(\fH_n^{\sss\triangle}(i)\geq \ell\big)
=
n\pr\big(\fH_n^{\sss\triangle}(v)\geq \ell\big)\, ,
\end{align}
where $v$ is uniformly distributed on $[n]$.
Now,
\begin{align}\label{eqn:177}
\pr\big(\fH_n^{\sss\triangle}(v)\geq \ell\big)
=
\pr\big( |\fG_{\ell}^{\sss\triangle}(v)|\geq 1\big)
\leq
\bE\big(|\fG_{\ell}^{\sss\triangle}(v)|\big)
=
\medint{\int}_0^1 T_n^{\ell}\mvone_{\circ}
\leq
\| T_n^{\ell}\mvone_{\circ}\|_2
\leq
\| T_n\|^{\ell}  ,
\end{align}
where the third step uses \eqref{eqn:105}.
By \eqref{eqn:655} and \eqref{eqn:765}, 
\begin{align}\label{eqn:16}
\| T_n \| \leq 1-C_{\ref{eqn:16}} n^{-\delta_0}\ \ \text{ for all }\ \ n\, ,
\end{align}
which combined \chsen{with} \eqref{eqn:177} and \eqref{eqn:215} yields the claim in  \eqref{eqn:981}.

We now turn to the proof of \eqref{eqn:982}.
Note that for any $\ell\geq 2$ and any $n$,
\begin{align}\label{eqn:18}
\| T_n^{\ell}\mvone_{\circ} \|_{\infty}
\leq 
C n^{\theta_0} \| T_n^{\ell - 1}\mvone_{\circ} \|_3
\leq 
C' n^{\theta_0} \| T_n^{\ell - 2}\mvone_{\circ} \|_2
\leq
C' n^{\theta_0} \| T_n\|^{\ell - 2} \, ,
\end{align}
where the first step uses Lemma~\ref{lem:12b}~\eqref{item:38}, and the second step uses \eqref{eqn:T-n-norm-bound}.
Define
\begin{align*}
\Phi_{n, k} : = 
\mvone_{\circ} 
+ 
\sum_{\ell=1}^{k}
\big( 1+\fb_n\big)^{\ell} T_n^{\ell}\mvone_{\circ}
\ \ \text{ and }\ \
\Phi_n : = 
\mvone_{\circ} 
+ 
\sum_{\ell=1}^{\infty}
\big( 1+\fb_n\big)^{\ell} T_n^{\ell}\mvone_{\circ}\, ,
\ \ \text{ where }\ \
\fb_n
=
\frac{C_{\ref{eqn:16}}}{2 n^{\delta_0}}\, .
\end{align*}
It follows from \eqref{eqn:16} and \eqref{eqn:18} that
\begin{align}\label{eqn:24}
\| \Phi_n\|_{\infty}
\leq
1 + 
\Big(1+\frac{C_{\ref{eqn:16}}}{2 n^{\delta_0}}\Big) \| T_n\mvone_{\circ}\|_{\infty}
+
C' n^{\theta_0}
\sum_{\ell = 2}^{\infty}
\Big(1+\frac{C_{\ref{eqn:16}}}{2 n^{\delta_0}}\Big)^{\ell}
\Big(1-\frac{C_{\ref{eqn:16}}}{n^{\delta_0}}\Big)^{\ell - 2}
\leq
C_{\ref{eqn:24}}n^{\theta_0 + \delta_0}\, .
\end{align}
Let
\begin{align}\label{eqn:25}
\eps_n 
:=
C_{\ref{eqn:25}} n^{-\theta_0 - 2\delta_0}\, ,
\ \ \text{ where }\ \
C_{\ref{eqn:25}}
=
\frac{C_{\ref{eqn:16}}}{2 e C_{\ref{eqn:24}}}\, .
\end{align}
Using \eqref{eqn:24}, we can choose $n_0\geq 1$ be such that 
$\eps_n \| \Phi_n\|_{\infty} \leq 1$ for all $n\geq n_0$.
Now we claim that for any $n\geq n_0, r\in [n]$, and $k\geq 1$,
\begin{align}\label{eqn:26}
\bE \exp\bigg(\eps_n\sum_{\ell = 0}^{k} |\fG_{\ell}^{\sss\triangle}(r)|\bigg)
\leq
\exp\bigg(\eps_n\Phi_{n, k} \big( \tfrac{r}{n} - \tfrac{1}{2n}\big)\bigg)
\, .
\end{align}
From this claim it follows, by letting $k\to\infty$, that for any $n\geq n_0$ and $r\in [n]$,
\begin{align*}
\bE \exp\big( \eps_n \fX_n^{\sss\triangle}(r)\big)
\leq
\exp\bigg(\eps_n\Phi_n \big( \tfrac{r}{n} - \tfrac{1}{2n}\big)\bigg)\, .
\end{align*}
This in turn implies that for any $C>0, n\geq n_0$, and $r\in [n]$,
\begin{align}\label{eqn:7}
\pr\big( \fX_n^{\sss\triangle}(r)\geq C(\log n)/\eps_n\big)
\leq
\exp(-C\log n)
\cdot
\exp\big(\eps_n \| \Phi_n \|_{\infty}\big)
\leq
\exp\big( -C\log n + 1 \big)\, ,
\end{align}
which combined with \eqref{eqn:coupling} and a union bound yields the claim in \eqref{eqn:982}.

It remains to prove \eqref{eqn:26}.
Fix $k\geq 1, n\geq n_0$, and $r\in [n]$.
Let us make note of the following simple inequality: $e^x - 1 \leq x (1+ xe)$ for all $x\in [0, 1]$.
Hence, for $x\in [0, 1]$ and $m \in \bZ_{\geq 0}$,
\begin{align}\label{eqn:28}
\exp\big(\eps_n\Phi_{n, m}(x)\big) - 1
&
\leq
\eps_n\Phi_{n, m}(x) \cdot \big( 1 + \eps_n\Phi_{n, m}(x) e \big)
\notag\\
&
\leq
\eps_n\Phi_{n, m}(x) \cdot \big( 1 + \eps_n\| \Phi_n\|_{\infty} e \big)
\leq
\eps_n\Phi_{n, m}(x)\cdot ( 1 + \fb_n )\, ,
\end{align}
where the last step uses \eqref{eqn:24} and \eqref{eqn:25}.
Consider a multitype branching process starting from a vertex of type $r\in [n]$ as described above \eqref{eqn:52}, and for any vertex $w$ in this branching process tree, let $i_w\in [n]$ denote its type, and let 
$x_w := n^{-1}(i_w - 1/2)$;
we will similarly write $x_j := n^{-1}(j - 1/2)$ for $j\in [n]$.
We will prove inductively that for $t=k-1, k-2, \ldots, 1, 0$,
on the event $\{ \fG_t^{\sss\triangle}(r) \neq \emptyset \}$,
\begin{align}\label{eqn:290}
\bE \bigg(
\exp\bigg(\eps_n\sum_{m = t }^{k} |\fG_m^{\sss\triangle} (r)|\bigg)
\, \bigg|\,
\fG_{t}^{\sss\triangle}(r)
\bigg)
\leq
\prod_{w\in\fG_{t}^{\sss\triangle}(r)}
\exp\bigg(\eps_n\Phi_{n, k-t} (x_w)\bigg)
\, ;
\end{align}
the inequality in \eqref{eqn:26} corresponds to $t = 0$.
Starting with $t = k-1$, we see that on the event $\{ \fG_{k-1}^{\sss\triangle}(r) \neq \emptyset \}$,
\begin{align*}
&
\bE \bigg(
\exp\big(\eps_n|\fG_k^{\sss\triangle} (r)|\big)
\, \bigg|\,
\fG_{k-1}^{\sss\triangle}(r)
\bigg)
=
\prod_{w\in\fG_{k-1}^{\sss\triangle}(r)}
\prod_{j\in [n]}
\bigg( 1 - \frac{\kappa_{i_w j}}{n} +  \frac{\kappa_{i_w j}}{n} e^{\eps_n} \bigg)
\\
&\hskip30pt
\leq
\prod_{w\in\fG_{k-1}^{\sss\triangle}(r)}
\prod_{j\in [n]}
\exp\bigg( \kappa_{i_w j} \big( e^{\eps_n} - 1\big) /n \bigg)
\leq
\prod_{w\in\fG_{k-1}^{\sss\triangle}(r)}
\prod_{j\in [n]}
\exp\bigg( \eps_n (1 + \fb_n) \kappa_{i_w j} /n \bigg)
\\
&\hskip60pt
=
\prod_{w\in\fG_{k-1}^{\sss\triangle}(r)}
\exp\bigg( \eps_n (1 + \fb_n) T_n \mvone_{\circ} (x_{w}) \bigg)
\, ,
\end{align*}
where the third step uses \eqref{eqn:28} with $m=0$.
This implies \eqref{eqn:290} for $t = k-1$.
Suppose \eqref{eqn:290} holds for $t = k-1, k-2, \ldots, \ell +1$ for some $\ell\geq 0$.
Then on the event $\{ \fG_{\ell}^{\sss\triangle}(r) \neq \emptyset \}$,
\begin{align*}
&
\bE \bigg(
\exp\bigg(\eps_n\sum_{m = \ell + 1 }^{k} |\fG_m^{\sss\triangle} (r)|\bigg)
\, \bigg|\,
\fG_{\ell}^{\sss\triangle}(r)
\bigg)
\leq
\bE \bigg( \prod_{\widetilde w\in\fG_{\ell + 1}^{\sss\triangle} (r) }
\exp\big( \eps_n \Phi_{n, k-\ell-1}(x_{\widetilde w})  \big)
\, \bigg|\,
\fG_{\ell}^{\sss\triangle}(r)
\bigg)
\\
&\hskip10pt
=
\prod_{w\in\fG_{\ell}^{\sss\triangle} (r)} \prod_{j \in [n]}
\bigg( 
1 - \frac{\kappa_{i_w j}}{n} +  \frac{\kappa_{i_w j}}{n} \exp\big(\eps_n \Phi_{n, k-\ell-1}(x_j ) \big) 
\bigg)
\\
&\hskip20pt
\leq
\prod_{w\in\fG_{\ell}^{\sss\triangle} (r)} \prod_{j \in [n]}
\exp\bigg(
\frac{\kappa_{i_w j}}{n} \bigg[ \exp\big(\eps_n \Phi_{n, k-\ell-1}(x_j ) \big) - 1 \bigg]
\bigg)
\\
&\hskip30pt
\leq
\prod_{w\in\fG_{\ell}^{\sss\triangle} (r)} \prod_{j \in [n]}
\exp\bigg(
\frac{\kappa_{i_w j}}{n} \bigg[ \eps_n\Phi_{n, k-\ell-1}(x_j) \cdot (1+\fb_n) \bigg]
\bigg)\\
&\hskip40pt
=
\prod_{w\in\fG_{\ell}^{\sss\triangle} (r)} 
\exp\bigg( \eps_n (1 + \fb_n) T_n\Phi_{n, k-\ell-1}( x_w ) 
\bigg)
=
\prod_{w\in\fG_{\ell}^{\sss\triangle} (r)} 
\exp\bigg( \eps_n \big(\Phi_{n, k-\ell}( x_w ) - 1 \big)
\bigg)\, ,
\end{align*}
where the fourth step uses \eqref{eqn:28}.
This completes the proof.
\qed

\subsection{The third susceptibility $s_3^{\bullet}$}\label{sec:g8}
Throughout this section, we work under Conditions~\ref{ass:graphon-1}, \ref{ass:graphon-2}, and \ref{ass:graphon-3a}.
We will further assume that \eqref{eqn:765} and \eqref{eqn:17} hold.
Our goal in this section is to prove the following lemma:

\begin{lem}\label{lem:5}
We have,
$
n^{-3\delta_0}\cdot s_3^{\bullet}
\weakc
\int_0^1\psi^3 \cdot \big( \int_0^1 \psi \big)^{-3}
$
as $n\to\infty$.
\end{lem}

Note that
\[
\var(s_3^{\bullet} )
\leq 
n^{-1}\| g_{n; 5} \|_{3/2}
\lesssim
n^{9\delta_0 - 1}
=
o(n^{6\delta_0})\, ,
\]
where the first step follows from \eqref{eqn:769}, the second step uses Lemma~\ref{lem:12}~\eqref{item:13}, and the last step uses the relation $\delta_0<1/3$.
Thus, the claim in Lemma~\ref{lem:5} follows from the next result and \eqref{eqn:555}.

\begin{lem}\label{lem:6}
We have,
$
0\leq 
\int_0^1 g_{n; 2} - \bE[s_3^{\bullet}] 
=
o(n^{3\delta_0})
$.
\end{lem}

\noindent{\bf Proof:}
Let $v$ be uniformly distributed over $[n]$.
Then
$\bE[s_3^{\bullet}] = \bE\big(|\cC^{\bullet} (v)|^2 \big)$.
Suppose $\cC^{\bullet}(v)$ is coupled with a multitype branching process started from a vertex of type $v$ as described above \eqref{eqn:coupling}.
In this coupling, $\fX_n^{\sss\triangle}(v)\geq |\cC^{\bullet}(v)|$ and consequently,
\[
\medint\int_0^1 g_{n; 2} - \bE[s_3^{\bullet}] 
=
\bE\big(\fX_n^{\sss\triangle}(v)^2\big) - \bE\big( |\cC^{\bullet} (v)|^2 \big)
=
\bE\big( \fX_n^{\sss\triangle}(v)^2 - |\cC^{\bullet} (v)|^2 \big)
\geq 
0\, .
\]
Further, for any $m, n\geq 1$,
\begin{align*}
&
\bE\big( \fX_n^{\sss\triangle}(v)^2 - |\cC^{\bullet} (v)|^2 \big)
=
\bE\big[ \big( \fX_n^{\sss\triangle}(v) + |\cC^{\bullet} (v)|\big) \cdot \big( \fX_n^{\sss\triangle}(v) - |\cC^{\bullet} (v)|\big) \big]
\\
&\hskip25pt
\leq
2\cdot\bE\big[ \fX_n^{\sss\triangle}(v) \cdot \big( \fX_n^{\sss\triangle}(v) - |\cC^{\bullet} (v)|\big)  \big]
\\
&\hskip50pt
\leq 
2m\cdot \bE\big( \fX_n^{\sss\triangle}(v) - |\cC^{\bullet} (v)|\big)  
+
2\cdot\bE\bigg( \fX_n^{\sss\triangle}(v)^2 \ind_{ \big\{ \fX_n^{\sss\triangle}(v) \geq m \big\}} \bigg)  
\\
&\hskip75pt
\leq
2m\cdot \bigg( \medint\int_0^1 g_n - \bE[s_2^{\bullet}] \bigg)
+
2\big(\medint\int_0^1 g_{n;4}\big)^{1/2} \pr\big(\fX_n^{\sss\triangle}(v) \geq m\big)^{1/2}
\\
&\hskip100pt
\leq
C\Big( m n^{4\delta_0 - 1} + n^{7\delta_0/2} \cdot \pr\big(\fX_n^{\sss\triangle}(v) \geq m\big)^{1/2}\Big)
\, ,
\end{align*}
where the last step uses  Lemma~\ref{lem:54}~\eqref{item:54a} and Lemma~\ref{lem:12}~\eqref{item:13}.
Using \eqref{eqn:7}, we conclude that
\[
0\leq
\medint\int_0^1 g_{n; 2} - \bE[s_3^{\bullet}] 
\, \lesssim\, 
n^{\theta_0 + 6\delta_0 - 1} \cdot \log n
=
o(n^{3\delta_0})\, ,
\]
since $3\delta_0 < 1-\theta_0$.
This completes the proof.
\qed

\subsection{Asymptotics for ${\cD}^{\bullet}$}\label{sec:g9}
Throughout this section, we work under Conditions~\ref{ass:graphon-1}, \ref{ass:graphon-2}, and \ref{ass:graphon-3a}.
Further, we assume that \eqref{eqn:765} and \eqref{eqn:17} hold.
The following lemma is the main result of this section.

\begin{lem}\label{lem:66}
We have, 
$n^{-2\delta_0} \cdot\cD^{\bullet}\weakc \big(\int_0^1\psi\big)^{-2}$
as $n\to\infty$.
\end{lem}

This result follows immediately from the next lemma and the fact that $\delta_0<1/3$.

\begin{lem}\label{lem:67}
The following hold:
\begin{enumeratea}
\item\label{item:aa}
$n^{-2\delta_0} \cdot \bE\big[\int_0^1\fZ_n \big]\to \big(\int_0^1\psi\big)^{-2}$ as $n\to\infty$;
\medskip
\item\label{item:bb}
$0 \leq \bE\big[\int_0^1\fZ_n \big] - \bE\big(\cD^{\bullet}\big) \lesssim n^{5\delta_0 -1}\cdot\log n\,$; and
\medskip
\item\label{item:cc}
$\var\big( \cD^{\bullet} \big) \lesssim n^{7\delta_0 -1}(\log n)^2$.
\end{enumeratea}
\end{lem}

\noindent{\bf Proof of Lemma~\ref{lem:67}~\eqref{item:aa}:}
Note that
\begin{align*}
\bE\big[\medint\int_0^1\fZ_n \big]
=
\medint{\int}_0^1 \sum_{\ell =1}^{\infty}\ell\cdot\bE\big(|\fG_{\ell}(x) |\big) \,  dx
&=
\frac{1}{(1-\|T_n\|)^2} \big\langle\psi_n\, ,\, \mvone_{\circ} \big\rangle^2
+
\sum_{j=2}^{n}
\frac{1}{(1-\lambda_{n;j})^2} \big\langle\psi_{n; j}\, ,\, \mvone_{\circ} \big\rangle^2
\notag\\
&=
\frac{1}{(1-\|T_n\|)^2} \big(\medint{\int_0^1} \psi_n \big)^2
+
O(1)\, ,
\end{align*}
where we have used \eqref{eqn:105}, \eqref{eqn:765}, and the fact that
$
\sum_{j=2}^n 
\big\langle\psi_{n; j}\, ,\, \mvone_{\circ} \big\rangle^2
\leq 1
$.
The desired result now follows from \eqref{eqn:655}.
\qed

\medskip

\noindent{\bf Proof of Lemma~\ref{lem:67}~\eqref{item:bb}:}
For $i\in [n]$, let
\begin{align}\label{eqn:85}
\cD^{\bullet}(i)
:=
\sum_{j\in\cC^{\bullet}(i)} d_{\cG_n^{\bullet}}(i, j)
=
\sum_{\ell=1}^{\infty} \ell\cdot \cY_{\ell}^{\bullet}(i)
\, ,
\end{align}
where $\cY_{\ell}^{\bullet}(i)$ denotes the number of vertices in $\cC^{\bullet}(i)$ at distance $\ell$ from $i$.
Let $v$ be uniformly distributed over $[n]$.
Then note that 
$
\bE\big(\cD^{\bullet}\big)
=
\bE\big(\cD^{\bullet}(v)\big)
$
and
$
\bE\big[\int_0^1\fZ_n \big] 
=
\bE\big[\fZ_n^{\sss\triangle}(v)\big]
$.
Assume that $\cC^{\bullet}(v)$ is coupled with a multitype branching process started from a vertex of type $v$ as explained above \eqref{eqn:coupling}, and recall the construction of the BFS tree of $\cC^{\bullet}(v)$ rooted at $v$.
Note that for any $j\in\cC^{\bullet}(v)$, the distance between $v$ and $j$ in the BFS tree equals $d_{\cG_n^{\bullet}}(v, j)$.
Consequently, in the above coupling, 
\begin{align}\label{eqn:341}
|\fG_{\ell}^{\sss\triangle}(v)| \geq \cY_{\ell}^{\bullet}(v)
\ \text{ for each }\ell\geq 1\ \text{ and }\
\fZ_n^{\sss\triangle}(v)\geq \cD^{\bullet}(v)\, .
\end{align}
In particular,
$
\bE\big[\int_0^1\fZ_n \big] 
=
\bE\big[\fZ_n^{\sss\triangle}(v)\big]
\geq
\bE\big(\cD^{\bullet}(v)\big)
=
\bE\big(\cD^{\bullet}\big)
$.
Further, 
\begin{align*}
\fZ_n^{\sss\triangle}(v)
-
\cD^{\bullet}(v)
=
\sum_{\ell\geq 1}\ell\big(
|\fG_{\ell}^{\sss\triangle}(v)| - \cY_{\ell}^{\bullet}(v)
\big)
\leq
\fH_n^{\sss\triangle}(v) 
\big(\fX_n^{\sss\triangle}(v) 
- 
|\cC^{\bullet}(v)|
\big)\, .
\end{align*}
Thus, for any $m, n\geq 1$,
\begin{align*}
&
\bE\big[\fZ_n^{\sss\triangle}(v)\big]
-
\bE\big(\cD^{\bullet}(v)\big)
\leq
m\bigg(\medint{\int}_0^1 g_n - \bE[s_2^{\bullet}]\bigg)
+
\bE\bigg(
\fH_n^{\sss\triangle}(v) 
\cdot
\fX_n^{\sss\triangle}(v) 
\cdot
\ind_{\{\fH_n^{\sss\triangle}(v)\geq m\}}
\bigg)
\\[2pt]
&\hskip80pt
\leq
C m n^{4\delta_0 - 1}
+
\bigg[\bE\bigg( \fH_n^{\sss\triangle}(v)^2\cdot\ind_{\{\fH_n^{\sss\triangle}(v)\geq m\}} \bigg)\bigg]^{1/2}\cdot
\bigg(\medint{\int}_0^1 g_{n; 2}\bigg)^{1/2}
,
\end{align*}
where the last step uses Lemma~\ref{lem:54}~\eqref{item:54a}.
The claimed result now follows from \eqref{eqn:177}, \eqref{eqn:16}, and Lemma~\ref{lem:12}~\eqref{item:12}.
\qed

\medskip

\noindent{\bf Proof of Lemma~\ref{lem:67}~\eqref{item:cc}:}
Recall the notation $\cD^{\bullet}(\cdot)$ from \eqref{eqn:85}.
Let $u$ and $v$ be uniformly distributed over $[n]$ such that $u, v$, and $\cG_n^{\bullet}$ are independent.
Then note that 
$
\bE\big[ \cD^{\bullet}(u) \, |\, \cG_n^{\bullet} \big]
=
\cD^{\bullet}
=
\bE\big[ \cD^{\bullet}(v) \, |\, \cG_n^{\bullet} \big]
$.
Hence,
\begin{align}\label{eqn:764}
\bE\big[ (\cD^{\bullet})^2 \big]
=
\bE\big[ \cD^{\bullet}(u) \cD^{\bullet}(v) \big]
=
\bE\bigg( \cD^{\bullet}(v) \cdot 
\bE\bigg[ \cD^{\bullet}(u) 
\, \big|\, 
\cC^{\bullet}(v) \bigg] 
\bigg)
\, .
\end{align}

For any nonempty $A\subseteq [n]$, let $\cG^{A-}$ be the graph induced by $\cG_n^{\bullet}$ on the vertex set $[n]\setminus A$, as introduced below \eqref{eqn:768}.
For $i\in [n]\setminus A$, let $\cD^{A-}(i) := \sum_{j\in\cC^{A-}(i)} d_{\cG^{A-}}(i, j)$ where $\cC^{A-}(i)$ denotes the component of $\cG^{A-}$ containing $i$.
Similarly, for $i\in [n]\setminus A$, let $\fZ^{A-}(i)$ be the analogue of $\fZ_n^{\sss\triangle}(i)$ (defined in \eqref{eqn:52}) for a multitype branching process with type space $[n]\setminus A$ started from a vertex of type $i$ where each vertex of type $j\in [n]\setminus A$ has $\text{Bernoulli}(\kappa_{jk}/n)$ many type $k$ children for each $k\in [n]\setminus A$.
Then similar to \eqref{eqn:341}, we have $\fZ^{A-}(i)\geq \cD^{A-}(i)$ in the obvious coupling.
Thus, on the event $\big\{ V\big(\cC^{\bullet}(v) \big) = A \big\}$,
\begin{align}\label{eqn:766}
&
\bE\bigg[ \cD^{\bullet}(u) 
\, \big|\, 
\cC^{\bullet}(v) \bigg] 
=
\frac{1}{n}\sum_{i\in A}\sum_{j\in A} d_{\cC^{\bullet}(v)}(i, j)
+
\frac{1}{n}\sum_{i\in [n]\setminus A}\bE\big(\cD^{A-}(i)\big)
\notag\\
&\hskip40pt
\leq 
\frac{1}{n}\sum_{i\in A}\sum_{j\in A} d_{\cC^{\bullet}(v)}(i, j)
+
\frac{1}{n}\sum_{i\in [n]\setminus A}\bE\big(\fZ^{A-}(i)\big)
\notag\\
&\hskip80pt
\leq
\frac{1}{n}\sum_{i\in A}\sum_{j\in A} d_{\cC^{\bullet}(v)}(i, j)
+
\frac{1}{n}\sum_{i\in [n]}\bE\big(\fZ_n^{\sss\triangle}(i)\big)
=:
F_1 + F_2\, ,
\end{align}
where the third step uses the relation $\fZ^{A-}(i) \leq \fZ_n^{\sss\triangle}(i)$ for each $i\in [n]\setminus A$ in the obvious coupling between the two branching processes.
Now,
\begin{align}\label{eqn:767}
F_1
&
\leq
\frac{1}{n}\sum_{i\in A}\sum_{j\in A} \big(
d_{\cC^{\bullet}(v)}(i, v) + d_{\cC^{\bullet}(v)}(j, v) 
\big)
=
\frac{2}{n} |\cC^{\bullet}(v)| \cdot \cD^{\bullet}(v)\, ,
\end{align}
and
\begin{align}\label{eqn:700}
F_2
=
\bE\big[\medint\int_0^1\fZ_n\big]
\leq
\bE\big[ \cD^{\bullet}(v) \big]
+
Cn^{5\delta_0 - 1}\log n\, ,
\end{align}
where the last step uses Lemma~\ref{lem:67}~\eqref{item:bb}.

Combining \eqref{eqn:764}, \eqref{eqn:766}, \eqref{eqn:767}, and \eqref{eqn:700}, we see that
\begin{align}
\bE\big[ (\cD^{\bullet})^2 \big]
&
\leq 
\frac{2}{n}\bE\big( |\cC^{\bullet}(v)|\cdot \cD^{\bullet}(v)^2\big)
+
\big(\bE[\cD^{\bullet}(v)]\big)^2
+
Cn^{5\delta_0 - 1}\log n \cdot \bE\big[\cD^{\bullet}(v)\big]
\notag\\
&\hskip30pt
\leq
\frac{2}{n}\bE\big[ \fX_n^{\sss\triangle}(v)\cdot \fZ_n^{\sss\triangle}(v)^2\big]
+
\big(\bE[\cD^{\bullet}]\big)^2
+
Cn^{7\delta_0 - 1}\log n\, ,\label{eqn:701}
\end{align}
where the last step uses Lemma~\ref{lem:67}~\eqref{item:aa}, \eqref{item:bb} and the fact that $\bE[\cD^{\bullet}] = \bE[\cD^{\bullet}(v)]$.
Using the relation 
$
\fZ_n^{\sss\triangle}(v)
\leq
\fH_n^{\sss\triangle}(v)\cdot\fX_n^{\sss\triangle}(v)
$, 
we see that for any $n, m\geq 1$,
\begin{align}
\bE\big[ \fX_n^{\sss\triangle}(v)\cdot \big(\fZ_n^{\sss\triangle}(v) \big)^2\big]
&
\leq
\bE\big[ \fX_n^{\sss\triangle}(v)^3\cdot \big(\fH_n^{\sss\triangle}(v) \big)^2\big]
\notag\\
&
\leq
m^2\bE\big[ \fX_n^{\sss\triangle}(v)^3\big]
+
\bE\big[ 
\fX_n^{\sss\triangle}(v)^3\cdot 
\big(\fH_n^{\sss\triangle}(v) \big)^2\cdot
\ind_{\{\fH_n^{\sss\triangle}(v)\geq m\}}
\big]
\notag\\[2pt]
&
\leq
m^2\big(\medint{\int}_0^1 g_{n; 3}\big)
+
\big(\medint{\int}_0^1 g_{n; 4}\big)^{3/4}
\cdot
\bigg[\bE\big(\fH_n^{\sss\triangle}(v)^8\cdot\ind_{\{\fH_n^{\sss\triangle}(v)\geq m\}}\big)\bigg]^{1/4}  \, .
\label{eqn:702}
\end{align}
Now the claimed bound on $\var\big(\cD^{\bullet}\big)$ follows upon combining \eqref{eqn:701} with \eqref{eqn:702}, and using Lemma~\ref{lem:12}~\eqref{item:13}, \eqref{eqn:177}, and \eqref{eqn:16}.
\qed

\subsection{Proof of Proposition~\ref{prop:graphon-intermediate}}\label{sec:g11}
Throughout this section, we will work under Conditions~\ref{ass:graphon-1}, \ref{ass:graphon-2}, and \ref{ass:graphon-3a}.
Further, we assume that \eqref{eqn:765} and \eqref{eqn:17} hold.
Recall the connection probabilities $p_{ij}^{\sss (n)}$ and $\widetilde p_{ij}^{\sss (n)}$ from \eqref{eqn:222}.
Then
\begin{align}\label{eqn:901}
\sum_{\substack{ i, j\in [n]\setminus B_n \\  i\neq j }}
\frac{\big( p_{ij}^{\sss (n)}- \widetilde p_{ij}^{\sss (n)} \big)^2}{p_{ij}^{\sss (n)}}
\leq
\frac{C}{n^3}\sum_{i, j\in [n]}\beta_{ij }^3
= O(1/n) \, ,
\end{align}
where the last step uses \eqref{eqn:20}.
Hence, by \cite[Corollary~2.12~(i)]{janson-equiv}, for all $n$ there exists a coupling of $\graphon_n$ and $\tildegraphon_n$ such that 
\begin{align}\label{eqn:902}
\pr\big(\graphon_n \neq \tildegraphon_n\big)\to 0\ \ \text{ as }\ \ n\to\infty\, .
\end{align}
Thus, it is enough to show that \eqref{eqn:65} and \eqref{eqn:66} hold for $\tildegraphon_n$.

Let $\widetilde\cG_n^{\bullet}$ be the random graph with vertex set $[n]$ where the edges $\{i, j\}$ appear independently with respective probabilities
$1 - \exp\big(-\kappa_{ij} /n\big)$, $1\leq i<j\leq n$.
By an argument similar to the one leading to \eqref{eqn:902}, we can couple $\widetilde\cG_n^{\bullet}$ with $\cG_n^{\bullet}$ for all $n$ such that 
\begin{align}\label{eqn:34a}
\pr\big(\widetilde\cG_n^{\bullet} \neq \cG_n^{\bullet}\big)\to 0\ \ \text{ as }\ \ n\to\infty\, .
\end{align}
Let
$
\widetilde D_{\max}^{\bullet} := \diam(\widetilde\cG_n^{\bullet})$,
$\widetilde\cD^{\bullet}:= \cD(\widetilde\cG_n^{\bullet})$, and 
$\widetilde s_k^{\bullet} := s_k(\widetilde\cG_n^{\bullet})$, $k\geq 1$.
Excluding the singleton vertices belonging to $B_n$ from $\widetilde\cG_n^{\bullet}$, arrange the remaining (say $m_n$ many) components in decreasing order of their sizes with ties broken following some deterministic rule;
let $\widetilde\cC_i^{\bullet}$ denote the $i$-th largest component, $1\leq i\leq m_n$.
Our plan is to apply Theorem \ref{thm:gen-2} with
$\vx=(x_1, \ldots, x_{m_n}), q$, and $\vM=(M_1, \ldots, M_{m_n})$, where
\begin{equation}\label{eqn:510}
	x_i = n^{-2/3}\chi^{1/3}|\widetilde\cC_i^{\bullet}| \, ,\ 
	q =\chi^{-2/3}n^{1/3-\delta_0} \, ,
	\ \text{ and }\
	M_i =  \scl\big(1, 1/|\widetilde\cC_i^{\bullet}| \big) \widetilde\cC_i^{\bullet}\, ,\ \ 1\leq i\leq m_n\, .
\end{equation}

Writing $\sigma_k=\sum_{i=1}^{m_n} x_i^k$, we have
\begin{equation}\label{eqn:991}
	\sigma_2 = \frac{\chi^{2/3}}{n^{1/3}}\bigg( \widetilde s_2^{\bullet} -  \frac{|B_n|}{n}\bigg) \, ,\
	\sigma_3 =  \frac{\chi}{n}\bigg( \widetilde s_3^{\bullet} -  \frac{|B_n|}{n}\bigg) \, ,\
	x_{\max} =\frac{\chi^{1/3}|\widetilde\cC_1^{\bullet}|}{n^{2/3}}\, ,\text{ and }
	x_{\min} \geq \frac{\chi^{1/3}}{n^{2/3}} \, .
\end{equation}
Using \eqref{eqn:34a} in conjunction with Lemma~\ref{lem:20a}, \eqref{eqn:982}, and Lemma~\ref{lem:5}, we see that Condition~\ref{ass:aldous-basic-assumption} holds with $\lambda=\zeta\chi^{-2/3}$. 
Next, we can clearly choose $r_0>0$ large so that $\sigma_2^{r_0}/x_{\min}\weakc 0$ as $n\to\infty$.
Further, \eqref{eqn:34a} and \eqref{eqn:982} show that for any $\eta_0>0$,
\[
x_{\max}/\sigma_2^{3/2+\eta_0}
=
O_P\Big(
n^{-\frac{1}{6} + \frac{\delta_0}{2} + \theta_0 + \eta_0\big(\frac{1}{3}-\delta_0\big)}
\cdot
\log n
\Big).
\]
Since $\delta_0<1/3 - 2\theta_0$, it follows that $x_{\max}/\sigma_2^{3/2+\eta_0}\weakc 0$ as $n\to\infty$ if we choose $\eta_0$ sufficiently small.
This shows that Condition~\ref{ass:gen-1} holds.
Turning to the verification of Condition~\ref{ass:gen-2}, note that
\begin{align*}
\tau
=
\sum_{\ell=1}^{m_n} x_\ell^2 u_{\ell}
= n^{-4/3}\chi^{2/3}
\sum_{\ell=1}^{m_n} 
\sum_{i,j \in\widetilde\cC_\ell^{\bullet}} d_{\widetilde\cC_\ell^{\bullet}}(i,j) 
= 
n^{-1/3}\chi^{2/3}\widetilde\cD^{\bullet}
\sim
\alpha\chi^{2/3}n^{2\delta_0-1/3}
\end{align*}
and $\diam_{\max}=\widetilde D_{\max}^{\bullet} = O_P(n^{\delta_0}\log n)$,
where we have used \eqref{eqn:34a}, Lemma~\ref{lem:66}, and \eqref{eqn:981}.
Using these relations together with the fact that 
$x_{\max} = O_P( n^{-2/3+2\delta_0+\theta_0} \log n  )$
and
$\sigma_2 \sim \chi^{2/3} n^{\delta_0 - 1/3}$, we see that
\[
\frac{\sigma_2^{\frac{3}{2}-\eta_0}\cdot\diam_{\max}}{\tau+\sigma_2}
=
O_P\Big(n^{-\big(\frac{1}{3}-\delta_0\big)\big(\frac{1}{2}-\eta_0\big)}\cdot \log n\Big)
\ \ \text{ and }\ \
\frac{\sigma_2 x_{\max}\diam_{\max} }{\tau}
=
O_P\Big(\frac{(\log n)^2}{n^{\frac{2}{3}-2\delta_0-\theta_0}}\Big)
\, .
\]
Since $\delta_0<1/3-2\theta_0$, it follows that Condition~\ref{ass:gen-2} holds for any choice of $\eta_0<1/2$.

Similar to the notation introduced at the beginning of Section~\ref{sec:smc-def}, we let $\cC_i$ be the $i$-th maximal component (with respect to the mass given by 
$\mass(\cC) = \sum_{i\in\cC}x_i$
for each component $\cC$) in $\cG(\vx,q)$.
Following the notation introduced around \eqref{eqn:uik-def}, let $\bar\cC_i$ be the component in $\bar\cG(\vx, q, \vM)$ corresponding to the component $\cC_i$ in  $\cG(\vx,q)$.
We will write $\bar\cC_i$ to denote both the graph and the metric measure space (where the metric is given by \eqref{eqn:44} and coincides with the graph distance, and the measure is as described in \eqref{eqn:barmu-on-full-met} and coincides with $n^{-2/3}\chi^{1/3}$ times the counting measure on the vertices).
Note that $\bar\cC_i$ is the $i$-th maximal component (in terms of the number of vertices) of 
$\bar\cG(\vx, q, \vM)$.
Since Condition~\ref{ass:aldous-basic-assumption} and Condition~\ref{ass:gen-2} are satisfied, the conclusions of Theorem~\ref{thm:aldous-review} and Theorem~\ref{thm:gen-2} hold for $(\cC_i\, ;\, i\geq 1)$ and $(\bar\cC_i\, ;\, i\geq 1)$.
Thus,
\begin{align}\label{eqn:601}
\bigg(
\bigg( n^{-2/3}\chi^{1/3} | \bar\cC_i |\, ,\, \spls(\cC_i) \bigg)
\, ;\, i \geq 1 
\bigg) 
\weakc 
\mvXi\big(\zeta\chi^{-2/3}\big)
\end{align}
as $n \to \infty$ with respect to the product topology and
\begin{align}\label{eqn:602}
\Big(
\scl\Big(\frac{\chi^{2/3}}{\alpha n^{1/3}}\, , \, 1 \Big) \bar \cC_i
\, ;\, i \geq 1 
\Big) 
\weakc 
\vCrit\big(\zeta\chi^{-2/3}\big)
\end{align}
as $n \to \infty$, where we have used the fact that $\sigma_2^2/(\sigma_2 + \tau) \sim \chi^{2/3}/(\alpha n^{1/3})$.

We will now use \eqref{eqn:601} and \eqref{eqn:602} to deduce the convergences \eqref{eqn:65} and \eqref{eqn:66} for $\tildegraphon_n$.
Note that $\tildegraphon_n$ can be constructed from $\widetilde\cG_n^{\bullet}$ as follows:
Consider all pairs $i, j\in [n]\setminus B_n$ with $i\neq j$ such that there is no edge between $i$ and $j$ in $\widetilde\cG_n^{\bullet}$, and then place edges between them independently with probability 
$1-\exp\big(-1/n^{1+\delta_0}\big)$.
On the other hand, 
\chsen{we construct $\bar\cG(\vx, q, \vM)$ from $\widetilde\cG_n^{\bullet}$ by ignoring the singleton vertices belonging to $B_n$ and repeating the following procedure independently for $1\leq i<j\leq m_n$:
place a single edge between $\widetilde\cC_i^{\bullet}$ and $\widetilde\cC_j^{\bullet}$ with probability $1-\exp(-|\widetilde\cC_i^{\bullet}|\cdot |\widetilde\cC_j^{\bullet}|/n^{1+\delta_0})$,
and if there is such an edge, then choose it to be one of the $|\widetilde\cC_i^{\bullet}|\cdot|\widetilde\cC_j^{\bullet}|$ possible edges uniformly.
There is a natural way of coupling $\widetilde\cG_n^{\bullet}$, $\tildegraphon_n$, $\bar\cG(\vx, q, \vM)$, and $\cG(\vx, q)$ as follows:
Construct $\tildegraphon_n$ from $\widetilde\cG_n^{\bullet}$ as described above;
suppose that during this construction, $k_{ij}$ many edges are added between $\widetilde\cC_i^{\bullet}$ and $\widetilde\cC_j^{\bullet}$, $1\leq i\leq j\leq m_n$.
Declare $\bar\cG(\vx, q, \vM)$ to be the space obtained by performing the following three operations on $\tildegraphon_n$:
(i) Delete the $k_{ii}$ edges added within $\widetilde\cC_i^{\bullet}$, $1\leq i\leq m_n$.
(ii) For each $1\leq i<j\leq m_n$, if $k_{ij}\geq 1$, then retain a uniformly chosen edge among the $k_{ij}$ edges added between $\widetilde\cC_i^{\bullet}$ and $\widetilde\cC_j^{\bullet}$, and delete the other $k_{ij}-1$ edges.
(iii) Delete the singleton vertices belonging to $B_n$ from $\tildegraphon_n$.
We declare $\cG(\vx, q)$ to be the following random graph on $[n]\setminus B_n$:
Place an edge between $i$ and $j$ if $k_{ij}\geq 1$, $1\leq i<j\leq m_n$.
}

\chsen{For simplicity, we only work with the maximal component $\cC_1(\tildegraphon_n)$.	
Note that in the above  coupling
$V\big(\cC_1(\tildegraphon_n)\big)=V(\bar\cC_1)$, 
and 
$E\big(\cC_1(\tildegraphon_n)\big)\supseteq E(\bar\cC_1)$.
Suppose we show that $E\big(\cC_1(\tildegraphon_n)\big) = E(\bar\cC_1)$ whp.
Then it will follow that $\cC_1(\tildegraphon_n) = \bar\cC_1$ whp, and \eqref{eqn:602} will imply that
\begin{align}\label{eqn:602-a}
\frac{\chi^{2/3}}{\alpha n^{1/3}} \cdot \cC_1(\tildegraphon_n) 
\weakc 
 \crit_1\big(\zeta \chi^{-2/3}\big)
\end{align}
as $n\to\infty$ with respect to the GHP topology.
Now, $\bar\cC_1$ is simply $\cC_1$ with its vertices replaced by some blobs.
Suppose, we show that the blobs $\widetilde\cC_i^{\bullet}$ that become part of $\bar\cC_1$ are all trees whp.
Then it will follow that $\spls( \bar\cC_1)  = \spls( \cC_1 )$ whp, which together with the relation $\cC_1(\tildegraphon_n) = \bar\cC_1$ whp will imply that 
$\spls(\cC_1) = \spls\big( \cC_1(\tildegraphon_n) \big)$ whp.
Then \eqref{eqn:601} will yield
\begin{align}\label{eqn:601-a}
\bigg( n^{-2/3}\chi^{1/3} | \cC_1(\tildegraphon_n) |\, ,\, \spls\big( \cC_1(\tildegraphon_n) \big) \bigg)
\weakc 
\bigg( 
\gamma_1\big( \zeta\chi^{-2/3} \big)\, ,\, \pois_1\big(  \zeta\chi^{-2/3} \big)
\bigg)
\end{align}
as $n\to\infty$.
The converegences in \eqref{eqn:602-a} and \eqref{eqn:601-a} can be strengthened to joint convergence for all the components in a straightforward way.
Thus, the following lemma completes the proof of  Proposition~\ref{prop:graphon-intermediate}.
}

\begin{lem}\label{lem:47}
The following hold in the above coupling:
\begin{enumeratea}
\item\label{item:47a}
With high probability, $E\big(\cC_1(\tildegraphon_n)\big) = E(\bar\cC_1)$.
\vskip2pt
\item\label{item:47b}
With high probability, $\widetilde\cC_i^{\bullet}$ is a tree for all $i$ such that $\widetilde\cC_i^{\bullet}\subseteq\bar\cC_1$.
\end{enumeratea}
\end{lem}

\noindent{\bf Proof of Lemma~\ref{lem:47}~\eqref{item:47a}:} 
Let $\cF$ denote the $\sigma$-field generated by 
$\widetilde\cG_n^{\bullet}$, $\cG(\vx, q)$, and $\bar\cG(\vx, q, \vM)$.
Write $\sum_{\dagger}$ for the sum over all $ k\in [m_n]$ such that $k\in V(\cC_1)$,
and
$\sum_{\dagger\dagger}$ for the sum over all $1\leq i<j\leq m_n$ such that  $\{i, j\}\in E(\cC_1)$.
Then conditional on $\cF$, 
$\# E\big(\cC_1(\tildegraphon_n)\big) - \# E(\bar\cC_1)$
is stochastically bounded by
\begin{align}\label{eqn:111}
\sum\displaystyle_{\dagger\dagger} \ (N_{ij} - 1)
+
\sum\displaystyle_{\dagger}\ N_{kk} \, ,
\end{align}
where
(i) for each  $i<j$ with $\{i, j\}\in E(\cC_1)$, $N_{ij}$ is distributed as a 
$\text{Binomial}\big( |\widetilde\cC_i^{\bullet}|\cdot|\widetilde\cC_j^{\bullet}|
\, ,\, 
1-\exp(-1/n^{1+\delta_0} )\big)$ 
random variable conditioned to be at least one,
(ii)  for each $k\in V(\cC_1)$, $N_{kk}$ is distributed as a 
$\text{Binomial}\big( |\widetilde\cC_k^{\bullet}|^2
\, ,\, 
1-\exp(-1/n^{1+\delta_0} )\big)$ 
random variable, and
(iii) all the random variables from (i) and (ii) are jointly independent.
Then conditional on $\cF$, on the event 
$\{\,
|\widetilde\cC_1^{\bullet} |
\leq 
n^{2\delta_0 + \theta_0}(\log n)^2
\}$,
the expectation of the expression in \eqref{eqn:111} is bounded by
\begin{align*}
Cn^{-1-\delta_0}\bigg(
\sum\displaystyle_{\dagger\dagger}\
 |\widetilde\cC_i^{\bullet}|\cdot|\widetilde\cC_j^{\bullet}|
+
\sum\displaystyle_{\dagger}\  |\widetilde\cC_k^{\bullet}|^2
\bigg)
\leq
Cn^{-1-\delta_0}\cdot | \widetilde\cC_1^{\bullet} |\cdot
\bigg(
|\bar\cC_1|
+
\spls(\cC_1) | \widetilde\cC_1^{\bullet} |
+
 |\bar\cC_1|
\bigg)
\, .
\end{align*}
It follows from \eqref{eqn:601} that $|\bar\cC_1| = O_P(n^{2/3})$ and $\spls(\cC_1) = O_P(1)$, and consequently, using the fact that $\delta_0<1/3 - \theta_0$, we see the right side of the last display is $o_P(1)$.
Since 
$\{\,
|\widetilde\cC_1^{\bullet} |
\leq 
n^{2\delta_0 + \theta_0}(\log n)^2
\}$ whp by \eqref{eqn:34a} and \eqref{eqn:982}, the proof is complete.
\qed

\medskip

\noindent{\bf Proof of Lemma~\ref{lem:47}~\eqref{item:47b}:} 
Recall the notation $\connects$, $\connectsell$, and $\leftrightarrow$ introduced below \eqref{eqn:898}.
Note that if $\spls\big(\cC^{\bullet}(i)\big)\geq 1$ for some $i\in [n]$, then either there exist $j, k\in [n]$ with $i, j, k$ pairwise distinct such that the event 
$
\{i\connects j\}\circ\{j\connectsell k\}\circ\{j\leftrightarrow k\}
$
occurs for some $\ell\geq 2$,
or there exists $k\neq i$ such that the event 
$
\{i\connectsell k\}\circ\{i\leftrightarrow k\}
$
occurs for some $\ell\geq 2$.
Let $v$ be uniformly distributed over $[n]$ independent of $\cG_n^{\bullet}$.
Using Lemma~\ref{lem:bk-inequality} and Lemma~\ref{lem:55}, and writing `p.d.' for `pairwise distinct,' we get
\begin{align}\label{eqn:538}
	\pr\big(\spls\big(\cC^{\bullet}(v)\big)\geq 1\big)
	\leq
	F_1 + F_2 \ ,
\end{align}
where
\begin{align}\label{eqn:538a}
	F_1
	&
	=
	n^{-1}\sum_{\ell\geq 2}\ \sum_{i, j, k \text{ (p.d.)}}
	\pr\big(i\connects j\big)  \pr\big(j\connectsell k\big)
    \pr\big(j\leftrightarrow k\big)
	\notag\\
	&\hskip20pt
	\leq
	n^{-1} \medint{\int}_{[0, 1]^3}
	\rho_n(x, y) \rho_n(y, z)\kappa_n(y, z)\, dx dy dz\, ,
\end{align}
and similarly
\begin{gather*}
	F_2
	\leq
	n^{-1}\medint{\int}_{[0, 1]^2}
	\rho_n(x, z)\kappa_n(x, z)\, dx dz\, .
	\label{eqn:538b}
\end{gather*}
It follows from \eqref{eqn:538a} that 
\begin{align}\label{eqn:539}
F_1
\leq
n^{-1}\big(\medint{\int}_{[0, 1]^2} \rho_n^3\big)^{1/3}
\cdot
\big(\medint{\int}_{[0, 1]^2} \rho_n^3\big)^{1/3}
\cdot
\big(\medint{\int}_{[0, 1]^2} \kappa_n^3\big)^{1/3}
\lesssim
n^{2\delta_0 - 1}\, ,
\end{align}
where we have used \eqref{eqn:125} and the fact $\int_{[0, 1]^2}\kappa_n^3 = O(1)$.
Similarly $F_2\lesssim n^{\delta_0 - 1}$, which combined with \eqref{eqn:539} and \eqref{eqn:538} yields
\begin{align}\label{eqn:640}
\bE\bigg(
\frac{1}{n}\sum_i |\cC_i^{\bullet}|\cdot a_i^{\bullet}
\bigg)
=
\pr\big(\spls\big(\cC^{\bullet}(v)\big)\geq 1\big)
\lesssim
n^{2\delta_0 - 1}\, ,
\end{align}
where $\cC_i^{\bullet}$ is $i$-th largest component of $\cG_n^{\bullet}$ (as defined below \eqref{eqn:898}) and 
$a_i^{\bullet}:=\ind_{\{\cC_i^{\bullet}\text{ is not a tree}\}}$.
For $1\leq i\leq m_n$, define 
$\widetilde a_i^{\bullet}:=\ind_{\{\widetilde\cC_i^{\bullet}\text{ is not a tree}\}}$.
Then using \eqref{eqn:34a}, \eqref{eqn:640}, and the fact that $|B_n| = o(n)$, we get
\begin{align}\label{eqn:641}
\frac{\sum_{i=1}^{m_n} x_i\widetilde a_i^{\bullet}}{\sum_{i=1}^{m_n} x_i}
=
\frac{\sum_{i=1}^{m_n} |\widetilde\cC_i^{\bullet}|\cdot\widetilde a_i^{\bullet}}{n-|B_n|}
=
o_P\big(n^{2\delta_0 - 1} \log n\big)
\, .
\end{align}

Consider the breadth-first exploration of $\cG(\vx, q)$ as described in Section~\ref{sec:def-size-bias}, and suppose $v(1),\ldots, v(m_n)$ is the order in which the blobs appear in this exploration.
Now, with $\vx$ as in \eqref{eqn:510}, we have $\sigma_1/\sigma_2=o_P(n^{2/3 - \delta_0}\log n)$.
Thus, by the observation made below \eqref{eqn:1438}, whp the exploration of the maximal component $\cC_1$ is completed in the first $n^{2/3 - \delta_0}\log n$ many steps, i.e.,
\begin{align}\label{eqn:642}
\pr\Big(
\big\{ i \in [m_n]\, :\, i\in V(\cC_1) \big\}
\subseteq
\big\{v(i)\, :\, 1\leq i\leq n^{2/3 - \delta_0}\log n \big\}
\Big)
\to
1
\end{align}
as $n\to\infty$.
Since  $v(1),\ldots, v(m_n)$ is a size-biased permutation of the blobs using the weight sequence $\vx$, we have, on the event 
\begin{align}\label{eqn:643}
\bigg\{ 
x_{\max}\leq n^{2\delta_0 + \theta_0 -2/3} (\log n)^2
\ \ \text{ and }\ \ 
\big(\sum_{i=1}^{m_n} x_i\widetilde a_i^{\bullet}\big)\big/\big(\sum_{i=1}^{m_n} x_i \big)
\leq
n^{2\delta_0 - 1}\log n
\bigg\}\, ,
\end{align}
\begin{align}
&\pr\big( \widetilde a^{\bullet}_{v(i)} = 1 \text{ for some }i\leq n^{2/3 - \delta_0}\log n
\, \big|\, 
\widetilde\cG_n^{\bullet}\big)
\leq
\frac{
n^{2/3 - \delta_0}\log n\cdot
\sum_{i=1}^{m_n} x_i\widetilde a_i^{\bullet}
}{
\sum_{i=1}^{m_n} x_i - x_{\max}\cdot n^{2/3 - \delta_0}\log n
}
\notag
\\
&\hskip50pt
\leq
2n^{2/3 - \delta_0}\log n\cdot
\frac{
\sum_{i=1}^{m_n} x_i\widetilde a_i^{\bullet}
}{
\sum_{i=1}^{m_n} x_i 
}
\leq
2n^{\delta_0 - 1/3}(\log n)^2
=
o(1)
\, ,
\label{eqn:644}
\end{align}
where the second step uses the relation $\delta_0<1/3 - \theta_0$ and holds for all large $n$.
By \eqref{eqn:34a} and \eqref{eqn:982},
$x_{\max}=o_P\big(n^{2\delta_0 + \theta_0 -2/3} (\log n)^2\big)$ 
which combined with \eqref{eqn:641} shows that the event in \eqref{eqn:643} occurs whp.
Thus, combining \eqref{eqn:642} and \eqref{eqn:644} completes the proof.
\qed

\subsection{Proofs of Theorem~\ref{thm:scaling-limit-graphon} and
Theorem~\ref{thm:component-limit-graphon}}\label{sec:g12}

We will first prove the following result:
 
\begin{lem}\label{lem:23}
Suppose Conditions~\ref{ass:graphon-1} and \ref{ass:graphon-2} hold, and in addition, with $\theta_0$ as in Condition~\ref{ass:graphon-2},
there exist $\delta_0\in (0\, ,\, 1/3-2\theta_0)$ and $\varpi_0\in [0\, ,\, 2\delta_0+1/3)$ such that the following hold:
For each $n\geq 1$, there exists $B_n\subseteq [n]$ satisfying \eqref{eqn:31}, and further,
$\beta_{ij} > n^{-\delta_0}$ whenever $i, j\in [n]\setminus B_n$ with $i\neq j$.
Then the convergences in \eqref{eqn:65} and \eqref{eqn:66} hold for $\graphon_n$ with the same constants as in \eqref{eqn:def-alpha-chi-zeta}.
\end{lem}
 
\noindent{\bf Proof:}
Let $\beta_{ii}^{\pm}:=0$ for all $i\in [n]$, and for $1\leq i\neq j\leq n$, define
\begin{align*}
\beta_{ij}^{-}
:=
\left\{
\begin{array}{l}
\beta_{ij }\, , \text{ if } i, j\in [n]\setminus B_n\, ,\\[2pt]
0,\text{ otherwise,}
\end{array}
\right.
\qquad \text{ and }\qquad
\beta_{ij}^{+}
:=
\left\{
\begin{array}{l}
\beta_{ij }\, , \text{ if } i, j\in [n]\setminus B_n\, ,\\[2pt]
\beta_{ij} + n^{-\delta_0},\text{ otherwise}.
\end{array}
\right.
\end{align*}
Let $\mvbeta_n^{\pm} := \big( \beta_{ij}^{\pm} \, ;\, i, j\in[n]\big)$.
Clearly, $\mvbeta_n^{+}$ and $\mvbeta_n^{-}$ satisfy Conditions~\ref{ass:graphon-2} and \ref{ass:graphon-3a}.
Define the functions $W_n^-$ and $W_n^+$ on $[0, 1]^2$ similar to \eqref{eqn:2} using the edge weights $\mvbeta_n^{-}$ and $\mvbeta_n^{+}$ respectively.
Then \eqref{eqn:20} is satisfied by both $W_n^-$ and $W_n^+$.
Next, note that
\begin{align*}
\bigg|
\big\| 
T_{H-n^{\sss 1/3}(W_n-W)}
\big\|_{2,2}
-
\big\| 
T_{H-n^{\sss 1/3}(W_n^--W)}
\big\|_{2,2}
\bigg|^2
\leq
n^{\frac{2}{3}}\medint{\int}_{[0, 1]^2}\big( W_n - W_n^-\big)^2
\leq
\frac{2}{n^{4/3}}\sum_{i\in B_n}\sum_{j=1}^n
\beta_{ij}^2
=
o(1)\, ,
\end{align*}
where the last step uses the second relation in \eqref{eqn:31}.
Thus, \eqref{eqn:21} is satisfied by the sequence $\big(W_n^-\, ;\, n\geq 1\big)$.
We can similarly show that \eqref{eqn:21} is satisfied by the sequence $\big(W_n^+\, ;\, n\geq 1\big)$ by using the first relation in \eqref{eqn:31}.
Hence, both the sequences $\big( \mvbeta_n^+ \big)$ and $\big( \mvbeta_n^- \big)$ satisfy Conditions~\ref{ass:graphon-1}, \ref{ass:graphon-2}, and \ref{ass:graphon-3a}, and consequently, Proposition~\ref{prop:graphon-intermediate} is applicable in both cases.

Let $\cG_n^{\pm}$ denote the random graph on $[n]$ where the edges $\{i, j\}$ appear independently with respective probabilities $1\wedge (\beta_{ij}^{\pm}/n)$, $1\leq i<j \leq n$.
Since $\beta_{ij}^- \leq \beta_{ij} \leq \beta_{ij}^+$ for all $i, j\in [n]$, there is a natural coupling in which \begin{align}\label{eqn:908}
\cG_n^-\subseteq \graphon_n \subseteq \cG_n^+\, .
\end{align}
Assume that $\cG_n^-, \graphon_n$, and $\cG_n^+$ are coupled in this way for all $n$.
For $i\geq 1$, let $\cC_i(\cG_n^{\pm})$ denote the $i$-th largest component of $\cG_n^{\pm}$, and view it as a metric measure space using the counting measure on the vertices.
In the coupling of \eqref{eqn:908}, we have 
\begin{align}\label{eqn:909}
|\cC_1(\cG_n^-)|
\leq
|\cC_1(\graphon_n)|
\leq
|\cC_1(\cG_n^+)|
\ \ 
\text{ for all }\ \ n.
\end{align}
Recall the random variables $\gamma_i(\cdot)$ and $\pois_i(\cdot)$ defined below \eqref{eqn:parabolic-bm}, and write
$\lambda = \chi^{-2/3}\zeta$.
It follows from Proposition~\ref{prop:graphon-intermediate} that
$n^{-2/3}\chi^{1/3}|\cC_1(\cG_n^{\pm})|\weakc\gamma_1(\lambda)$,
which combined with \eqref{eqn:909} yields 
$n^{-2/3}\big( |\cC_1(\cG_n^+)| -  |\cC_1(\cG_n^-)| \big)\weakc 0$.
Since, $\gamma_1(\lambda)>\gamma_2(\lambda)$ a.s., we get
$
\cC_1(\cG_n^-)
\subseteq
\cC_1(\graphon_n)
\subseteq
\cC_1(\cG_n^+)
$
whp.
This implies that 
$
\spls\big(\cC_1(\cG_n^-)\big)
\leq
\spls\big(\cC_1(\graphon_n)\big)
\leq
\spls\big(\cC_1(\cG_n^+)\big)
$
whp.
Now, by Proposition~\ref{prop:graphon-intermediate},
$
\spls\big(\cC_1(\cG_n^{\pm})\big)
\weakc
\pois_i(\lambda)
$
as $n\to\infty$,
and consequently, we have
$
\spls\big(\cC_1(\cG_n^-)\big)
=
\spls\big(\cC_1(\graphon_n)\big)
=
\spls\big(\cC_1(\cG_n^+)\big)
$
whp.
Proceeding inductively, we can show that for each $i\geq 1$,
\begin{gather}
\cC_i(\cG_n^-)
\subseteq
\cC_i(\graphon_n)
\subseteq
\cC_i(\cG_n^+)
\ \text{ whp, }\ 
n^{-2/3}\big(
|\cC_i(\cG_n^+)|
-
|\cC_i(\cG_n^-)|
\big)
\weakc 
0\, , 
\ \text{ and}
\label{eqn:910}\\
\spls\big(\cC_i(\cG_n^-)\big)
=
\spls\big(\cC_i(\graphon_n)\big)
=
\spls\big(\cC_i(\cG_n^+)\big)
\ \text{ whp.} 
\label{eqn:911}
\end{gather}
This shows that \eqref{eqn:66} holds.

Turning to the proof of \eqref{eqn:65}, 
let 
\[
\cE_n
:=
\big\{
\cC_1(\cG_n^-)
\subseteq
\cC_1(\graphon_n)
\subseteq
\cC_1(\cG_n^+)
\, ,\
\spls\big(\cC_1(\cG_n^-)\big)
=
\spls\big(\cC_1(\graphon_n)\big)
=
\spls\big(\cC_1(\cG_n^+)\big)
\big\}.
\]
Note that by Proposition~\ref{prop:graphon-intermediate}, 
\begin{align}\label{eqn:913}
\scl\Big( \frac{\chi^{2/3} }{\alpha n^{1/3}}\, ,\, \frac{\chi^{1/3}}{n^{2/3}} \Big) 
\cC_1(\cG_n^{\pm})
\weakc 
\crit_1(\lambda)
\ \ \text{ as }\ \ n\to\infty
\end{align}
with respect to the GHP topology, which combined with 
\cite[Proposition~6.5]{addarioberry-sen} yields
\begin{align}\label{eqn:912}
n^{-1/3}\cdot
\ind_{\cE_n}
\cdot
d_{\text{H}}\big(
\cC_1(\cG_n^{-}) \, ,\, 
\cC_1(\cG_n^{+})
\big)
\weakc 
0\, .
\end{align}
On the event $\cE_n$, 
$\cC_1(\cG_n^{+})$ can be obtained from $\cC_1(\cG_n^{-})$ by attaching finitely many, say $k_n$ trees $\fT_1, \ldots, \fT_{k_n}$; define $k_n=0$ on $\cE_n^c$.
Then \eqref{eqn:912} implies that 
$
\max_{1\leq i\leq k_n}\diam(\fT_i) = o_P(n^{1/3})
$.
Now, on the event $\cE_n$, $\cC_1(\graphon_n)$ is the space $\cC_1(\cG_n^{-})$ with some additional trees attached each of which is a subtree of $\fT_i$ for some $i\leq k_n$.
We thus have 
\begin{align}\label{eqn:914}
\ind_{\cE_n}
\cdot
d_{\text{H}}\big(
\cC_1(\cG_n^{-}) \, ,\, 
\cC_1(\graphon_n)
\big)
=
o_P(n^{1/3})
\ \text{ and \ }
|\cC_1(\graphon_n)|
-
|\cC_1(\cG_n^-)|
=
o_P(n^{2/3})
\, ,
\end{align}
where the second relation uses \eqref{eqn:910}.
Using \eqref{eqn:914} and the fact that $\pr(\cE_n)\to 1$, it is easy to show that
\[
d_{\GHP}
\bigg(
\scl\big( n^{-1/3}\, ,\, n^{-2/3} \big) 
\cC_1(\cG_n^{-})
\, ,\
\scl\big( n^{-1/3}\, ,\, n^{-2/3} \big) 
\cC_1(\graphon_n)
\bigg)
\weakc 0\, .
\]
It thus follows that \eqref{eqn:913} holds if we replace $\cC_1(\cG_n^{\pm})$ by $\cC_1(\graphon_n)$.
This argument can be extended in a straighforward way to establish the convergence in \eqref{eqn:65} with respect to the product topology.
\qed

\medskip

\noindent{\bf Completing the proofs of Theorems~\ref{thm:scaling-limit-graphon} and \ref{thm:component-limit-graphon}:}
Suppose $(\mvbeta_n)$ satisfies Conditions~\ref{ass:graphon-1}, \ref{ass:graphon-2}, and \ref{ass:graphon-3}.
For each $n$, let
$\beta_{ii}^{\odot}:=0$ for all $i\in [n]$, and for $1\leq i\neq j\leq n$, define
\begin{align*}
\beta_{ij}^{\odot}
:=
\left\{
\begin{array}{l}
2n^{-\delta_0}\, , \text{ if } i, j\in [n]\setminus B_n\ \text{ and }\ \beta_{ij}\leq n^{-\delta_0},\\[2pt]
\beta_{ij},\text{ otherwise.}
\end{array}
\right.
\end{align*}
Define the function $W_n^{\odot}$ similar to \eqref{eqn:2} using the edge weights $\beta_{ij}^{\odot}$, $i, j\in [n]$.
Then 
\begin{align*}
\big\| 
T_{n^{1/3}(W_n^{\odot} - W_n)}
\big\|_{\sss 2,2}^2
\leq
n^{2/3}\medint{\int}_{[0, 1]^2}\big(W_n^{\odot}-W_n\big)^2
\leq
n^{2/3}\big(2/n^{\delta_0}\big)^2\cdot n^{-2}\cdot o(n^{1+\delta_0})
=
o(1)\, ,
\end{align*}
where the second step uses \eqref{eqn:33}.
Hence, the assumptions in Lemma~\ref{lem:23} are satisfied by $\big(\beta_{ij}^{\odot}\, ;\, i,j\in [n]\big)$, $n\geq 1$, and consequently, the convergences in \eqref{eqn:65} and \eqref{eqn:66} hold for the sequence of random graphs with these edge weights.
Since
\[
\sum_{i, j\in [n]\setminus B_n,\ i\neq j}
\frac{\big( n^{-1}(\beta_{ij} - \beta_{ij}^{\odot} ) \big)^2}{(n^{-1}\beta_{ij}^{\odot} ) }
\leq
\frac{\big(2/n^{1+\delta_0}\big)^2}{2/n^{1+\delta_0}}\cdot o(n^{1+\delta_0})
=
o(1)\, ,
\]
where in the first step we have usesd\eqref{eqn:33}, an application of \cite[Corollary~2.12~(i)]{janson-equiv} shows that the convergences in \eqref{eqn:65} and \eqref{eqn:66} hold for $\graphon_n$.
Finally, the same result can be transferred to $\tildegraphon_n$ by another application of \cite[Corollary~2.12~(i)]{janson-equiv}.
\qed

\subsection{Proofs of Theorem~\ref{thm:graphon-Ui-Uj} and Theorem~\ref{thm:l-2-approx}}\label{sec:g13}
We will verify that Conditions~\ref{ass:graphon-1}, \ref{ass:graphon-2}, and \ref{ass:graphon-3} hold in the settings of both these theorems.

\medskip

\noindent{\bf Proof of Theorem~\ref{thm:graphon-Ui-Uj}:}
Let $U_1, \ldots, U_n$ and $V_1^{\sss (n)}, \ldots, V_n^{\sss (n)}$ be as defined before \eqref{eqn:23}.
Let $W_n$ be as in \eqref{eqn:2} corresponding to the edge weights in \eqref{eqn:23}.
Similarly define
$
F_n(x, y)
=
W\big( V_i^{\sss (n)} ,\, V_j^{\sss (n)} \big) 
$
and 
$
G_n(x, y)
=
H\big( V_i^{\sss (n)} ,\, V_j^{\sss (n)} \big) 
$
if $(x, y)\in Q_{ij}$ for some $1\leq i\neq j\leq n$;
let 
$
F_n(x, y)
=
G_n(x, y)
=
0
$
if $(x, y)\in [0, 1]^2\setminus\big( \bigcup_{1\leq i\neq j\leq n}Q_{ij}\big)$.
Clearly,
$
W_n
=
\big( F_n + n^{-1/3} G_n\big)^+
$.
Further,
$\max_{i\in [n]}
\big| 
V_i^{\sss (n)} - i/n 
\big|
=
o_P\big( n^{-1/2} \cdot \log n\big)
$
whp, and consequently,
\begin{align}\label{eqn:454}
\cE_n
:=
\big\{
\text{for all }1\leq i\neq j\leq n, \
| V_i^{\sss (n)} - x |
+
| V_j^{\sss (n)} - y |
\leq
n^{-1/2} \cdot \log n
\text{ for all }
(x, y)\in Q_{ij}
\big\}
\end{align}
occurs whp.

\medskip

\noindent{\textit{\textbf{Verification of Condition~\ref{ass:graphon-1}}}:}
By the main theorem of \cite{hoeffding-u-statistics}, we have
\begin{align}\label{eqn:332}
\medint{\int}_{[0, 1]^2} F_n^3
\weakc 
\medint\int_{[0, 1]^2} W^3
\ \ \text{ as }\ \ n\to\infty\, .
\end{align}
Next, writing $\eps_n = 1/n$, we note that
\begin{align}\label{eqn:333}
\medint{\int}_{[0, 1]^2\setminus D_n}\big( W - F_n\big)^2
=
\medint{\int}_{[0, 1]^2\setminus(D_n\cup\cA_{\eps_n})}\big( W - F_n\big)^2
+
\medint{\int}_{\cA_{\eps_n}\setminus D_n}\big( W - F_n\big)^2
=:
J_1^{\sss (n)} + J_2^{\sss (n)}\, ,
\end{align}
where $D_n:=\bigcup_{i=1}^n Q_{ii}$ as defined below \eqref{eqn:21}.
Now, 
$
J_1^{\sss (n)}
\leq
\int_{[0, 1]^2\setminus\cA_{\eps_n}}
\osc_{\eps_n; W}^2(x, y)\, dx dy
$
on the event $\cE_n$.
Using Condition~\ref{ass:W-Ui-Uj}~\eqref{item:1}, we get
\begin{align}\label{eqn:334}
n^{2/3}\cdot J_1^{\sss (n)} 
\weakc 
0
\ \ \text{ as }\ \ n\to\infty\, .
\end{align}
To bound $J_2^{\sss(n)}$, note that 
\begin{align*}
&
\medint\int_{\cA_{\eps_n}} 
F_n^2
=
\sum_{1\leq i\neq j\leq n}
W^2\big( V_i^{\sss (n)} ,\, V_j^{\sss (n)} \big) 
\cdot
\mu_{\sss \square} \big(\cA_{\eps_n}\cap Q_{ij} \big)
\notag\\
&\hskip20pt
=
\sum_{1\leq i\neq j\leq n}
W^2\big( V_i^{\sss (n)} ,\, V_j^{\sss (n)} \big) 
\cdot
\mu_{\sss \square} \big(\cA_{\eps_n}\cap Q_{ij} \big)
\cdot
\ind_{\big\{ \big( V_i^{\sss (n)} , V_j^{\sss (n) } \big) \in \widetilde \cA_{\eps_n} \big\} }
\notag\\
&\hskip60pt
\leq
\frac{1}{n^2}
\sum_{1\leq i\neq j\leq n}
W^2\big( V_i^{\sss (n)} ,\, V_j^{\sss (n)} \big) 
\cdot
\ind_{\big\{ \big( V_i^{\sss (n)} , V_j^{\sss (n) } \big) \in \widetilde \cA_{\eps_n} \big\} }
\notag\\
&\hskip80pt
=
\frac{1}{n^2}
\sum_{1\leq i\neq j\leq n}
W^2\big( U_i ,\, U_j \big) 
\cdot
\ind_{\{ ( U_i , U_j ) \in \widetilde \cA_{\eps_n} \} }
\ ,
\end{align*}
where the second equality holds on the event $\cE_n$.
Thus,
$
\bE\big(
\ind_{\cE_n}\cdot \int_{\cA_{\eps_n}}  F_n^2
\big)
\leq
\int_{\widetilde\cA_{\eps_n}} W^2
=
o(n^{-{\sss 2/3}})
$,
where we have again used Condition~\ref{ass:W-Ui-Uj}~\eqref{item:1}.
This last observation combined with the fact that
$
J_2^{\sss (n)}
\leq
2\int_{\cA_{\eps_n}}
(W^2 + F_n^2)
$, 
and \eqref{eqn:333} and \eqref{eqn:334} yields
\begin{align}\label{eqn:335}
n^{2/3}\medint{\int}_{[0, 1]^2\setminus D_n}\big( W - F_n\big)^2
\weakc 
0
\ \ \text{ as }\ \ n\to\infty\, .
\end{align}
If $H=\lambda W$ for some $\lambda\in\bR$, 
combining \eqref{eqn:332} and \eqref{eqn:335} with the fact that \eqref{eqn:21-prime} implies \eqref{eqn:21}, we see that the functions $W_n$, $n\geq 1$, can be constructed on the same probability space where Condition~\ref{ass:graphon-1} is satisfied almost surely.

If $H$ satisfies the condition in Theorem~\ref{thm:graphon-Ui-Uj}~\eqref{item:5}, then it suffices to show that
\begin{gather}\label{eqn:336}
\medint{\int}_{[0, 1]^2\setminus D_n}\big( H - G_n \big)^2
\weakc
0
\qquad \text{ and }\qquad
n^{2/3}\medint{\int}_{[0, 1]^2} \bigg( \big(F_n + n^{-1/3} G_n \big)\wedge 0\bigg)^2 
\weakc 
0\, .
\end{gather}
To prove the first convergence, fix $\eta >0$ and let $S_{\eta}$ be a closed subset of $[0, 1]^2$ such that $\mu_{\sss \square}(S_{\eta}^c)\leq \eta$ and the restriction of $H$ to $S_{\eta}$ is continuous (and hence uniformly continuous).
Then for all large $n$,
\begin{align}
&
\medint\int_{[0, 1]^2\setminus D_n}\pr\big( \big| H(x, y) - G_n(x, y) \big| \geq \eta\big) dx dy
\notag
\\
&\hskip30pt
\leq
\mu_{\sss \square}(S_{\eta}^c) + \pr(\cE_n^c)
+
\sum_{1\leq i\neq j\leq n} \medint{\int}_{Q_{ij}\cap S_{\eta}} 
\pr\big(
(V_i^{\sss (n)}, V_j^{\sss(n)}) \notin S_{\eta}
\big) dx dy
\notag
\\
&\hskip60pt
\leq
\eta + \pr(\cE_n^c) 
+\frac{1}{n^2} \sum_{1\leq i\neq j\leq n} \pr\big(
(U_i, U_j) \notin S_{\eta}
\big)
\leq
2\eta + o(1) \, .
\label{eqn:337}
\end{align}
Now \eqref{eqn:337} combined with the facts 
$\int_{[0, 1]^2} |H^3|<\infty$ and
\begin{align}\label{eqn:338}
\bE\big( \medint\int_{[0, 1]^2} |G_n^3|\big) 
= 
n^{-1}(n-1)\medint\int_{[0, 1]^2} |H^3|
=
O(1)
\end{align}
yields 
$\bE\big[ \int_{[0, 1]^2\setminus D_n}\big( H - G_n \big)^2\big] \to 0$, 
which in particular implies the first convergence in \eqref{eqn:336}.

To prove the second convergence in \eqref{eqn:336}, note that
\begin{align}\label{eqn:339}
&
n^{2/3}\medint{\int}_{[0, 1]^2} \bigg( \big(F_n + n^{-1/3} G_n \big)\wedge 0\bigg)^2 
\leq
\medint\int_{F_n \leq n^{-\delta_0}} G_n^2
+
\medint\int_{|G_n|\geq n^{1/3-\delta_0}} G_n^2
\notag\\[2pt]
&
\leq
\big(\medint\int_{[0, 1]^2} |G_n|^3\big)^{2/3}
\cdot
\mu_{\sss \square}\big(\big\{
(x, y)\, :\, 
F_n(x, y)\leq n^{-\delta_0}
\big\}
\big)^{1/3}
+
n^{\delta_0 - 1/3}
\medint\int_{[0, 1]^2}|G_n|^3
\ .
\end{align}
Let $B_n:=\big\{
i\in [n]\, :\, V_i^{\sss(n)}\in \cB_{2^{  1/\delta_0 }\cdot n^{-1}}
\big\}$,
where $\cB_{\sss\bullet}$ is as in Condition~\ref{ass:W-Ui-Uj}~\eqref{item:3}.
Then
\begin{align}\label{eqn:340}
\bE\big(|B_n|\big)
=
o(n^{\varpi_0})
\ \ \text{ and }\ \
\bE\big(
\#\big\{
(i, j)\, :\, i, j\notin B_n
\ \text{ and }\
W(V_i^{\sss (n)}, V_j^{\sss(n)})
\leq 
2 n^{-\delta_0}
\big\}
\big)
=
o(n^{1+\delta_0})\, ,
\end{align}
where we have used Condition~\ref{ass:W-Ui-Uj}~\eqref{item:3}.
Now the second convergence in \eqref{eqn:336} follows upon combining \eqref{eqn:339}, \eqref{eqn:338}, and \eqref{eqn:340}.

\medskip

\noindent{\textit{\textbf{Verification of Condition~\ref{ass:graphon-2}}}:}
Define
\begin{align*}
f(x):=\medint\int_0^1 W^{3/2}(x, y)dy\, ,\ x\in[0, 1]\, ,
\ \text{ and }\
W_n^{\sss(\wedge)}(x, y) := W(x, y)\wedge n^{2/3}\, ,\ x, y\in [0, 1]\, .
\end{align*}
Let $p_1$ be as in Condition~\ref{ass:W-Ui-Uj}~\eqref{item:2}, and let $\theta_1 = 1/p_1$.
It follows from Condition~\ref{ass:W-Ui-Uj}~\eqref{item:2} that for each $\eps>0$,
\begin{align}\label{eqn:582}
n\cdot\mu_{\circ}(S_{n; \eps})
\to 
0
\ \ \text{ as }\ \ n\to\infty
\ \ \text{ where }\ \ 
S_{n; \eps}:=\{x\in [0, 1]\, :\, f(x) \geq \eps n^{\theta_1}\}\, .
\end{align}
Let $r\in\bZ_{>0}$ be such that $r \leq p_1< 2r$.
Then for each $\eps>0$,
\begin{align}\label{eqn:583}
&
\pr\bigg(
\max_{x\in [0, 1]}\medint{\int_0^1}F_n^{3/2}(x, y)dy
\geq
2 \eps n^{\theta_1}
\bigg)
\notag\\
&\hskip15pt
=
\pr\bigg(
\sum_{ j\in [n]\setminus\{i\} }
\frac{1}{n}\big(W_n^{\sss(\wedge)}(U_i, U_j)\big)^{\frac{3}{2}}
\geq 
2\eps n^{\theta_1}
\text{ for some } i\in [n]
\bigg)
+
o(1)
\notag\\
&\hskip30pt
\leq
n\cdot\pr\bigg(
\sum_{j=2}^n \frac{1}{n}\big(W_n^{\sss(\wedge)}(U_1, U_j)\big)^{3/2}
\geq 2\eps n^{\theta_1}
\bigg)
+ o(1)
\notag\\
&\hskip45pt
\leq
n\cdot\pr(U_1\in S_{n; \eps})
+n\cdot\pr\bigg(
U_1\notin S_{n; \eps}\, ,\ 
\sum_{j=2}^n \frac{1}{n}\big(W_n^{\sss(\wedge)}(U_1, U_j)\big)^{3/2}
\geq 2\eps n^{\theta_1}
\bigg)
+o(1)
\notag\\
&\hskip60pt
\leq
o(1)
+
n\cdot
\medint{\int}_{x\in [0, 1]\setminus S_{n; \eps}}\pr\bigg(
\sum_{j=2}^n \frac{1}{n}\Delta_j(x) \geq \eps n^{\theta_1} 
\bigg)\ dx
\notag\\
&\hskip75pt
\leq
o(1)
+
n^{1-2r\cdot\theta_1}\cdot \eps^{-2r}
\medint{\int}_{x\in [0, 1]\setminus S_{n; \eps}}\bE\bigg[
\bigg(\sum_{j=2}^n \frac{1}{n}\Delta_j(x) \bigg)^{2r}
\bigg]\ dx
\, ,
\end{align}
where the first step uses the fact that $\{W(U_i, U_j)\leq n^{2/3}$ for all $1\leq i\neq j\leq n\}$ occurs whp, and the fourth step uses the notation 
$\Delta_j(x) := 
\big(W_n^{\sss(\wedge)}(x, U_j)\big)^{3/2} 
- 
\bE\big[\big(W_n^{\sss(\wedge)}(x, U_j)\big)^{3/2}\big]$, $j=2,\ldots, n$,
together with \eqref{eqn:582} and the relation
$\bE\big[\big(W_n^{\sss(\wedge)}(x, U_j)\big)^{3/2}\big]
\leq
f(x)
\leq 
\eps n^{\theta_1}$
when $x\in [0, 1]\setminus S_{n; \eps}$.
We claim that
\begin{align}\label{eqn:585}
\bE\bigg[
\bigg(\sum_{j=2}^n \frac{1}{n}\Delta_j(x) \bigg)^{2r}
\bigg]
\leq
C(r)\cdot \big( 1 + f(x)^r\big) 
\end{align}
for a constant $C(r)$ depending only on $r$.
For simplicity, we prove \eqref{eqn:585} when $r=5$.
Writing $\Delta_j$ for $\Delta_j(x)$, we see that
\begin{align}\label{eqn:584}
\bE\bigg[\bigg(\sum_{j=2}^n \Delta_j \bigg)^{10}\bigg]
\leq
C\cdot\bE\bigg[\bigg(\sum_{j=2}^n \Delta_j^2 \bigg)^5\bigg]
&
\leq
C'\cdot\bE\big[
n\Delta_2^{10} + n^2\Delta_2^8\Delta_3^2 + n^2 \Delta_2^6\Delta_3^4
+ n^3\Delta_2^4\Delta_3^4\Delta_4^2
\notag\\
&\hskip40pt
+ n^4\Delta_2^4\Delta_3^2\Delta_4^2\Delta_5^2
+ n^5\Delta_2^2\Delta_3^2\Delta_4^2\Delta_5^2\Delta_6^2
\big] \, ,
\end{align}
where the first step uses the Burkholder-Davis-Gundy inequality \cite{burkholder-davis-gundy}.
Note that
$\bE\big(\Delta_2^2\big)
\leq
\bE\big[\big(W_n^{\sss(\wedge)}(x, U_2)\big)^3\big]
\leq
nf(x)
$, 
and consequently
$
\bE\big[ n^5\Delta_2^2\Delta_3^2\Delta_4^2\Delta_5^2\Delta_6^2 \big]
\leq
n^{10}f(x)^5
$.
Similarly, using the bound $|\Delta_2|\leq 2n$, we get
$
\bE\big[ n^4\Delta_2^4\Delta_3^2\Delta_4^2\Delta_5^2\big]
\leq
4n^6\big(\bE\big[\Delta_2^2\big]\big)^4
\leq 
4n^{10}f(x)^4
$.
Treating the other terms on the right side of \eqref{eqn:584} in a similar manner, we get
$
\bE\big[\big(\sum_{j=2}^n \Delta_j \big)^{10}\big]
\leq 
Cn^{10}\big(1 + f(x)^5\big)
$.
The claim in \eqref{eqn:585} follows similarly for a general $r$.
Now \eqref{eqn:583} combined with \eqref{eqn:585} yields
\[
\max_{x\in [0, 1]}\medint{\int_0^1}F_n^{3/2}(x, y)dy
=
o_P(n^{\theta_1})\, .
\]
This completes the verification of Condition~\ref{ass:graphon-2} when $H=\lambda W$.

If $H$ satisfies the condition in Theorem~\ref{thm:graphon-Ui-Uj}~\eqref{item:5}, then for any sequence $\omega_n\to\infty$,
\begin{align*}
&
\pr\bigg(
\max_{x\in [0, 1]}\medint{\int_0^1}\bigg|\frac{G_n(x, y)}{n^{1/3}}\bigg|^{3/2}dy 
\geq 
\omega_n
\bigg)
=
\pr\bigg(\exists i\in [n]\text{ such that }
\sum_{j\in [n]\setminus\{i\}} |H(U_i, U_j)|^{3/2} 
\geq
n^{3/2}\omega_n
\bigg)
\\
&
\leq
n\cdot\pr\bigg(
\sum_{j =2 }^n |H(U_1, U_j)|^{3/2} 
\geq
n^{3/2}\omega_n
\bigg)
\leq
n\cdot\pr\bigg(
\sum_{j =2 }^n |H(U_1, U_j)|^3
\geq
n^2\omega_n^2
\bigg)
\leq
\frac{1}{\omega_n^2}
\medint\int_{[0, 1]^2}|H|^3
=o(1)\, ,
\end{align*}
which completes the verification of Condition~\ref{ass:graphon-2}.

\medskip

\noindent{\textit{\textbf{Verification of Condition~\ref{ass:graphon-3}}}:}
Note that with $B_n$ as defined below \eqref{eqn:339}, we have
\begin{align*}
\bE\bigg[
\sum_{i\in B_n}\sum_{j\in [n]\setminus\{i\}}
W\big(V_i^{\sss(n)} , V_j^{\sss (n)}\big)^2
\bigg]
=
n(n-1)\medint{\int}_{ \cB_{2^{  1/\delta_0 }\cdot n^{-1}}} \medint{\int}_{ [0, 1] } W(x, y)^2 dx dy
=
o(n^{4/3})\, ,
\end{align*}
where the last step uses the second relation in \eqref{eqn:34}.
This combined with \eqref{eqn:340} completes the verification of Condition~\ref{ass:graphon-3} when $H=\lambda W$.

If $H$ satisfies the condition in Theorem~\ref{thm:graphon-Ui-Uj}~\eqref{item:5}, then
\begin{align}\label{eqn:639}
\bE\bigg[
\sum_{i\in B_n}\sum_{j\in [n]\setminus\{i\}}
\big(n^{-1/3} H(V_i^{\sss(n)} , V_j^{\sss (n)})\big)^2
\bigg]
=
n(n-1) n^{-2/3}\medint{\int}_{ [0, 1]\times\cB_{2^{  1/\delta_0 }\cdot n^{-1}} } H^2 
=
o(n^{4/3})\, .
\end{align}
Further, if $\beta_{ij}$ is as in \eqref{eqn:23}, then 
$\beta_{ij}\leq n^{-\delta_0}$ implies that either 
$W(V_i^{\sss (n)}, V_j^{\sss (n)}) \leq 2n^{-\delta_0}$ or 
$H(V_i^{\sss (n)}, V_j^{\sss (n)}) \leq -n^{1/3 - \delta_0}$.
Thus, in view of \eqref{eqn:340}, the following bound completes the verification of Condition~\ref{ass:graphon-3}:
\begin{align*}
&\bE\big[
\#\big\{
(i, j)\, :\, 1\leq i\neq  j\leq n\, ,\
H(V_i^{\sss (n)}, V_j^{\sss(n)})
\leq 
-n^{1/3 - \delta_0}
\big\}
\big]
\\
&\hskip30pt
=
n(n-1)\cdot\pr\big( H^-(U_1, U_2) \geq n^{1/3 - \delta_0 }\big)
\leq
\frac{n^2}{n^{(1/3 - \delta_0)q} }
\medint{\int}_{H^- \geq n^{1/3 - \delta_0}} (H^-)^q
=
o(n^{1+\delta_0})\, .
\end{align*}

Finally, the claimed result follows in an identical way when $W$ is bounded away from zero and 
$H\in L^3([0, 1]^2)$; we do not repeat the argument.
\qed

\medskip

\noindent{\bf Proof of Theorem~\ref{thm:l-2-approx}:}
Let $W_n$ be as in \eqref{eqn:2} corresponding to the edge weights in \eqref{eqn:56}.
Similarly define
$
F_n(x, y)
=
\dashint_{\sss Q_{ij}} W
$
and 
$
G_n(x, y)
=
\dashint_{\sss Q_{ij}} H
$
if $(x, y)\in Q_{ij}$ for some $1\leq i\neq j\leq n$;
let 
$
F_n(x, y)
=
G_n(x, y)
=
0
$
if $(x, y)\in [0, 1]^2\setminus\big( \bigcup_{1\leq i\neq j\leq n}Q_{ij}\big)$.
Clearly,
$
W_n
=
\big( F_n + n^{-1/3} G_n\big)^+
$.
We make note of the following elementary inequality: 
for any Borel measurable $f:[0, 1]^2\to [0, \infty)$, $p\geq 1$, and $i, j\in [n]$,
\begin{align}\label{eqn:435}
\bigg(\
\dashint_{\sss Q_{ij}} f\bigg)^p
\leq
\dashint_{\sss Q_{ij}} f^p\, .
\end{align}

\medskip

\noindent{\textit{\textbf{Verification of Condition~\ref{ass:graphon-1}}}:}
That \eqref{eqn:20} holds follows immediately from \eqref{eqn:435}.
Further, it follows from \eqref{eqn:21-prime} and \eqref{eqn:41} that \eqref{eqn:21} holds when $H=\lambda W$ for some $\lambda\in\bR$.
If $H$ is as in Theorem~\ref{thm:l-2-approx}~\eqref{item:5a}, to complete the verification of Condition~\ref{ass:graphon-1}, it suffices to prove \eqref{eqn:336} in the present setting.

To prove the first convergence in \eqref{eqn:336}, note that
$
\lim_{n \to \infty}G_n(x, y) = H(x, y)
$
for $\mu_{\sss \square}$-a.e. $(x, y)$ (see, e.g., \cite[Theorem~7.10]{rudin-real-and-complex-analysis}).
Further,
$
(H - G_n)^2
\leq
2 H^2 + 2 (\cM H)^2
$,
where $\cM$ denotes the Hardy-Littlewood maximal function given by 
$\cM H(x, y)
=
\sup_{Q\ni (x, y)} \big(\mu_{\sss \square}(Q)\big)^{-1}\int_Q |H|$,
$x, y\in [0, 1]$,
with the supremum taken over all squares $Q$ containing $(x, y)$.
Now, $\cM H\in L^2([0, 1]^2)$ since $H\in L^2([0, 1]^2)$ 
(see. e.g., \cite[Theorem~2.5]{duoandikoetxea-fourier}), 
and consequently the first convergence in \eqref{eqn:336} holds.

To prove the second convergence in \eqref{eqn:336}, using an argument similar to the one in \eqref{eqn:339}, it suffices to show that 
$
\mu_{\sss \square}\big(\big\{
(x, y)\, :\, 
F_n(x, y)\leq n^{-\delta_0}
\big\}
\big)
=
o(1)
$.
To this end, note that \eqref{eqn:35-l2} implies that we can assume, by choosing a slighly larger $\delta_0$ if necessary, that
\begin{align}\label{eqn:35-l2-prime}
\mu_{\sss \square}\big(
\big\{
(x, y) \in \big([0, 1]\setminus\cB_{1/n}\big)^2\, :\, 
W(x, y)\leq 4n^{-\delta_0}
\big\}
\big)
=o(n^{\delta_0-1}) \ \ \text{ as }\ \ n\to\infty\, .
\end{align}
Now, if
$
\dashint_{\sss Q_{ij}}W \leq 2n^{-\delta_0}
$,
then
$
\mu_{\sss \square}\big(
\big\{
(x, y)\in Q_{ij} \, :\, W(x, y)\leq 4n^{-\delta_0}
\}
\big)
\geq
1/(2n^2)
$.
Using \eqref{eqn:35-l2-prime}, we see that
\begin{align}\label{eqn:436}
\#\bigg\{
(i, j)\, :\, i, j\in [n]\setminus B_n\, ,\ \dashint_{\sss Q_{ij}}W \leq 2n^{-\delta_0}
\bigg\}
=
o(n^{1+\delta_0})\, .
\end{align}
Consequently,
$
\mu_{\sss \square}\big(\big\{
(x, y)\, :\, 
F_n(x, y)\leq 2n^{-\delta_0}
\big\}
\big)
\leq
2|B_n|/n + o(n^{\delta_0 - 1})
=
o(1)
$,
where we have used \eqref{eqn:34-l2} and \eqref{eqn:436}.
Thus, the second convergence in \eqref{eqn:336} follows.

\medskip

\noindent{\textit{\textbf{Verification of Condition~\ref{ass:graphon-2}}}:}
Write $\theta_1=1/p_1$, and note that
\begin{align*}
&
\max_{x\in [0, 1]}\medint{\int_0^1}F_n(x, y)^{3/2}dy
\leq
n\cdot\max_{i\in [n]}
\medint\int_{(i-1)/n}^{i/n}\bigg(
\medint\int_0^1 W(x, y)^{3/2}dy
\bigg)dx
\\
&\hskip40pt
\leq
n\cdot\bigg(
\max_{i\in [n]}
\medint\int_{(i-1)/n}^{i/n}\bigg(
\medint\int_0^1 W(x, y)^{3/2}dy
\bigg)^{p_1}dx
\bigg)^{\theta_1}
\cdot
\big(1/n\big)^{1-\theta_1}
=
o(n^{\theta_1})\, ,
\end{align*}
where the first step uses \eqref{eqn:435}, and the last step uses Condition~\ref{ass:W-l2-approx}~\eqref{item:2-l2}.
This shows that Condition~\ref{ass:graphon-2} holds when $H=\lambda W$.
If $H$ is as in Theorem~\ref{thm:l-2-approx}~\eqref{item:5a}, then proceeding similarly, we get
\begin{align*}
\max_{x\in [0, 1]}\medint{\int_0^1}\bigg|\frac{G_n(x, y)}{n^{1/3}}\bigg|^{3/2}dy 
\leq
\bigg(
\max_{i\in [n]}
\medint\int_{(i-1)/n}^{i/n}\, 
\medint\int_0^1 | H(x, y) |^3
dy dx
\bigg)^{1/2}
=
o(1)\, ,
\end{align*}
which completes the verification of  Condition~\ref{ass:graphon-2}.

\medskip

\noindent{\textit{\textbf{Verification of Condition~\ref{ass:graphon-3}}}:}
If $H=\lambda W$, then combining \eqref{eqn:435} with \eqref{eqn:34-l2} ensures that \eqref{eqn:31} holds,
and further, \eqref{eqn:33} follows from \eqref{eqn:436}.
If $H$ is as in Theorem~\ref{thm:l-2-approx}~\eqref{item:5a}, then analogous to the argument below \eqref{eqn:639}, it would suffice to show that
\begin{align}\label{eqn:678}
\#\big\{
(i, j)\, :\, i, j\in [n],\ \dashint_{\sss Q_{ij}} H \leq - n^{1/3 - \delta_0}
\big\}
= 
o(n^{1+\delta_0})\, .
\end{align}
Note that the left side of \eqref{eqn:678} is bounded by
$
m_n
:=
\#\big\{
(i, j)\, :\, i, j\in [n],\ \dashint_{\sss Q_{ij}} (H^-)^q \geq n^{1 - \delta_0}
\big\}
$.
Further, letting 
$
S_n:= \bigcup\big\{
Q_{ij}
\, :\, 
i, j\in [n]\, ,\ 
\dashint_{\sss Q_{ij}} (H^-)^q \geq n^{1 - \delta_0}
\big\}$, we have
\begin{align}\label{eqn:679}
m_n
\leq
n^{1+\delta_0} \medint{\int}_{S_n}(H^-)^q
\leq
n^{1+\delta_0} \medint{\int}_{[0, 1]^2}(H^-)^q\, .
\end{align}
This in turn implies that $\mu_{\sss \square}(S_n) = m_n/n^2 = o(1)$, which combined with the first bound in \eqref{eqn:679} shows that $m_n = o(n^{1+\delta_0})$.
This completes the proof of \eqref{eqn:678}, and consequently, the verification of Condition~\ref{ass:graphon-3}.

Finally, when $W$ is bounded away from zero and $H\in L^3([0, 1]^2)$, the desired result follows in a similar way, so we do not repeat the argument.
\qed

\subsection{Proof of Theorem~\ref{thm:rgiv-scaling-limit}}\label{sec:g14}
Let $t_{n, \lambda}$ be as in \eqref{eqn:1a}.
Then the number of vertices in $\rgiv_n(\lambda)$ is given by 
\begin{align}\label{eqn:721}
N_n = n t_{n,\lambda} + O_P(\sqrt{n})\, .
\end{align}
Let $\rgiv_{n, \otimes}(\lambda)$ be $\rgiv_n(\lambda)$ where we ignore multiple edges between any two vertices.

Now, conditional on $N_n=m$, and labeling the vertices as $m, m-1, \ldots, 2, 1$ (i.e., the last vertex to appear is labeled $1$),
the arrival times of the vertices have the same joint distribution as 
$
t_{n,\lambda}\big(1 - V_i^{\sss (m) }\big)
$,
$1\leq i\leq m$, where $V_1^{\sss (m)} < \ldots < V_m^{(m)}$ are $m$ i.i.d. Uniform$[0, 1]$ random variables arranged in increasing order.
Further, the edge set of $\rgiv_{n, \otimes}(\lambda)$ can be generated by placing edges independently with respective probabilities
$
1-\exp\big(-t_{n, \lambda} (V_i^{\sss (m) }\wedge V_j^{\sss (m) } ) /n \big)
$
between the vertices $i$ and $j$, $1\leq i<j\leq m$.

Thus, conditional on $N_n$, $\rgiv_{n, \otimes}(\lambda)$ has the same distribution as $\tildegraphon_{N_n}$ constructed using the edge weights
\begin{align}\label{eqn:722}
&
\beta_{ij}^{\sss (N_n)}
:=
\frac{N_n}{n}\cdot t_{n, \lambda} 
\cdot 
\big( V_i^{\sss (N_n) }\wedge V_j^{\sss (N_n) } \big)
=
\big( t_{n,\lambda} + O_P(n^{-1/2})\big) \cdot 
t_{n, \lambda} \cdot 
\big( V_i^{\sss (N_n) }\wedge V_j^{\sss (N_n) } \big)
\\
&
=
\left( 
\frac{\pi^2}{4} + \frac{\lambda\pi}{n^{1/3}} + O_P\big(\frac{1}{\sqrt{n}}\big)
\right)
\cdot 
\big( V_i^{\sss (N_n) }\wedge V_j^{\sss (N_n) } \big)
=
\left( 
\frac{\pi^2}{4} 
+ 
\frac{\lambda\pi}{N_n^{1/3}}\bigg(\frac{\pi}{2}\bigg)^{\frac{1}{3}}
+ 
O_P\big(\frac{1}{\sqrt{n}}\big)
\right)
\cdot 
\big( V_i^{\sss (N_n) }\wedge V_j^{\sss (N_n) } \big)
\, ,\notag
\end{align}
where we have used \eqref{eqn:1a} and \eqref{eqn:721}.
Let
\begin{align}\label{eqn:723}
W(x, y) := \pi^2 (x\wedge y)/4
\ \ \text{ and }\ \
H(x, y) := \lambda\pi^{4/3}(x\wedge y)/2^{1/3}\, , \ \ x, y\in [0, 1]\, .
\end{align}
It can be shown by a direct calculation that $\| T_W\|_{\sss 2,2 }=1$, and the corresponding nonnegative eigenfuction normalized to have $L^2$ norm one is given by
\begin{align}\label{eqn:724}
\psi(x) : = \sqrt{2}\cdot \sin\big(\pi x/2\big)\, , \ \ x\in [0, 1]\, .
\end{align}
By the discussion around \eqref{eqn:45}, $W$ given by \eqref{eqn:723} satisfies Condition~\ref{ass:W-Ui-Uj}.
As seen in the proof of Theorem~\ref{thm:graphon-Ui-Uj}, this implies that the edge weights in \eqref{eqn:23} corresponding to the pair $(W, H)$ in \eqref{eqn:723} satisfy Conditions~\ref{ass:graphon-1}, \ref{ass:graphon-2}, and \ref{ass:graphon-3}.
Consequently, these conditions are also satisfied by the edge weights in \eqref{eqn:722} 
(i.e., the $O_P(n^{-1/2})$ term appearing on the right side of \eqref{eqn:722} does not make any difference).
Hence, an application of Theorem~\ref{thm:scaling-limit-graphon} yields the claimed result in Theorem~\ref{thm:rgiv-scaling-limit} for $\rgiv_{n, \otimes}(\lambda)$.
Finally, it is routine to show that for each $i\geq 1$, whp there are no multiple edges in $\cC_i\big(\rgiv_n(\lambda)\big)$, and consequently, $\cC_i\big(\rgiv_n(\lambda)\big)=\cC_i\big(\rgiv_{n, \otimes}(\lambda)\big)$ whp.
This completes the proof.

\section{Discussion}\label{sec:disc}
An important class of edge weights to which Theorems~\ref{thm:scaling-limit-graphon} and \ref{thm:component-limit-graphon} do not apply directly are edge weights of the form \eqref{eqn:444}, \eqref{eqn:23}, and \eqref{eqn:56} when the kernel $W$ does not satisfy Condition~\ref{ass:graphon-3}, Condition~\ref{ass:W-Ui-Uj}~\eqref{item:3}, and Condition~\ref{ass:W-l2-approx}~\eqref{item:3-l2} respectively, but is otherwise well-behaved (so that it satisfies the other required conditions).
The underlying kernels here are not bounded away from zero and the mass distribution on the set where it assumes small values cannot be handled using the above assumptions.
However, as discussed around Conjecture~\ref{conj:x-max-y}, Theorems~\ref{thm:scaling-limit-graphon} and \ref{thm:component-limit-graphon} may be applicable in such cases if we modify the edge weights by taking maximum with $n^{-\delta}$ for some $\delta\in (1/6, 1/3)$.
So, if one can prove a convergence result with the original edge weights similar to \eqref{eqn:91} that only concerns the component sizes and the surplus edges (which are usually easier to deal with compared to the metric space structure), then that would yield the metric scaling limit as well using a coupling argument.
Thus, a general result of this form (dealing only with the component sizes and the surplus edges) applicable to a broad collection of kernels that are not bounded away from zero would significantly extend the class of edge weights for which the convergence in \eqref{eqn:65} holds.

Consider the setting in Theorem~\ref{thm:graphon-Ui-Uj}.
As can be seen from its proof, Condition~\ref{ass:W-Ui-Uj}~\eqref{item:1} was only needed to ensure that
\eqref{eqn:335} holds.
Now, the left side of \eqref{eqn:335} equals
\begin{align}\label{eqn:550}
n^{2/3}\bigg[
\medint{\int}_{[0, 1]^2\setminus D_n} W^2
+
\frac{1}{n(n-1)}\sum_{1\leq i\neq j\leq n}W^2(U_i, U_j)
-
\frac{4}{n^2}\sum_{1\leq i < j\leq n} 
W\big(V_i^{\sss (n)},\, V_j^{\sss (n)} \big)
\cdot
\dashint_{\sss Q_{ij}} W
\bigg]
\, ,
\end{align}
where $U_1, \ldots, U_n$ and $V_1^{\sss (n)}, \ldots, V_n^{\sss (n)}$ are as defined before \eqref{eqn:23}.
Note that the third term appearing in \eqref{eqn:550} is a weighted $U$-statistic.
Thus, using the expression in \eqref{eqn:550}, it should be possible to prove that \eqref{eqn:335} holds for a class of kernels much bigger than the collection of kernels that satisfy Condition~\ref{ass:W-Ui-Uj}~\eqref{item:1}, and consequently, the assumption in Theorem~\ref{thm:graphon-Ui-Uj} can be relaxed.
However, we are not aware of a result on $U$-statistics that can be directly applied here.

Fix $d\in\bZ_{\geq 3}$, and consider a random graph process on $[n]$ that starts with the empty graph on $[n]$ and evolves as follows: at each step, an available pair of vertices is selected uniformly and an edge is added between them, where a pair of vertices is designated `available' if they are non-adjacent and have degrees strictly smaller than $d$.
This dynamic random graph model with degree restrictions was introduced by Ruci\'nski and Wormald \cite{rucinski-wormald}.
For $t\geq 0$, let $\cG_n^{\text{dr}}(t)$ denote the above random graph process at time $\lfloor nt \rfloor$.
The proof of a phase transition in this process has recently been announced in \cite{warnke-wormald}.
The evolution rule for $\cG_n^{\text{dr}}(\cdot)$ is a variation of that of the dynamic configuration model.
However, the process $\cG_n^{\text{dr}}(\cdot)$ lacks the nice exchangeability properties of the configuration model \cite[Chapter~7]{van2009random}, which makes it harder to analyze.
In particular, it is not clear to us if the near-critical behavior of $\cG_n^{\text{dr}}(\cdot)$ can be studied using the exploration process technique or the spanning tree technique mentioned in Section~\ref{sec:intro}.
However, we believe that Theorem~\ref{thm:gen-2} can be applied to this process to get the metric scaling limit inside the critical window.
To do this, one would require the following result:
Let $|\cC_1^{\text{dr}}(t)|$ denote the size of the largest component of $\cG_n^{\text{dr}}(t)$, $t\geq 0$.
Then there exists $C>0$ and $\delta\in(1/6,\, 1/3)$ such that
\begin{align*}
\lim_{n\to\infty}
\pr\bigg(
|\cC_1^{\text{dr}}(t)|
\leq 
\frac{C\log n}{(t_c - t)^2}
\ \ \text{ and }\ \
\diam\big(\cG_n^{\text{dr}}(t)\big)
\leq 
\frac{C\log n}{(t_c - t)}
\ \ \text{ for all }\ \
t\leq t_c - n^{-\delta}
\bigg)
=
1\, ,
\end{align*}
where $t_c$ denotes the point of phase transition for this process.
It seems to us that the rest of the program can be carried out relatively easily.

We close this section with the result on the critical scaling limit of the stochastic block model mentioned in Section~\ref{sec:stochastic-block-model}.
Fix $k\in\bZ_{>0}$, and let $\kappa_n: [k] \times [k] \to  (0, \infty)$ be symmetric kernels for each $n\geq 1$.
Suppose each vertex $i\in [n]$ has a ``type'' $x_i\in[k]$. 
For $1\leq i\neq j\leq n$, define $\beta_{ij}^{\sss (n)} := \kappa_n(x_i, x_j)$, and let $\cG_n^{\text{sbm}}$ and $\widetilde\cG_n^{\text{sbm}}$ respectively be the random graphs $\graphon_n$ and $\tildegraphon_n$ constructed using these edge weights.
For each $x\in [k]$ and each $n$, let
$
\mu_n(x)
:=
n^{-1}\cdot\#\big\{ i\in [n]\, :\, x_i = x \big\}
$.

\begin{ass}\label{ass:irg-strong}
{\upshape (a)} {\bf Convergence of the kernels:} 
There exists a kernel $\kappa(\cdot, \cdot): [k]\times [k]\to (0, \infty)$ and a matrix $A=(a_{xy})_{x, y\in [k]}$ with real-valued (not necessarily positive) entries such that for each $x,y\in [k]$, 
$ \lim_{n \to \infty} n^{1/3}\left(\kappa_{n}(x,y)-\kappa(x,y)\right)= a_{xy}$.
	
\vskip3pt
	
\noindent{\upshape (b)} {\bf Convergence of the empirical measures:}
There exists a probability measure $\mu$ on $[k]$ and a vector $\vb=(b_1,\hdots, b_k)^t$ such that
$\min_{x\in [k]}\mu(x)>0$ and for each $x\in [k]$, 
\[
n^{1/3}\left(\mu_n(x) - \mu(x)\right)\to b_x\ \ \text{ as }\ \ n \to \infty\, .
\]

\vskip3pt
	
\noindent{\upshape (c)} {\bf Criticality:} 
The operator norm of $T_\kappa$ in $L^2([k], \mu)$ equals one.
\end{ass}

We will simply write $\kappa$ for the $k\times k$ matrix with entries $\kappa(i, j)$. 
Define $\mvmu$ (resp. $\mvone$) to be the vector $(\mu(1),\hdots,\mu(k))^t$ (resp. the $k\times k$ column vector where each coordinate is one).
Let $D:=\mathrm{Diag}(\bmu)$--the 
$k\times k$ diagonal matrix with diagonal entries $\mu(1), \hdots, \mu(k)$. 
Similarly define $B:=\mathrm{Diag}(\vb)$.
Define $m_{ij}:=\mu(j)\kappa(i,j)$ for $i, j\in [k]$ and let $M=\left(m_{ij}\right)_{i,j \in [K]}$.
Note that Condition~\ref{ass:irg-strong}~(c) is equivalent to the Perron root of $M$ being equal to one. 
Let $\vu$ and $\vv$
respectively be the right and left eigenvectors of $M$ corresponding to the eigenvalue one normalized so that $\vv^t\vu=1$
and $\vu^t\vone=1$, \ie
\begin{align}\label{eqn:defn-u-v}
M\vu=\vu\, ,\ \ \vv^t M=\vv^t,\ \ \vu^t\vone=1\, , \ \ \text{ and }\ \ \vv^t\vu=1.
\end{align}
Writing $\vu=(u_1,\hdots, u_K)^t$ and $\vv=(v_1,\hdots, v_K)^t$, define
\begin{align}\label{eqn:defn-alpha-beta-zeta}
\alpha:=\frac{1}{(\vv^t\vone)\cdot(\bmu^t\vu)}\, ,\;\;
\chi:=\frac{\sum_{x\in[K]}v_x u_x^2}{(\vv^t\vone)\cdot(\bmu^t\vu)^2}\, , \ \text{ and }\
\zeta:=\alpha\cdot\left[\vv^t(AD+\kappa B)\vu\right].
\end{align}
\begin{thm}\label{thm:scaling-limit-irg}
Under Condition~\ref{ass:irg-strong}, the convergences in \eqref{eqn:65} and \eqref{eqn:66} continue to hold if we replace $\graphonn_n$ by $\cG_n^{\text{sbm}}$ or $\widetilde\cG_n^{\text{sbm}}$, where the constants are as in \eqref{eqn:defn-alpha-beta-zeta}.
\end{thm}

\section*{Acknowledgments}
We thank two anonymous reviewers for their helpful comments and suggestions on a previous version of the paper.
SB has been partially supported by NSF-DMS grants 
1613072, 1606839, and ARO grant W911NF-17-1-0010.
SS has been supported in part by MATRICS grant MTR/2019/000745 from SERB, 
the DST FIST program - 2021 [TPN--700661], 
and by the  Infosys Foundation, Bangalore.

\bibliographystyle{imsart-number}
\bibliography{scaling}

\end{document}